\documentclass[12pt,reqno]{amsart}

\usepackage{amssymb}
\usepackage{bbm}
\usepackage{mathrsfs}
\usepackage{extarrows}
\usepackage{float}
\usepackage[all]{xy}
\usepackage{pict2e}

\newtheorem{theo}{Theorem}[section]
\newtheorem{lem}[theo]{Lemma}
\newtheorem{fac}[theo]{Fact}
\newtheorem{facs}[theo]{Facts}
\newtheorem{prop}[theo]{Proposition}
\newtheorem{coro}[theo]{Corollary}

\theoremstyle{definition}
\newtheorem{defi}[theo]{Definition}
\newtheorem{alg}[theo]{Algorithm}

\newtheorem*{ackn}{Acknowledgements}

\theoremstyle{remark}
\newtheorem{rem}[theo]{Remark}
\newtheorem{rems}[theo]{Remarks}
\newtheorem{ex}[theo]{Example}
\newtheorem{exs}[theo]{Examples}
\newtheorem{nota}[theo]{Notation}
\newtheorem{ttt}[theo]{}

\newtheorem*{conv}{Convention}
\newtheorem*{cau}{Caution}

\renewcommand{\atop}[2]{\genfrac{}{}{0pt}{}{#1}{#2}}

\newcommand{\br}{ }
\newcommand{\brr}{, }

\newcommand{\Spec}{\mathop{\text{\rm Spec}}\nolimits}
\newcommand{\Aut}{\mathop{\text{\rm Aut}}\nolimits}
\newcommand{\Dic}{\mathop{\text{\rm Dic}}\nolimits}
\newcommand{\Div}{\mathop{\text{\rm Div}}\nolimits}
\newcommand{\Pic}{\mathop{\text{\rm Pic}}\nolimits}
\newcommand{\AGL}{\mathop{\text{\rm AGL}}\nolimits}
\newcommand{\NS}{\mathop{\text{\rm NS}}\nolimits}
\newcommand{\CoInd}{\mathop{\text{\rm CoInd}}\nolimits}
\newcommand{\Hom}{\mathop{\text{\rm Hom}}\nolimits}
\newcommand{\Gal}{\mathop{\text{\rm Gal}}\nolimits}
\newcommand{\Frob}{\mathop{\text{\rm Frob}}\nolimits}
\newcommand{\Br}{\mathop{\text{\rm Br}}\nolimits}

\newcommand{\im}{\mathop{\text{\rm im}}\nolimits}
\newcommand{\res}{\mathop{\text{\rm res}}\nolimits}
\newcommand{\coker}{\mathop{\text{\rm coker}}\nolimits}
\newcommand{\cores}{\mathop{\text{\rm cores}}\nolimits}
\newcommand{\sgn}{\mathop{\text{\rm sgn}}\nolimits}
\newcommand{\pr}{\mathop{\text{\rm pr}}\nolimits}
\newcommand{\cl}{\mathop{\text{\rm cl}}\nolimits}
\renewcommand{\div}{\mathop{\text{\rm div}}\nolimits}
\newcommand{\inv}{\mathop{\text{\rm inv}}\nolimits}
\newcommand{\expo}{\mathop{\text{\rm ex}}\nolimits}

\newcommand{\sep}{\text{\rm sep}}
\newcommand{\cs}{\text{\rm cs}}
\newcommand{\ev}{\text{\rm ev}}
\newcommand{\ab}{\text{\rm ab}}
\newcommand{\st}{\text{\rm st}}
\newcommand{\nr}{\text{\rm nr}}
\newcommand{\et}{\text{\rm \'et}}

\newcommand{\bbF}{{\mathbbm F}}
\newcommand{\bbG}{{\mathbbm G}}

\newcommand{\bbQ}{{\mathbbm Q}}
\newcommand{\bbR}{{\mathbbm R}}
\newcommand{\bbZ}{{\mathbbm Z}}

\newcommand{\bbx}{{\mathbbm x}}

\newcommand{\calC}{{\mathscr{C}}}
\newcommand{\calH}{{\mathscr{H}}}
\newcommand{\calI}{{\mathscr{I}}}
\newcommand{\calK}{{\mathscr{K}}}
\newcommand{\calL}{{\mathscr{L}}}
\newcommand{\calO}{{\mathscr{O}}}
\newcommand{\calU}{{\mathscr{U}}}
\newcommand{\calW}{{\mathscr{W}}}
\newcommand{\calX}{{\mathscr{X}}}

\newcommand{\frakm}{\mathfrak{m}}
\newcommand{\frakp}{\mathfrak{p}}
\newcommand{\frakq}{\mathfrak{q}}
\newcommand{\frakG}{\mathfrak{G}}
\newcommand{\frakP}{\mathfrak{P}}

\newcommand{\Ab}{{\text{\bf A}}}
\newcommand{\Pb}{{\text{\bf P}}}

\renewcommand{\atop}[2]{\genfrac{}{}{0pt}{}{#1}{#2}}

\newcounter{abc}
\newenvironment{abc}{\begin{list}{\rm \alph{abc}) }%
{\usecounter{abc} \leftmargin=0.0pt \labelsep=0.0pt %
\listparindent=0.0pt \labelwidth=0.0pt \parsep=\smallskipamount%
 \itemsep=0.0pt \topsep=0.0pt \partopsep=\smallskipamount}}{\end{list}}

\newcounter{iii}
\newenvironment{iii}{\begin{list}{\rm \roman{iii}) }%
{\usecounter{iii} \leftmargin=0.0pt \labelsep=0.0pt %
\listparindent=0.0pt \labelwidth=0.0pt \parsep=\smallskipamount%
 \itemsep=0.0pt \topsep=0.0pt \partopsep=\smallskipamount}}{\end{list}}

\makeatletter
\def\hsmash{\relax 
  \ifmmode\def\next{\mathpalette\mathhsm@sh}\else\let\next\makehsm@sh
  \fi\next}
\def\makehsm@sh#1{\setbox\z@\hbox{#1}\finhsm@sh}
\def\mathhsm@sh#1#2{\setbox\z@\hbox{$\m@th#1{#2}$}\finhsm@sh}
\def\finhsm@sh{\wd\z@\z@ \box\z@}
\makeatother

\def\rightend#1#2{{%
 \leavevmode\nobreak\hskip .5em plus 1fil
 \penalty600 \hskip 0pt plus -1filll
 \vadjust{}\nobreak\hskip 0pt plus 1filll%
 #1\parfillskip=#2\relax \par}}

\def\eop{\ifmmode\rule[-22pt]{0pt}{1pt}\ifinner\tag*{$\square$}\else\eqno{\square}\fi\else\rightend{$\square$}{0pt}\fi}

\thanks{}

\title[Algebraic Brauer classes on open del Pezzo surfaces]{On the algebraic Brauer classes on open \\degree four del Pezzo surfaces}

\begin{document}

\renewcommand{\thefootnote}{\fnsymbol{footnote}}
\author[J\"org Jahnel]{J\"org Jahnel${}^\ddagger$}

\address{\mbox{Department Mathematik\\ \!Univ.\ \!Siegen\\ \!Walter-Flex-Str.\ \!3\\ \!D-57068 \!Siegen\\ \!Germany}}
\email{jahnel@mathematik.uni-siegen.de}
\urladdr{http://www.uni-math.gwdg.de/jahnel}

\author[Damaris Schindler]{Damaris Schindler${}^{*,\ddagger}$}

\address{Mathematisch Instituut\\ \!Universiteit \!Utrecht\\ \!Budapestlaan~6\\ \!NL-3584 \!CD \!Utrecht\\ The Netherlands}
\email{d.schindler@uu.nl}
\urladdr{http://www.uu.nl/staff/DSchindler}


\date{January~10,~2019.}



\begin{abstract}
We study the algebraic Brauer classes on open del Pezzo surfaces of degree
$4$.
I.e., on the complements of geometrically irreducible hyperplane sections of del Pezzo surfaces of degree
$4$.
We show that the
\mbox{$2$-torsion}
part is generated by classes of two different types. Moreover, there are two types of
\mbox{$4$-torsion}
classes. For each type, we discuss methods for the evaluation of such a class at a rational point over a
\mbox{$p$-adic}
field.
\end{abstract}

\footnotetext[1]{The second author was supported by the {\em NWO\/} Veni Grant No.\ 016.Veni.173.016.}

\footnotetext[3]{All computations are with {\tt magma}~\cite{BCP}.}

\maketitle
\thispagestyle{empty}

\section{Introduction}

In this article, we systematically study the algebraic Brauer classes that may arise on an open del Pezzo surface of degree~four. By~which we mean the
complement~$U$
of a geometrically irreducible hyperplane section of a del Pezzo surface of degree~four.
When~$U$
is defined over a
field~$k$,
the algebraic part of the Brauer group
of~$U$
is contained in the Galois cohomology group
$H^1(\Gal(k^\sep/k), \Pic U_{k^\sep})$.
It~equals that group, at least when
$k$
is a number~field.

Let~$X$
be a del Pezzo surface of degree four and
$U \subset X$
be an open degree four del Pezzo~surface.
Then,~on
$X$,
there are exactly 16 lines, which generate the geometric Picard lattice
$\Pic X_{k^\sep}$,
and therefore
$\Pic U_{k^\sep}$.
The~automorphism group of the geometric Picard lattice is isomorphic to the Weyl group
$W(D_5) \cong (\bbZ/2\bbZ)^4 \!\rtimes\! S_5$
of
order~$1920$.
We~thus have to study the cohomology groups
$H^1(G, \frakP)$,
for
$G \subseteq (\bbZ/2\bbZ)^4 \!\rtimes\! S_5$
any subgroup and
$\frakP$
the
\mbox{$(\bbZ/2\bbZ)^4 \!\rtimes\! S_5$-module}
$\Pic U_{k^\sep}$,
which is the same for all open del Pezzo surfaces of degree~four.

The~group
$(\bbZ/2\bbZ)^4 \!\rtimes\! S_5$
has 197 conjugacy classes of subgroups~\cite{BCP}. Among them, there are five maximal subgroups, which are of indices
$2$,
$5$,
$6$,
$10$,
and~$16$,
respectively. The~first four of these are just the preimages under the natural projection
$p\colon (\bbZ/2\bbZ)^4 \!\rtimes\! S_5 \to S_5$
of the maximal subgroups
$A_5$,
$S_4$,
$\AGL_1(\bbF_{\!5})$,
and~$S_3 \!\times\! S_2$
of~$S_5$.
Note~that there are, in fact, four conjugacy classes of subgroups of
index~$10$
in~$(\bbZ/2\bbZ)^4 \!\rtimes\! S_5$,
but only one of them, that of
$p^{-1}(S_3 \!\times\! S_2)$,
is~maximal.

It~has been known for a while as an experimental result that
$H^1(G, \frakP)$
may only be
$(\bbZ/2\bbZ)^e$,
for
$0 \leq e \leq 4$,
or
$\bbZ/4\bbZ \times (\bbZ/2\bbZ)^e$,
for
$0 \leq e \leq 2$~\cite{JS,BL}. In~Section~\ref{sec_Brauer} of this article, we give a formal proof, not relying on any computer work, for the fact that every class is annihilated
by~$4$.
Moreover,~we classify the nontrivial classes that may arise. It~turns out that
$H^1(G, \frakP)_2$
is generated by classes of two different kinds, which we call Brauer classes of types~I and~II.
A~\mbox{$2$-torsion}
class of type~I occurs when
$G$
is the index~5 subgroup, and every 
\mbox{$2$-torsion}
class of type~I is a restriction of~that. Similarly,~a
\mbox{$2$-torsion}
class of type~II occurs when
$G$
is the maximal index~10 subgroup, and every 
\mbox{$2$-torsion}
class of type~II~is its restriction. Moreover,~there are two types of
\mbox{$4$-torsion}
classes possible, which occur when
$G$
is contained in specific subgroups of orders
$24$
and~$64$,~respectively.

\subsubsection*{The Brauer--Manin obstruction}\leavevmode

\noindent
Let~$k$
be a number~field and let
$\Sigma_k$
denote the set of all its~places. The~Grothendieck--Brauer group is a contravariant functor from the category of~schemes to the category of abelian groups. In~particular, given a scheme
$U$
and a
\mbox{$k_\nu$-rational}
point
$x \colon \Spec k_\nu \to U$,
for
$k_\nu$
the completion
of~$k$
with respect to any
$\nu \in \Sigma_k$,
there~is the restriction~homomorphism
$\smash{x^*\colon \Br(X) \to \Br(k_\nu)}$.
For~a Brauer class
$\alpha \in \Br(U)$,
one~calls
$$\ev_{\alpha,\nu}\colon U(k_\nu) \longrightarrow \bbQ/\bbZ \, , \quad x \mapsto \inv(x^*(\alpha)) \, ,$$
the {\em local evaluation map,} associated
with~$\alpha$.
Here,~$\inv\colon \Br(k_\nu) \hookrightarrow \bbQ/\bbZ$
denotes the invariant map, as~usual. The~local evaluation map
$\ev_{\alpha,\nu}$
is continuous with respect to the
\mbox{$\nu$-adic}
topology
on~$U(k_\nu)$~\cite[Chapter~IV, Proposition~2.3.a.ii)]{Ja}.

Let~$U$
be a separated scheme that is smooth and of finite type
over~$k$
and let
$\calU$
be a {\em model\/}
of~$U$
(cf.~Definition \ref{model}) over the ring
$\calO_k$
of integers
in~$k$.
The~choice of the model
$\calU$
determines a set
$\calU(\calO_k)$
of
{\em \mbox{$\calO_k$-integral\/}}
points, together with a canonical injection
$\calU(\calO_k) \hookrightarrow U(k)$.
Similarly,~$\calU$
determines a set
$\calU(\calO_{k,\nu}) \hookrightarrow U(k_\nu)$
of
{\em \mbox{$\calO_{k,\nu}$-integral\/}}
points, for every non-archimedean
place~$\nu$,
as well as a set
$\calU(\calO_{k,S}) \hookrightarrow U(k)$
of
{\em \mbox{$S$-integral\/}}
points, for
$S \subset \Sigma_k$
any finite set of places that includes all archimedean~ones. Here,~by
$\calO_{k,\nu} \subset k_\nu$,
we denote the ring of all integral elements
in~$k_\nu$,
and by
$\calO_S \subset k$
the ring of all elements that are integral
outside~$S$
(cf.~\cite[\S1]{CX}).

It~is well known that, for every
$\alpha \in \Br(U)$,
there exists a finite set
$\smash{S_{\calU ,\alpha} \subset \Sigma_k}$
\mbox{\cite[\S5.2]{Sk}}
of places such that the restriction
$$\ev_{\alpha,\nu} |_{\calU(\calO_{k,\nu})} \colon \calU(\calO_{k,\nu}) \longrightarrow \bbQ/\bbZ$$
is the zero map for each
$\smash{\nu \in \Sigma_k \!\setminus\! S_{\calU ,\alpha}}$.
Consequently,~for the set of all
\mbox{$S$-integral}
points, one has the inclusions (cf.\ \mbox{\cite[\S1]{CX}})
\begin{align*}
\calU(\calO_{k,S}) \subset \Big( &\prod_{\nu \in S} U(k_\nu) \times \!\!\!\prod_{\nu \in \Sigma_k \setminus S}\!\!\!\! \calU(\calO_{k,\nu}) \Big)^{\Br(U)} \subset \prod_{\nu \in S} U(k_\nu) \times \!\!\!\prod_{\nu \in \Sigma_k \setminus S}\!\!\!\! \calU(\calO_{k,\nu}) \,, \\
\intertext{for}
\Big( &\prod_{\nu \in S} U(k_\nu) \times \!\!\!\prod_{\nu \in \Sigma_k \setminus S}\!\!\!\! \calU(\calO_{k,\nu}) \Big)^{\Br(U)} := \\
\Big\{ (x_\nu)_\nu \in &\prod_{\nu \in S} U(k_\nu) \times \!\!\!\prod_{\nu \in \Sigma_k \setminus S}\!\!\!\! \calU(\calO_{k,\nu}) \;\Big|\; \forall \alpha \in \Br(U)\colon \sum_{\nu \in \Sigma_k}\! \ev_{\alpha,\nu}(x_\nu) = 0 \Big\} \,.
\end{align*}
Here,
$S \subset \Sigma_k$
may be any finite set of places including all archimedean~ones.

This~phenomenon is called the Brauer--Manin obstruction and it is, at least from our point of view, the most important application of the Brauer~group. In~the form described here, for
\mbox{$S$-integral}
points, it is due to J.-L.\ Colliot-Th\'el\`ene and~F.~Xu. In~\cite{CX}, many classical counterexamples to the integral Hasse principle or strong approximation (cf.~\cite[Definition~2.3]{DW} or~\cite[\S7.1]{PR}) off certain primes were explained by the Brauer--Manin~obstruction.
In~another direction, the Brauer-Manin obstruction has been pursued by J.-L.\ Colliot-Th\'el\`ene and O.~Wittenberg~\cite{CW} for families of affine cubic surfaces, such as
the representation problem of an integer by the sum of three~cubes. Moreover,~we advise the reader to consult the work of J.~Berg~\cite{Be} on affine Ch\^atelet~surfaces and that of M.\ Bright and J.\ Lyczak~\cite{BL} concerning certain
$\log K3$
surfaces. Let~us note at this point that open del Pezzo surfaces of degree four, too, form a particular type of
$\log K3$~surfaces.

For~\mbox{$k$-rational}
points in the case of a proper scheme over a number
field~$k$,
the Brauer--Manin obstruction is due to Yu.\,I.~Manin~\cite[Chapter~VI]{Ma}  
and takes the form that a nontrivial Brauer class may exclude certain adelic points from being approximated by
\mbox{$k$-rational}
points. In~particular, the Brauer--Manin~obstruction may explain counterexamples to the Hasse~principle.

\subsubsection*{Methods for explicit evaluation}\leavevmode

\noindent
The~two final sections of this article are devoted to methods for the evaluation of the Brauer classes on open degree four del Pezzo~surfaces. In~doing so, we distinguish between the four~types.

First~of all,
\mbox{$2$-torsion}~classes
of type~I have been evaluated before~\cite{JS}. These~are, in fact, {\em cyclic\/} classes, which means that there exists a normal subgroup
$G' \subset G$
with cyclic quotient
$G/G'$
such that the class vanishes under restriction to
$H^1(G', \frakP)$.
In~this case, Manin's original~\cite[\S45]{Ma} class field theoretic approach to the evaluation~applies.
A~\mbox{$2$-torsion}
class of type~I occurs when one of the degenerate quadrics
in~$\Pb^4$
defining~$U$
is
\mbox{$k$-rational}
and the corresponding two linear systems of conics are interchanged
by~$\Gal(k^\sep/k)$.
For~these concepts related to the geometry
of~$U$,
cf.~Figure~\ref{fig_1} below and the discussion prior to~it.

We~show that
\mbox{$4$-torsion}
classes of type~I are cyclic,~too. They~allow a beautiful geometric interpretation, as~well. There~is a quadrilateral responsible for such a class, which is cut out
of~$U$
by a
\mbox{$k$-rational}
hyperplane. The~edges of this quadrilateral are acted upon cyclically
by~$\Gal(k^\sep/k)$.

However, it happens that
\mbox{$2$-torsion}
classes of type~II are usually non-cyclic, a phenomenon that has no analogue, e.g., for proper cubic surfaces~\cite{EJ10} and shows that the present case is more~difficult. Worse,~$H^1(G, \frakP)_2$
does not even need to be generated by cyclic~classes. Nonetheless,~the explicit evaluation of
\mbox{$2$-torsion}
classes of type~II is not hard. Indeed,~they turn out to be corestrictions of
\mbox{$2$-torsion}
classes of type~I.

\subsubsection*{The generic algorithm}\leavevmode

\noindent
Finally,~the
\mbox{$4$-torsion}
classes of type~II are the most mysterious to~us. For~their evaluation, in general, only a generic algorithm~helps. We~describe such an algorithm (Algorithm~\ref{alg_gen}) as well as our implementation in the final section of this~article. It~seems to us that nothing of this kind has ever been tried before, although similar ideas occur in the work of T.~Preu~\cite{Pr}.

The~basic idea is to describe the Brauer class under consideration
by an explicit
\mbox{$2$-cocycle}
$c$
with values in the function
field~$l(U)^*$,
for~$l$
the field of definition of the 16~lines. In~order to evaluate at a
point~$\xi$
over a
\mbox{$p$-adic}
field~$k_\nu$,
$c$
has to be restricted, at~first. The~result is a
\mbox{$2$-cocycle}
$\smash{c_{\nu,\xi}}$
describing a Brauer class
of~$k_\nu$.

It~is, however, not entirely obvious how to compute
from~$\smash{c_{\nu,\xi}}$
the invariant as an element
in~$\bbQ/\bbZ$.
The~difficulty occurs when the
\mbox{$2$-cocycle}
$\smash{c_{\nu,\xi}}$
does not come via inflation from a cyclic quotient
of~$\smash{\Gal(\overline{k}_\nu/k_\nu)}$.
We~solve this problem by a computer-algebraic approach, cf.~Algorithm~\ref{eval_loc}. In~the case when the groups involved are sufficiently close to being cyclic, such as dihedral groups, it is known that the same may be achieved essentially without using a computer, cf.\ \mbox{\cite[Example~4.4]{Pr}} or \mbox{\cite[\S4]{Be}}.

\begin{conv}
Throughout~this article, we assume that
$k$
is a field of
characteristic~$\neq\! 2$.
\end{conv}

\begin{ackn}
We wish to thank Andreas-Stephan Elsenhans (W\"urzburg) for several hints on how to handle in~{\tt magma} the relatively large number fields that occurred in the computations related to this~project.\smallskip

\noindent
We are grateful to the anonymous referee, who read an earlier version of this article very carefully, for his/her suggestions. They were enormously helpful to improve the presentation of our~results.
\end{ackn}

\section{The Picard group and its automorphism group}

\subsubsection*{Generalities on del Pezzo surfaces of degree four}\leavevmode

\noindent
A~del Pezzo surface is a non-singular, proper algebraic
surface~$X$
over a field
$k$
with an ample anti-canonical sheaf
$\calK^{-1}$.
Every~non-singular complete intersection of two quadrics in
$\Pb^4$
is del Pezzo, according to the adjunction formula, and clearly of degree four~\cite[\S I.7]{Ha}. The converse is true, as~well. For every del Pezzo surface of degree four, its anticanonical image is the complete intersection of two quadrics
in~$\Pb^4$
\cite[Theorem~8.6.2]{Do}.

Thus,~associated with a degree four del Pezzo
surface~$X$,
there is a pencil
$\smash{(\lambda q_1 + \mu q_2)_{(\lambda:\mu) \in \Pb^1}}$
of quadrics
in~$\Pb^4$,
the base locus of which
is~$X$.
This pencil is uniquely determined up to an automorphism
of~$\Pb^4$.
It~contains exactly five degenerate fibres, each of which is exactly of rank
$4$~\cite[Proposition~3.26.iv)]{Wi}.

Over an algebraically closed field, a degree four del Pezzo surface contains exactly 16~lines, which generate the Picard~group. The~same is still true if the base field is only separably closed~\cite[Theorem~1.6]{Va}. Of~more fundamental importance for us, however, are the ten one-dimensional linear systems of conics~\cite[\S2.3]{VAV}, which are obtained as~follows. Take~the five degenerate quadrics in the pencil, which, again, are defined
over~$k$
as soon as
$k$
is only separably closed~\cite[Proposition~3.26.ii)]{Wi}. As~they are of
rank~$4$,
each of them contains two one-dimensional linear systems of planes. Intersecting~with a second quadric from the pencil, one finds one-dimensional linear systems of conics.
In~particular, the ten linear systems of conics naturally break into five~pairs. Numbering~the degenerate quadrics from
$1$
to~$5$,
we denote the linear systems of conics as~follows.

\begin{figure}[H]
\begin{center}
\begin{picture}(200,60)
\linethickness{1pt}
\setlength\unitlength{2pt}
\put(0,3){\circle{1.5}}
\put(0,26){\circle{1.5}}
\put(25,3){\circle{1.5}}
\put(25,26){\circle{1.5}}
\put(50,3){\circle{1.5}}
\put(50,26){\circle{1.5}}
\put(75,3){\circle{1.5}}
\put(75,26){\circle{1.5}}
\put(100,3){\circle{1.5}}
\put(100,26){\circle{1.5}}
\put(-2,30){$e_1$}\put(0,14.5){\oval(5,27)}
\put(23,30){$e_2$}\put(25,14.5){\oval(5,27)}
\put(48,30){$e_3$}\put(50,14.5){\oval(5,27)}
\put(73,30){$e_4$}\put(75,14.5){\oval(5,27)}
\put(98,30){$e_5$}\put(100,14.5){\oval(5,27)}
\put(-2,-4){$e'_1$}
\put(23,-4){$e'_2$}
\put(48,-4){$e'_3$}
\put(73,-4){$e'_4$}
\put(98,-4){$e'_5$}
\end{picture}
\end{center}\vskip-1mm

\caption{The combinatorial structure of the ten linear systems of conics}
\label{fig_1}
\end{figure}
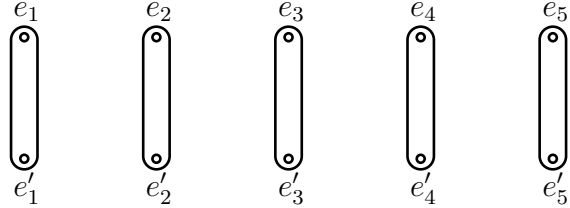\vskip-3\smallskipamount

The~picture might suggest that the automorphism group
of~$\Pic X$
is the wreath product
$S_2 \wr S_5 = (\bbZ/2\bbZ)^5 \!\rtimes\! S_5$.
There~is, however, an additional fine structure that reduces this group by a further index
of~$2$.

\begin{lem}[Automorphisms of the Picard lattice]
Let\/~$X$
be a del Pezzo surface of degree four over a separably closed~field.

\begin{iii}
\item
Then~the
group\/~$W$
of all automorphisms
of\/~$\Pic X$
respecting the canonical
class\/~$K$
and the intersection pairing is isomorphic to the Weyl
group\/~$W(D_5) \cong (\bbZ/2\bbZ)^4 \!\rtimes\! S_5$.
Here,
$(\bbZ/2\bbZ)^4 \subset (\bbZ/2\bbZ)^5$
means the subgroup
\begin{equation}
\label{sum_0}
\{ (a_0, \ldots, a_4) \in (\bbZ/2\bbZ)^5 \mid a_0 + \cdots + a_4 = 0 \}
\end{equation}
being acted upon
by\/~$S_5$
in the natural~way.
\item
The~operation
of\/~$W$
permutes the ten classes of the linear systems of conics faithfully, according to the natural embedding\/
$(\bbZ/2\bbZ)^4 \!\rtimes\! S_5 \hookrightarrow A_{10}$.
\end{iii}\smallskip

\noindent
{\bf Proof.}
{\em
The~isomorphism 
$W(D_5) \cong (\bbZ/2\bbZ)^4 \!\rtimes\! S_5$
is well-known~\cite[\S12]{Hu} and the main part of assertion~i) is \cite[Theorem 23.9]{Ma}. Observe~that
$\Pic X = \Pic X_{\overline{k}}$,
as the latter is generated by the classes of the 16~lines, which occur in
$\Pic X$~already.
Assertion~ii) is explained in \cite[pp.~8--10]{KST}.
}
\eop
\end{lem}

\begin{rem}
$(\bbZ/2\bbZ)^4 \!\rtimes\! S_5$
is a subgroup of the wreath product
$S_2 \wr S_5 = (\bbZ/2\bbZ)^5 \!\rtimes\! S_5$
in~a natural~way, and we work in its standard imprimitive permutation representation
$\iota\colon S_2 \wr S_5 \hookrightarrow S_{10}$
in degree
$10$
\cite[\S2.6]{DM}. Note~that Condition~(\ref{sum_0}) yields~that
\begin{equation}
\label{A10}
(\bbZ/2\bbZ)^4 \!\rtimes\! S_5 = \iota^{-1}(A_{10}) \, . \vspace{2mm}
\end{equation}
\end{rem}\pagebreak[3]

\begin{facs}
\label{PicX}
Let\/~$X$
be a del Pezzo surface of degree four over a separably closed~field.

\begin{iii}
\item
Then,~together with the canonical class\/
$K = -H$,
the classes of the linear systems\/
$[e_1], \ldots, [e_5]$
generate a subgroup of
index\/~$2$
in the Picard~group. More precisely, after a suitable rearrangement of the type
$e_i \leftrightarrow e'_i$,
the group\/
$\Pic X \cong \bbZ^6$
is freely generated by\/
$[e_1], \ldots, [e_5]$,
and\/~$\smash{\frac12(K - [e_1] - \cdots - [e_5])}$.
\item
One~has\/
$[e_i] + [e'_i] = -K$,
for\/~$i = 1, \ldots, 5$.
\item
The classes in\/
$\Pic X$
of the 16 lines are all~the\/
$\smash{\frac12(-3K - [e_1^{(\prime)}] - \cdots - [e_5^{(\prime)}])}$,
in which
$\smash{[e_j^{(\prime)}] \in \{[e_j], [e'_j]\}}$,
for
$j = 1,\ldots,5$,
and the total number of terms
$[e'_j]$
is~even.
\end{iii}\smallskip

\noindent
{\bf Proof.}
{\em
i) is shown in \cite[Proposition~2.2]{VAV}, while ii)~is clear from the construction of the linear systems
$e_i$
and~$e'_i$.\smallskip

\noindent
iii)
Put~$\smash{D := \frac12(-3K - [e_1] - \cdots - [e_5])}$.
Then~direct calculations show that
$\smash{DK = \frac12(-12 + 2\!\cdot\!5) = -1}$
and
$\smash{D^2 = \frac14(36 - 6\!\cdot\!5\!\cdot\!2 + 2(\atop{5}{2})) = -1}$,
such that the Riemann--Roch Theorem~yields
$$\textstyle h^0(X,D) + h^0(X,K-D) \geq 1 + \frac12(D-K)D = 1 \, .$$
However,~$h^0(X,K-D) = 0$,
since
$(K-D)(-K) = -5$
although
$(-K)$
is~ample.
Thus,~$D$
is represented by an effective curve of degree
$D(-K) = 1$,
a~line. The~operation of
$W$
sends
$D$
to an orbit of
length~$16$,
consisting exactly of the divisor classes~described.%
}%
\eop
\end{facs}

\begin{ttt}[The blown-up model]
\label{blowup}
\begin{abc}
\item
A del Pezzo surface of degree four over a separably closed field is isomorphic
to~$\Pb^2$,
blown up in five points
$x_1, \ldots, x_5$
in general position~\cite[Theorem~1.6]{Va}. In~the blown-up model, the 16 lines are given as follows, cf.~\cite[Theorem~26.2]{Ma}.
\begin{iii}
\item
$E_i$,
for~$i=1, \ldots, 5$,
the exceptional curve lying above the blow-up point
$x_i$.
\item
$L_{ij}$,
for~$1 \leq i < j \leq 5$,
the strict transform of the line through
$x_i$
and~$x_j$.
The~class
of~$\smash{L_{ij}}$
in
$\Pic X$
is~$L-E_i-E_j$,
for~$L$
the inverse image of the class of a general line
in~$\Pb^2$.
\item
$C$,
the strict transform of the conic through all five blow-up points. The~class
of~$C$
in
$\Pic X$
is~$2L-E_1-\cdots-E_5$.
\end{iii}\smallskip

\item
In~addition, the ten linear systems of conics are easily identified to lie in the~classes
\begin{iii}
\item
$[e_i] = L - E_i$
and
\item
$[e'_i] = 2L - E_1 - \cdots - E_5 + E_i$,
for~$i=1, \ldots, 5$.
\end{iii}
In~fact, it is well known that lines through a point and conics through four points in general position form pencils of effective curves
in~$\Pb^2$.
Furthermore,~the intersection numbers with the hyperplane section
$H = 3L - E_1 - \cdots - E_5$
are directly seen to be equal
to~$2$,
so that the curves are indeed~conics.
Moreover,~the arrangement of the ten divisor classes given above is such that
$\smash{[e_i] + [e'_i] = H = -K}$,
for~$i=1, \ldots, 5$.\smallskip

\item
One~may now identify by a direct calculation the 16~lines in the above form
$\smash{\frac12(-3K - [e_1^{(\prime)}] - \cdots - [e_5^{(\prime)}])}$
with those in the blown-up~model. It~turns out that

\begin{iii}
\item
(No term
$[e'_\cdot]$)
$$\textstyle \frac12(-3K - [e_1] - \cdots - [e_5]) = [C] \,,$$
\item
(Two terms
$[e'_\cdot]$)
$$\textstyle \frac12\big(-3K - [e_1] - \cdots - [e_5] + ([e_i] - [e'_i]) + ([e_j] - [e'_j])\big) = [L_{ij}] \,,$$
for~$1 \leq i < j \leq 5$,
and
\item
(Four terms
$[e'_\cdot]$)
$$\textstyle \frac12\big(-3K - [e'_1] - \cdots - [e'_5] + ([e'_i] - [e_i])\big) = [E_i] \,,$$
for~$i=1, \ldots, 5$.
\end{iii}
\end{abc}
\end{ttt}

\subsubsection*{Open del Pezzo surfaces of degree four}\leavevmode

\noindent
The open degree four del Pezzo surfaces, which are the subject of this note, are the~following.

\begin{defi}
\label{open_dP4}
By an {\em open del Pezzo surface of degree four,} we mean the complement
$U = X \setminus H$
of a geometrically irreducible hyperplane section
$X \cap H$
of a del Pezzo
surface~$X$
of degree~four.
\end{defi}

\begin{conv}
In~order to have a clear terminology, from this point on, a del Pezzo~surface of degree four in the usual sense is called a {\em proper\/} degree four del~Pezzo surface.
\end{conv}

\begin{lem}
\label{Pic_open}
Let\/~$U = X \setminus H$
be an open del Pezzo surface of degree four over a separably closed~field.

\begin{iii}
\item
Then~the classes of the linear systems\/
$[e_1], \ldots, [e_5]$
freely generate a subgroup of
index\/~$2$
in the Picard
group\/~$\Pic U$.
The~group\/
$\Pic U$
itself is generated by\/
$[e_1], \ldots, [e_5]$,
and\/
$\smash{\frac12([e_1] + \cdots + [e_5])}$.
\item
One~has\/
$[e_i] = -[e'_i]$,
for\/~$i = 1, \ldots, 5$.
\item
After~a suitable rearrangement of the type
$e_i \leftrightarrow e'_i$,
the classes in\/
$\Pic U$
of the 16 lines are all~the\/
$\smash{\frac12(\pm[e_1] \pm \cdots \pm[e_5])}$,
where the total number of plus signs is~even.
\end{iii}\smallskip

\noindent
{\bf Proof.}
{\em
As~$X \cap H$
is geometrically irreducible, \cite[Proposition~II.6.5]{Ha} shows that
$\Pic U = \Pic X/\langle H\rangle = \Pic X/\langle K\rangle$.
Thus,~assertions ii) and~iii) are immediate consequences of their counterparts formulated in Facts~\ref{PicX}. Concerning~i), the direct implication of~Fact~\ref{PicX}.i) is a sixth generator of the form
$\smash{\frac12\big((-1)^{i_1}[e_1] + \cdots + (-1)^{i_5}[e_5]\big)}$,
which, however, yields the same~group.
}
\eop
\end{lem}

\section{The Brauer group}
\label{sec_Brauer}

\subsubsection*{A purely algebraic fact.}\leavevmode

\noindent
The following elementary fact is repeatedly used in this~section.

\begin{lem}
\label{inv_norm}
Let\/~$G$
be a finite group and\/
$M$
a\/
\mbox{$G$-module}
that is finitely generated and torsion-free as a\/
\mbox{$\bbZ$-module}. Then~an element
$m \in M$
is\/
\mbox{$G$-invariant}
if and only if there exists some integer
$c \neq 0$
such that
$cm \in M$
is a norm,
$$cm = \sum_{g\in G} g(m_0)$$
for some
$m_0 \in M$.\medskip

\noindent
{\bf Proof.}
{\em
Norms are clearly
\mbox{$G$-invariant}.
As~$M$
is supposed torsion-free, this implies that the elements described are
\mbox{$G$-invariant},
too.
On~the other hand, the quotient
$\smash{\widehat{H}^0(G,M) = M^G/N_G(M)}$
of invariant elements modulo norms is known to be an abelian group that is annihilated by the order
of~$G$~\cite[Corollary~6.1]{AW}.
The claim immediately follows from~this.
}
\eop
\end{lem}

\subsubsection*{Generalities.}\leavevmode

\noindent
The cohomological Grothendieck--Brauer group
$\Br(U)$
of a scheme
$U$
is, by definition~\cite[Remarque~2.7]{Gr}, the \'etale cohomology group
$H^2_\et(U, \bbG_m)$.
If~$U$
is defined over a field
$k$
then the Hochschild--Serre spectral~sequence~\cite[Exp.\ VIII, Proposition~8.4]{SGA4}
$$H^p(\Gal(k^\sep/k), H^q_\et (U_{k^\sep}, \bbG_m)) \Longrightarrow H^{p+q}_\et (U, \bbG_m)$$
yields a three-step filtration
$$0 \subseteq \Br_0(U) \subseteq \Br_1(U) \subseteq \Br(U) \, .$$

In~the case that
$U$~is
an open degree four del Pezzo surface (Definition~\ref{open_dP4}), one has that
$\smash{\Gamma(U_{k^\sep}, \bbG_m) = (k^\sep)^*}$,
cf.~\cite[Lemma~4.1]{JS}. Hence,
$\Br_0(U)$
is a quotient group of
$\smash{H^2(\Gal(k^\sep/k), (k^\sep)^*) = \Br(k)}$.
It~is true that
$\Br_0(U) \cong \Br(k)$
if
$U$
has a
\mbox{$k$-rational}
point, or, in the case that
$k$
is a number field, an adelic point~\cite[Proposition~1.3.4.1]{Co}. (Note~that the assumption of projectivity made in \cite{Co} is not used in the proof of this particular~statement.)

The subquotient
$\Br_1(U)/\Br_0(U) = \Br_1(U)/\im(\Br(k))$
is called the {\em algebraic\/} part of the Brauer group. It~is, in general, isomorphic to
$$\ker d_2^{1,1}\colon H^1(\Gal(k^\sep/k), \Pic U_{k^\sep}) \longrightarrow H^3(\Gal(k^\sep/k), \Gamma_\et (U_{k^\sep}, \bbG_m)) \, ,$$
for
$d_2^{1,1}$
the differential in the Hochschild--Serre spectral~sequence. Because~of
$\smash{\Gamma_\et (U_{k^\sep}, \bbG_m) = (k^\sep)^*}$,
the right hand side simplifies to
$\smash{H^3(\Gal(k^\sep/k), (k^\sep)^*)}$.
Moreover,~if
$k$
is a number field then, as 
a by-product of class field theory~\cite[section~11.4]{Ta}, it is known that the latter group~vanishes.

Let~us note here that open del Pezzo surfaces of degree four may well have {\em transcendental\/} Brauer classes~\cite[Examples~6.1 and~7.1]{JS}, i.e.\ such not contained
in~$\Br_1(U)$.

\begin{facs}
\label{Lichtenbaum}
\begin{iii}
\item
The composition
$$r\colon H^1(\Gal(k^\sep/k), \Pic U_{k^\sep}) \longrightarrow \Br_1(U)/\Br_0(U) \xlongrightarrow{\overline\res}\Br(k(U))/\im\Br(k)
$$
of the homomorphism defined by the spectral sequence with that induced by the restriction to the generic point has a more direct description as~follows.

The~homomorphism\/
$r$
factors via the quotient~of
$$H^2(\Gal(k^\sep/k), k(U)^\sep{}^*) \subset H^2(\Gal(k(U)^\sep/k(U)), k(U)^\sep{}^*) = \Br(k(U))$$
modulo\/
$\im\Br(k)$,
which itself is clearly contained in\/
$H^2(\Gal(k^\sep/k), k(U)^\sep{}^*/k^\sep{}^*)$.
And~the resulting homomorphism is the boundary map in cohomology associated with the short exact~sequence
$$0 \longrightarrow k(U)^\sep{}^*/k^\sep{}^* \longrightarrow \Div U_{k^\sep} \longrightarrow \Pic U_{k^\sep} \longrightarrow 0 \,.$$
\item\looseness-1
For~example, when
$U$
is an open del Pezzo surface of degree four (Definition~\ref{open_dP4}), the~homomorphism
$H^1(\Gal(k^\sep/k), \Pic U_{k^\sep}) \to \Br_1(U)/\Br_0(U)$
is characterised by this description
of\/~$r$.
\end{iii}\smallskip

\noindent
{\bf Proof.}
{\em
i)~is shown in~\cite[\S2]{Li}. Concerning~ii), one just needs to note that the subsequent homomorphism
$\overline\res$
is~injective. This~is an immediate consequence of the injectivity of
$\res\colon \Br_1(U) \to \Br(k(U))$,
which is true under some minor assumptions~\cite[Corollaire~1.8]{Gr} that are clearly satisfied in our~situation.%
}%
\eop
\end{facs}

\subsubsection*{The algebraic part of the Brauer group.}\leavevmode

\noindent
The observations made in the previous section show that, as an abelian group,
$\Pic U_{k^\sep}$
is always the same, independently of the concrete open del Pezzo
surface~$U$
of degree~four. Thus, let us fix once and for all
that~$\frakP$
is the free abelian group of
rank~$5$,
generated by the symbols
$[e_1], \ldots, [e_5]$,
together with the element
$\smash{\frac12([e_1] + \cdots + [e_5])}$,
and let
$\frakP$
be acted upon
by~$(\bbZ/2\bbZ)^4 \!\rtimes\! S_5$
in the natural way,
$S_5$~permuting
the indices and
$(\bbZ/2\bbZ)^4$
operating by reversing~signs.

\begin{nota}
Let~$G \subseteq (\bbZ/2\bbZ)^4 \!\rtimes\! S_5$
be any~subgroup. Then~we write
$p\colon G \to S_5$
for the natural projection, and put
$T := \ker p$
and~$S := \im p$.
One~thus has
$S \subseteq S_5$,
$T \subseteq (\bbZ/2\bbZ)^4$,
and there is a short exact sequence
$0 \to T \to G \to S \to 0$.
\end{nota}

\begin{nota}[The submodule generated by the linear systems of conics]
We write
$P$
for the subgroup
of~$\frakP$
generated by the symbols
$[e_1], \ldots, [e_5]$.
Then~$P$
is actually a
sub-$G$-module
of~$\frakP$,
and one has the short exact sequence
\begin{equation}
0 \longrightarrow P \longrightarrow \frakP \longrightarrow \bbZ/2\bbZ \longrightarrow 0
\end{equation}
of~$G$-modules,
in~which
$\bbZ/2\bbZ$
represents a module with trivial
\mbox{$G$-operation}.
The~associated long exact sequence in cohomology reads
\begin{equation}
\label{long}
0 \to P^G \to \frakP^G \to \bbZ/2\bbZ \stackrel{\delta}{\to} H^1(G,P) \stackrel{j}{\to} H^1(G,\frakP) \to \Hom(G,\bbZ/2\bbZ) \, .
\end{equation}
\end{nota}

\begin{lem}
\label{ann_2}
Let\/~$G \subseteq (\bbZ/2\bbZ)^4 \!\rtimes\! S_5$
be an arbitrary~subgroup. Then~every element of\/
$H^1(G, P)$
is annihilated
by\/~$2$.\medskip

\noindent
{\bf Proof.}
{\em
{\em First step.}
Inflation.

\noindent
In the inflation-restriction sequence
$$0 \longrightarrow H^1(S,P^T) \longrightarrow H^1(G,P) \longrightarrow H^1(T,P) \,,$$
let us first consider the term
$H^1(S,P^T)$
to the~left.
One~has
$P^T = \bbZ[X]$,
for~$X$
a subset
of~$\{1,\ldots,5\}$
that is invariant
under~$S \subseteq S_5$,
hence
$$P^T \cong \bigoplus_x \bbZ[S/S_x] \cong \bigoplus_x \CoInd_{S_x}^S \bbZ \,, \vspace{-3mm}$$
where the direct sum runs over a system of representatives of the
\mbox{$S$-orbits}
of~$X$,
and~$S_x$
denotes the stabiliser
of~$x$.
Consequently, according to Shapiro's~lemma,
$$H^1(S,P^T) = \bigoplus_x H^1(S, \CoInd_{S_x}^S \bbZ) = \bigoplus_x H^1(S_x,\bbZ) = 0 \,. \vspace{-2.5mm}$$
Observe~here that the operation
of~$S_x$
on~$\bbZ$
is~trivial.
Thus, 
$H^1(G,P)$
injects into
$H^1(T,P)$
and it suffices to show that
$H^1(T,P)$
is annihilated
by~$2$.\smallskip

\noindent
{\em Second step.}
Restriction.

\noindent
As~a
\mbox{$T$-module},
$P$
splits into the direct sum
$P \cong \bigoplus_{i=1}^5 \bbZ\!\cdot\![e_i]$,
so that we only have to show that each
$H^1(T,\bbZ\!\cdot\![e_i])$
is annihilated
by~$2$.
But,~clearly,
$$H^1(T,\bbZ\!\cdot\![e_i]) = H^1(T/T_{[e_i]},\bbZ\!\cdot\![e_i])$$
for
$T_{[e_i]} \subset T$
the stabiliser
of~$[e_i]$.
As~$\#(T/T_{[e_i]}) \leq 2$,
the proof is~complete.
}
\eop
\end{lem}

\begin{theo}
\label{lim4}
Let\/~$G \subseteq (\bbZ/2\bbZ)^4 \!\rtimes\! S_5$
be an arbitrary~subgroup. Then~every element of\/
$H^1(G, \frakP)$
is annihilated
by\/~$4$.\medskip

\noindent
{\bf Proof.}
{\em
In~the exact sequence~(\ref{long}), the term to the far right is obviously annihilated
by~$2$,
while
$H^1(G, P)$
is annihilated
by~$2$,
according to Lemma~\ref{ann_2}. The assertion~follows.
}
\eop
\end{theo}

\begin{rem}
This~result has been known to us before as an experimental finding, just calculating
$H^1(G, \frakP)$
in {\tt magma} for all 197 conjugacy classes of subgroups
of~$(\bbZ/2\bbZ)^4 \!\rtimes\! S_5$
\cite[Remark~4.6.(i)]{JS}. Calculations~of the same kind have recently been carried out more generally for open del Pezzo surfaces of an arbitrary degree
$d \leq 7$
by M.~Bright and J.~Lyczak~\cite{BL}.
\end{rem}

\subsubsection*{$2$-torsion}\leavevmode

\begin{lem}
\label{H1_expl}
Let\/~$G \subseteq (\bbZ/2\bbZ)^4 \!\rtimes\! S_5$
be an arbitrary~subgroup.

\begin{abc}
\item
Then~there are the natural~isomorphisms
\begin{iii}
\item
$\smash{\iota_P\colon H^1(G, P) \stackrel{\cong}{\longleftarrow} (P/2P)^G / (P^G/2P^G)}$
and
\item
$\smash{\iota_{\Pic}\colon H^1(G, \frakP)_2 \stackrel{\cong}{\longleftarrow} (\frakP/2\frakP)^G / (\frakP^G/2\frakP^G)}$.
\end{iii}

Under\/
$\iota_P$
and\/~$\iota_{\Pic}$,
the homomorphism\/
$\smash{j\colon H^1(G, P) \to H^1(G, \frakP)_2}$
(cf.~(\ref{long})) goes over into the homomorphism\/
$(P/2P)^G / (P^G/2P^G) \to (\frakP/2\frakP)^G / (\frakP^G/2\frakP^G)$
induced by~the inclusion\/
$\smash{P \stackrel{\subset}{\longrightarrow} \frakP}$.

\item
For\/~$G' \subseteq G$
another subgroup, the isomorphisms\/
$\iota_P$
and~$\iota_{\Pic}$
are compatible with the restriction\/
$\smash{\res^G_{G'}}$,
while the corestriction\/
$\smash{\cores^G_{G'}}$
goes over into the norm~map.
\item
The boundary homomorphism
$\delta\colon \bbZ/2\bbZ \to H^1(G,P)$
(cf.~(\ref{long})) satisfies
$$\iota_P^{-1}(\delta(\overline{1})) = \overline{[e_1] + \cdots + [e_5]} \,.$$
\end{abc}
\smallskip

\noindent
{\bf Proof.}
{\em
a)
By~Lemma \ref{ann_2},
$H^1(G, P)$
is annihilated
by~$2$.
Moreover,~the commutative~diagram\vskip-6mm
$$
\xymatrix@R=10pt{
0 \ar@{->}[r]& P \ar@{->}[r]^{\cdot 2}\ar@{->}[d]& P \ar@{->}[r]\ar@{->}[d]& P/2P \ar@{->}[r]\ar@{->}[d]& 0 \\
0 \ar@{->}[r]& \frakP \ar@{->}[r]^{\cdot 2}& \frakP \ar@{->}[r]& \frakP/2\frakP \ar@{->}[r]& 0
}
$$
of short exact sequences induces the commutative diagram
$$
\xymatrix@R=10pt{
0 \ar@{->}[r]& P^G \ar@{->}[r]^{\cdot 2}\ar@{->}[d]& P^G \ar@{->}[r]\ar@{->}[d]& (P/2P)^G \ar@{->>}[r]\ar@{->}[d]& H^1(G, P)\ar@{->}[d] \\
0 \ar@{->}[r]& \frakP^G \ar@{->}[r]^{\cdot 2}& \frakP^G \ar@{->}[r]& (\frakP/2\frakP)^G \ar@{->>}[r]& H^1(G, \frakP)_2
}
$$
of long exact sequences in~cohomology.\smallskip

\noindent
b)
$\res^G_{G'}$,
and~$\cores^G_{G'}$
commute with arbitrary boundary homomorphisms, in particular with
$\iota_P$
and~$\iota_{\Pic}$.
The~assertion therefore follows from the naive description of restriction and corestriction on zeroth cohomology~\cite[Chapitre~VII, \S7]{Se}.\smallskip

\noindent
c)
We~observe that the diagram
$$
\xymatrix@R=10pt{
0 \ar@{->}[r]& P \ar@{->}[r]^{\subset\;\;}\ar@{=}[d]& \frakP \ar@{->}[r]\ar@{->}[d]^{\cdot2}& \bbZ/2\bbZ \ar@{->}[r]\ar@{->}[d]^{\overline{1} \mapsto \overline{[e_1] + \cdots + [e_5]}}& 0 \\
0 \ar@{->}[r]& P \ar@{->}[r]^{\cdot2}& P \ar@{->}[r]& P/2P \ar@{->}[r]& 0
}
$$
of short exact sequences is~commutative. Therefore,
$\delta(\overline{1}) = \iota_P(\overline{[e_1] + \cdots + [e_5]})$.
}
\eop
\end{lem}

As a consequence, we observe that all
\mbox{$2$-torsion}
classes in
$H^1(G, \frakP)$
are induced by
$H^1(G, P)$,
and hence can be obtained by only considering the
submodule~$P$
of~$\frakP$
that is generated by
$[e_1], \ldots, [e_5]$.

\begin{coro}
\label{surj}
Let\/~$G \subseteq (\bbZ/2\bbZ)^4 \!\rtimes\! S_5$
be an arbitrary~subgroup. Then~the natural homomorphism\/
$j\colon H^1(G, P) \to H^1(G, \frakP)_2$
is~surjective.\medskip

\noindent
{\bf Proof.}
{\em
According~to Lemma~\ref{H1_expl}.a), we have to show that
$(\frakP/2\frakP)^G$
is generated by the images of
$(P/2P)^G$
and~$\frakP^G$.
For~this, let
$\alpha \in (\frakP/2\frakP)^G$
be an arbitrary element and choose a representative
$\widetilde\alpha \in \frakP$.
The~assumption
on~$\alpha$
then means that
$\sigma(\widetilde\alpha) - \widetilde\alpha \in 2\frakP$,
for
each~$\sigma \in G$.\smallskip

\noindent
{\em First case.}
$\widetilde\alpha \in P$.

\noindent
Without restriction,
$\widetilde\alpha = c_1[e_1] + \cdots + c_5[e_5]$,
for some
$c_1,\ldots,c_5 \in \{0,1\}$.
Then~one clearly has
$\sigma(\widetilde\alpha) - \widetilde\alpha = d_1[e_1] + \cdots + d_5[e_5]$,
for certain
$d_1,\ldots,d_5 \in \{-2,-1,0,1\}$.
Moreover,~as
$\sigma$
operates as a signed permutation,
$d_1+\cdots+d_5$
is~even. However,~such an element can lie in
$2\frakP$
only when every coefficient is even, which implies that
$\widetilde\alpha$
represents a class
in~$(P/2P)^G$.\smallskip

\noindent
{\em Second case.}
$\widetilde\alpha \not\in P$.

\noindent
Without restriction,
$\widetilde\alpha = c_1[e_1] + \cdots + c_5[e_5]$
for
$\smash{c_1,\ldots,c_5 \in \{-\frac12,\frac12\}}$.
Then~one has
$\sigma(\widetilde\alpha) - \widetilde\alpha = d_1[e_1] + \cdots + d_5[e_5]$,
for certain
$d_1,\ldots,d_5 \in \{-1,0,1\}$.
Such~a linear combination lies in
$2\frakP$
when either every coefficient is even or every coefficient is~odd. Therefore, unless
$\smash{\sigma(\widetilde\alpha) - \widetilde\alpha = 0}$
for
each~$\sigma \in G$
and hence
$\widetilde\alpha \in \frakP^G$,
there must be some
$\sigma \in G$
such
that~$\sigma(\widetilde\alpha) = -\widetilde\alpha$.
Then,~in particular, 
$\sigma$
changes the parity of the number of coefficients
of~$\widetilde\alpha$
that are~positive. This~is, however, a contradiction, as elements of
$G \subseteq (\bbZ/2\bbZ)^4 \!\rtimes\! S_5$
operate as permutations followed by evenly many sign~reversals.
}
\eop
\end{coro}

\begin{defi}
Let~$O = \{i_1, \ldots, i_n\} \subseteq \{1, \ldots, 5\}$
be an
\mbox{$S$-orbit}.
Then~the subset
$\{[e_{i_1}], \ldots, [e_{i_n}], -[e_{i_1}], \ldots, -[e_{i_n}]\} \subset \frakP$
is acted upon
by~$G$
either transitively or in such a way that it has exactly two
\mbox{$G$-orbits}.
In~the former case, we say that
$O$~is
is a {\em non-split\/}
\mbox{$S$-orbit}.
In~the latter case, it is called a {\em split\/}
\mbox{$S$-orbit}.
\end{defi}

\begin{theo}[Explicit description of
$H^1(G, \frakP)_2$]
\label{expl_2}
Let\/~$G \subseteq (\bbZ/2\bbZ)^4 \!\rtimes\! S_5$
be an arbitrary subgroup and
let\/~$M$
be the free\/
\mbox{$\bbZ/2\bbZ$-module\/}
generated by the set of all\/
$O_j = \{i_{j,1}, \ldots, i_{j,n_j}\}$,
for\/
$j=1,\ldots,k$,
the orbits
of\/~$S \subseteq S_5$.\smallskip

\noindent
Then\/~$H^1(G, \frakP)_2$
is the quotient
of\/~$M$
modulo\/
$(O_1+\cdots+O_k)$
and the split\/
\mbox{$S$-orbits}.\medskip

\noindent
{\bf Proof.}
{\em
One~clearly has
$\smash{P/2P \cong \bigoplus\limits_{i=1}^5} \,\bbZ/2\bbZ\!\cdot\![e_i]$,~hence
$$(P/2P)^G \cong \bigoplus_{j=1}^k \bbZ/2\bbZ\!\cdot\!([e_{i_{j,1}}] + \cdots + [e_{i_{j,n_j}}]) =: \bigoplus_{j=1}^k \bbZ/2\bbZ\!\cdot\!O_j \, .$$
On~the other hand, by Lemma~\ref{inv_norm},
$P^G$
is generated by all elements obtained as the sum of a
\mbox{$G$-orbit}.
In~the non-split case, such a sum~vanishes. Otherwise,~it is mapped
to~$O_j$
under the projection
to~$(P/2P)^G$.
Thus,~Lemma~\ref{H1_expl}.a.i) shows that
$H^1(G, P)$
is the quotient of
$M$~modulo
all split
\mbox{$S$-orbits}.
Corollary~\ref{surj}
together with Lemma~\ref{H1_expl}.c) proves the~claim.
}
\eop
\end{theo}

In~particular, let us observe the consequence~below.

\begin{coro}
\label{Br_nontriv}
Let\/~$G \subseteq (\bbZ/2\bbZ)^4 \!\rtimes\! S_5$
be an arbitrary subgroup.

\begin{abc}
\item
Then\/
$H^1(G, \frakP) \neq 0$
if and only if\/
$S := p(G)\subseteq S_5$
is intransitive, and at least two of the orbits are non-split.

In this case,
$\iota_{\Pic}(\overline{[e_{i_1}] + \cdots + [e_{i_n}]}) \neq 0$
whenever
$\{i_1, \ldots, i_n\}$
is a non-split~orbit.
\item
In particular,
$H^1(G, \frakP) \neq 0$
implies that
$S$
is~intransitive.
\end{abc}\smallskip

\noindent
{\bf Proof.}
{\em
a)
By Theorem~\ref{lim4}, 
$H^1(G, \frakP) \neq 0$
is possible only when
$H^1(G, \frakP)_2 \neq 0$.
The~explicit description of
$H^1(G, \frakP)_2$
given in Theorem~\ref{expl_2} therefore implies the~claim.\smallskip

\noindent
b)~is a direct consequence of~a).
}
\eop
\end{coro}

Like all results in this section, the corollary above immediately translates into a result about open del Pezzo surfaces of degree four (cf.~Definition~\ref{open_dP4}). For~example, part~b) yields the~following.

\begin{theo}
Let\/~$U \subset X = V(q_1, q_2) \subset \Pb^4$
be an open del Pezzo surface of degree four over a
field\/~$k$
and suppose that\/
$H^1(\Gal(k^\sep/k), \Pic U_{k^\sep}) \neq 0$.
Then

\begin{abc}
\item
the~five degenerate quadrics in the pencil
defining\/~$X$
are acted upon intransitively by the Galois
group\/~$\Gal(k^\sep/k)$.
\item
The~binary quintic form\/
$\det (\lambda q_1 + \mu q_2) \in k[\lambda,\mu]$
is reducible
over\/~$k$.
\end{abc}\smallskip

\noindent
{\bf Proof.}
{\em
One~has
$\Pic U_{k^\sep} \cong \frakP$,
according to Lemma~\ref{Pic_open}.i). Moreover,~the natural operation
of~$\Gal(k^\sep/k)$
yields a homomorphism
$\Gal(k^\sep/k) \to (\bbZ/2\bbZ)^4 \!\rtimes\! S_5$,
whose image we denote
by~$G$,
and whose kernel
by~$\frakG$.
Then~$\frakG$
is a pro-finite group acting trivially on
$\Pic U_{k^\sep}$,
which is torsion-free. Hence,
$H^1(\frakG, \Pic U_{k^\sep}) = 0$.
Thus,~the inflation-restriction exact sequence \cite[Proposition~1.6.6]{NSW}
shows~that
$$H^1(G, \frakP) = H^1(\Gal(k^\sep/k)/\frakG, \Pic U_{k^\sep}) \cong H^1(\Gal(k^\sep/k), \Pic U_{k^\sep}) \,.$$

\noindent
a)
now immediately follows from Corollary~\ref{Br_nontriv} and b)~is just a reformulation of~a).
}
\eop
\end{theo}

\subsubsection*{$2$-torsion classes of types~I and~II}

\begin{defi}
\label{types_2tors}
Let~$\alpha \in H^1(G, \frakP)_2$
be a nonzero~element.

\begin{iii}
\item
Suppose~that
$i \in \{1,\ldots,5\}$
is
\mbox{$S$-invariant}. I.e.,~that
$\{i\}$
is an
\mbox{$S$-orbit}
of length~one. 
If~$\smash{\alpha = \iota_{\Pic}(\overline{[e_i]})}$
then we call
$\alpha$
a 
{\em \mbox{$2$-torsion}
class of type~I}.
\item
Analogously,~if
$\{i,j\} \subset \{1,\ldots,5\}$
is an
\mbox{$S$-orbit}
of length two then we call
$\smash{\alpha = \iota_{\Pic}(\overline{[e_i] + [e_j]})}$
a 
{\em \mbox{$2$-torsion}
class of type~II}.
\end{iii}
\end{defi}

\begin{lem}
Let\/~$G \subseteq (\bbZ/2\bbZ)^4 \!\rtimes\! S_5$
be an arbitrary~subgroup. Then\/
$H^1(G, \frakP)_2$
is generated by
\mbox{$2$-torsion}
classes of types~I and~II.\medskip

\noindent
{\bf Proof.}
{\em
Relying~on the explicit description of
$H^1(G, \frakP)_2$,
given in Theorem~\ref{expl_2}, one sees that an orbit of length~5 defines the zero~class. Furthermore,~orbits of lengths 3 and~4 may be replaced by the complementary
\mbox{$S$-invariant}~sets,
which are of sizes 2 and~1, respectively.
This~completes the~proof.
}
\eop
\end{lem}

\begin{ex}
There~may be an overlap between the two~types. For~example, let
$G$~be
the cyclic group of
order~$4$,
generated by an element fixing
$[e_4]$
and~$[e_5]$
and operating~as
$[e_1] \mapsto [e_2] \mapsto -[e_1]$
and~$[e_3] \mapsto -[e_3]$,~otherwise.

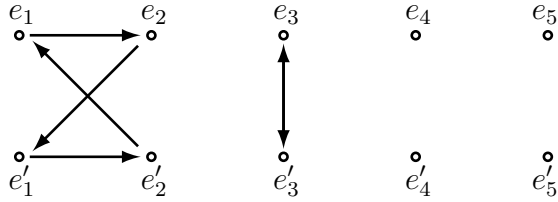
\begin{figure}[H]
\begin{center}	
\begin{picture}(200,60)
\linethickness{1pt}
\setlength\unitlength{2pt}
\put(0,3){\circle{1.5}}
\put(0,26){\circle{1.5}}
\put(25,3){\circle{1.5}}
\put(25,26){\circle{1.5}}
\put(50,3){\circle{1.5}}
\put(50,26){\circle{1.5}}
\put(75,3){\circle{1.5}}
\put(75,26){\circle{1.5}}
\put(100,3){\circle{1.5}}
\put(100,26){\circle{1.5}}
\put(-2,29){$e_1$}\put(2,26){\vector(1,0){21}}
\put(23,29){$e_2$}\put(22.5,24){\vector(-1,-1){20}}
\put(48,29){$e_3$}\put(50,22){\vector(0,-1){17.5}}
\put(73,29){$e_4$}
\put(97,29){$e_5$}
\put(-2,-3){$e'_1$}\put(2,3){\vector(1,0){21}}
\put(23,-3){$e'_2$}\put(22.5,5){\vector(-1,1){20}}
\put(48,-3){$e'_3$}\put(50,7){\vector(0,1){17.5}}
\put(73,-3){$e'_4$}
\put(97,-3){$e'_5$}
\end{picture}
\end{center}\vskip-\medskipamount

\caption{The generator of the cyclic group $G$}\vskip-\bigskipamount
\end{figure}

\noindent
Then~$H^1(G, \frakP)_2 \cong \bbZ/2\bbZ$
and a generator is provided by the~element
$$\iota_{\Pic}(\overline{[e_3]}) = \iota_{\Pic}(\overline{[e_1]+[e_2]}) \,,$$
which may be considered as being of type~I, as well as~II. Let~us note that one has, in fact,
$H^1(G, \frakP) \cong \bbZ/2\bbZ$
in this~example.
\end{ex}

\begin{coro}[A~nontriviality criterion for type~I]
\label{nontriv}
Let\/~$G \subseteq (\bbZ/2\bbZ)^4 \!\rtimes\! S_5$
be an arbitrary subgroup and\/
$\{i\} \subset \{1,\ldots,5\}$
an
\mbox{$S$-orbit}
of length one that is non-split.
I.e.,~such that
$\sigma([e_i]) = -[e_i]$
for some\/~$\sigma \in G$.
Then\/~$\smash{\iota_{\Pic}(\overline{[e_i]}) \neq 0}$.\medskip

\noindent
{\bf Proof.}
{\em
According~to Corollary~\ref{Br_nontriv}, one has to show that there is a second non-split orbit
in~$\{1,\ldots,5\}$.
Assume,~to the contrary, that all the remaining orbits would be~split.
Then~an element
$\sigma \in G$
as above operates
on~$P \subset \frakP$
as~$\sigma([e_i]) = -[e_i]$
and via~cycles
$$[e_j] \mapsto \pm[e_{j_1}] \mapsto \pm[e_{j_2}] \mapsto \cdots \mapsto \pm[e_{j_{l-1}}] \mapsto [e_j]$$
where
the~$j, j_1, \ldots, j_{l-2}$,
and~$j_{l-1}$
are all distinct
and~$\neq\! i$.
Each~such cycle, however, comes in a pair together with a second one that has all signs~reversed.
Hence,~$\sigma$
in total operates as an odd permutation on all the
$\pm[e_1], \ldots, \pm[e_5]$,
which is a contradiction to formula~(\ref{A10}).
}
\eop
\end{coro}

\begin{rem}
Immediately from Definition~\ref{types_2tors}, one has the following two~observations.

\begin{iii}
\item
The~largest subgroup
of~$(\bbZ/2\bbZ)^4 \!\rtimes\! S_5$
that gives rise to a
\mbox{$2$-torsion}
class of type~I is
$p^{-1}(S_4)$,
i.e.\ that of
index~$5$.
This~is clearly a maximal subgroup
of~$(\bbZ/2\bbZ)^4 \!\rtimes\! S_5$,
of
order~$384$.
\item
The~largest subgroup
of~$(\bbZ/2\bbZ)^4 \!\rtimes\! S_5$
that gives rise to a
\mbox{$2$-torsion}
class of type~II is
$p^{-1}(S_2 \!\times\! S_3)$,
of
index~$10$.
This~is a maximal subgroup
of~$(\bbZ/2\bbZ)^4 \!\rtimes\! S_5$,
too, of
order~$192$.
\end{iii}
\end{rem}

\begin{rem}
An experiment in {\tt magma}, running in a loop over the 197 conjugacy classes of subgroups of
$(\bbZ/2\bbZ)^4 \!\rtimes\! S_5$,
shows that
$H^1(G, \frakP)_2$
is isomorphic~to

\begin{iii}
\item
$0$
in 59~cases,
\item
$\bbZ/2\bbZ$
in 71~cases,
\item
$(\bbZ/2\bbZ)^2$
in 47~cases,
\item
$(\bbZ/2\bbZ)^3$
in 17~cases, and
\item
$(\bbZ/2\bbZ)^4$
in three~cases. These~are the group
$G = (\bbZ/2\bbZ)^4$
and its subgroups of orders
$8$
and~$4$
that are still large enough not to split any of the pairs
$\{[e_i], [e'_i]\}$,
for~$i = 1, \ldots, 5$.
The~order~$4$
group may be thought of as being generated by an element reversing the sign
of~$[e_i]$,
for~$i=1,2$,
and that doing the same
for~$i=2,3,4,5$.

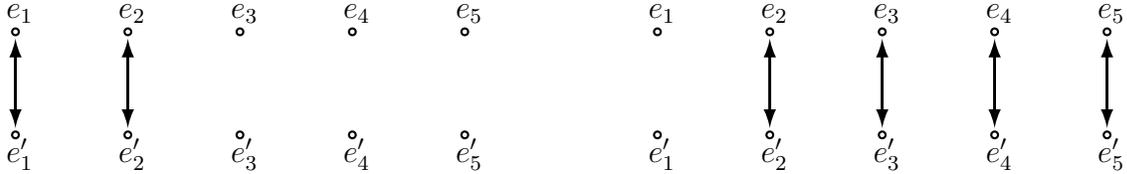
\begin{figure}[H]
\begin{center}	
\begin{picture}(165,55)
\linethickness{.8pt}
\setlength\unitlength{1.7pt}
\put(0,3){\circle{1.5}}
\put(0,26){\circle{1.5}}
\put(25,3){\circle{1.5}}
\put(25,26){\circle{1.5}}
\put(50,3){\circle{1.5}}
\put(50,26){\circle{1.5}}
\put(75,3){\circle{1.5}}
\put(75,26){\circle{1.5}}
\put(100,3){\circle{1.5}}
\put(100,26){\circle{1.5}}
\linethickness{1pt}
\put(-2,29){$e_1$}\put(0,22){\vector(0,-1){17.5}}
\put(23,29){$e_2$}\put(25,22){\vector(0,-1){17.5}}
\put(48,29){$e_3$}
\put(73,29){$e_4$}
\put(98,29){$e_5$}
\put(-2,-3){$e'_1$}\put(0,7){\vector(0,1){17.5}}
\put(23,-3){$e'_2$}\put(25,7){\vector(0,1){17.5}}
\put(48,-3){$e'_3$}
\put(73,-3){$e'_4$}
\put(98,-3){$e'_5$}
\end{picture}\hspace{2.6cm}
\begin{picture}(165,55)
\linethickness{0.8pt}
\setlength\unitlength{1.7pt}
\put(0,3){\circle{1.5}}
\put(0,26){\circle{1.5}}
\put(25,3){\circle{1.5}}
\put(25,26){\circle{1.5}}
\put(50,3){\circle{1.5}}
\put(50,26){\circle{1.5}}
\put(75,3){\circle{1.5}}
\put(75,26){\circle{1.5}}
\put(100,3){\circle{1.5}}
\put(100,26){\circle{1.5}}
\linethickness{1pt}
\put(-2,29){$e_1$}
\put(23,29){$e_2$}\put(25,22){\vector(0,-1){17.5}}
\put(48,29){$e_3$}\put(50,22){\vector(0,-1){17.5}}
\put(73,29){$e_4$}\put(75,22){\vector(0,-1){17.5}}
\put(98,29){$e_5$}\put(100,22){\vector(0,-1){17.5}}
\put(-2,-3){$e'_1$}
\put(23,-3){$e'_2$}\put(25,7){\vector(0,1){17.5}}
\put(48,-3){$e'_3$}\put(50,7){\vector(0,1){17.5}}
\put(73,-3){$e'_4$}\put(75,7){\vector(0,1){17.5}}
\put(98,-3){$e'_5$}\put(100,7){\vector(0,1){17.5}}
\end{picture}
\end{center}

\caption{Generators of the order $4$ group}\vskip-\bigskipamount
\end{figure}
\end{iii}
\end{rem}

\subsubsection*{$4$-torsion}\leavevmode

\begin{lem}
\label{4_tors}
Let\/~$G \subseteq (\bbZ/2\bbZ)^4 \!\rtimes\! S_5$
be an arbitrary~subgroup.

\begin{abc}
\item
Then~there is an~isomorphism\/
$\smash{\iota_{\Pic,4}\colon H^1(G, \frakP) \stackrel{\cong}{\longleftarrow} (\frakP/4\frakP)^G / (\frakP^G/4\frakP^G)}$.
\item
If\/
$c \in \frakP$
is\/
\mbox{$G$-invariant}
modulo\/~$2\frakP$
then\/
$\iota_{\Pic,4}(\overline{2c}) = \iota_{\Pic}(\overline{c})$.
In~particular, if\/
$c$
is\/
\mbox{$G$-invariant}
modulo\/~$4\frakP$
then\/
$2\iota_{\Pic,4}(\overline{c}) = \iota_{\Pic}(\overline{c})$.
\end{abc}\smallskip

\noindent
{\bf Proof.}
{\em
By~Theorem \ref{lim4},
$H^1(G, \frakP)$
is annihilated
by~$4$.
Moreover,~the commutative~diagram
$$
\xymatrix@R=10pt{
0 \ar@{->}[r]& \frakP \ar@{->}[r]^{\cdot2}\ar@{=}[d]& \frakP \ar@{->}[r]\ar@{->}[d]^{\cdot2}& \frakP/2\frakP \ar@{->}[r]\ar@{->}[d]^{\overline{\cdot2}}& 0 \\
0 \ar@{->}[r]& \frakP \ar@{->}[r]^{\cdot4}& \frakP \ar@{->}[r]& \frakP/4\frakP \ar@{->}[r]& 0
}
$$
of short exact sequences induces the commutative diagram
$$
\xymatrix@R=10pt{
0 \ar@{->}[r]& \frakP^G \ar@{->}[r]^{\cdot2}\ar@{=}[d]& \frakP^G \ar@{->}[r]\ar@{->}[d]^{\cdot2}& (\frakP/2\frakP)^G \ar@{->>}[r]\ar@{->}[d]^{\overline{\cdot2}}& H^1(G, \frakP)_2\ar@{^{(}->}[d] \\
0 \ar@{->}[r]& \frakP^G \ar@{->}[r]^{\cdot4}& \frakP^G \ar@{->}[r]& (\frakP/4\frakP)^G \ar@{->>}[r]& H^1(G, \frakP)
}
$$
of long exact sequences in~cohomology. This~immediately establishes part~a) and the first assertion of~b). Observe~that
$2c$
is
\mbox{$G$-invariant}
modulo~$4\frakP$,
since
$c$
is
\mbox{$G$-invariant}
modulo~$2\frakP$.
The~second assertion of~b) is a direct~consequence.
\smallskip
}
\eop
\end{lem}

There are two obvious ways to explicitly construct a
\mbox{$4$-torsion}
class.

\begin{ttt}[$4$-torsion classes of type~I]
\label{type1}
Suppose that
$G \subseteq S_3 \times \bbZ/4\bbZ$,
which occurs as a subgroup of
$(\bbZ/2\bbZ)^4 \!\rtimes\! S_5$
in the manner~below.

\begin{iii}
\item
The generator
$\overline{1} \in \bbZ/4\bbZ$
operates as
$\tau\colon [e_4] \mapsto [e_5] \mapsto -[e_4] \mapsto -[e_5] \mapsto [e_4]$
and
$[e_i] \leftrightarrow -[e_i]$,
for~$i = 1$,
$2$,~$3$.
\item
An~element
$\sigma \in S_3$
sends
$[e_i]$
to~$[e_{\sigma(i)}]$,
for~$i = 1,2,3$,
and fixes
$[e_4]$,
as well
as~$[e_5]$.\vskip-\smallskipamount

\begin{figure}[H]
\begin{center}	
\begin{picture}(165,55)
\linethickness{.8pt}
\setlength\unitlength{1.7pt}
\put(0,3){\circle{1.5}}
\put(0,26){\circle{1.5}}
\put(25,3){\circle{1.5}}
\put(25,26){\circle{1.5}}
\put(50,3){\circle{1.5}}
\put(50,26){\circle{1.5}}
\put(75,3){\circle{1.5}}
\put(75,26){\circle{1.5}}
\put(100,3){\circle{1.5}}
\put(100,26){\circle{1.5}}
\linethickness{1pt}
\put(-2,29){$e_1$}\put(0,22){\vector(0,-1){17.5}}
\put(23,29){$e_2$}\put(25,22){\vector(0,-1){17.5}}
\put(48,29){$e_3$}\put(50,22){\vector(0,-1){17.5}}
\put(73,29){$e_4$}\put(77,26){\vector(1,0){21}}
\put(98,29){$e_5$}\put(97.5,24){\vector(-1,-1){20}}
\put(-2,-3){$e'_1$}\put(0,7){\vector(0,1){17.5}}
\put(23,-3){$e'_2$}\put(25,7){\vector(0,1){17.5}}
\put(48,-3){$e'_3$}\put(50,7){\vector(0,1){17.5}}
\put(73,-3){$e'_4$}\put(77,3){\vector(1,0){21}}
\put(98,-3){$e'_5$}\put(97.5,5){\vector(-1,1){20}}
\end{picture}\hspace{2.6cm}
\begin{picture}(165,55)
\linethickness{0.8pt}
\setlength\unitlength{1.7pt}
\put(0,3){\circle{1.5}}
\put(0,26){\circle{1.5}}
\put(25,3){\circle{1.5}}
\put(25,26){\circle{1.5}}
\put(50,3){\circle{1.5}}
\put(50,26){\circle{1.5}}
\put(75,3){\circle{1.5}}
\put(75,26){\circle{1.5}}
\put(100,3){\circle{1.5}}
\put(100,26){\circle{1.5}}
\linethickness{1pt}
\put(-2,29){$e_1$}\put(2,26){\vector(1,0){21}}
\put(23,29){$e_2$}\put(27,26){\vector(1,0){21}}
\put(48,29){$e_3$}
\put(73,29){$e_4$}
\put(98,29){$e_5$}
\put(-2,-3){$e'_1$}\put(2,3){\vector(1,0){21}}
\put(23,-3){$e'_2$}\put(27,3){\vector(1,0){21}}
\put(48,-3){$e'_3$}\cbezier(3,5)(18,10)(32,10)(47,5)\put(3,5){\vector(-10,-4){1}}
\put(73,-3){$e'_4$}\cbezier(3,24)(18,19)(32,19)(47,24)\put(3,24){\vector(-10,4){1}}
\put(98,-3){$e'_5$}
\end{picture}
\end{center}

\caption{The element $\tau$ and one of the possible elements $\sigma \in S_3$}%
\label{fig_4}%
\end{figure}
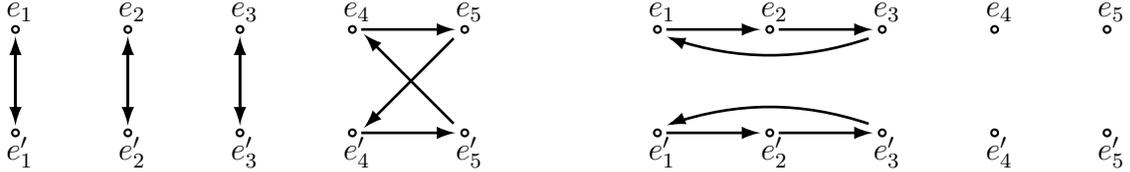%
\end{iii}\vskip-4\smallskipamount

\noindent
Observe~that
$\tau$
clearly commutes with the operation
of~$S_3$
and that every element
of~$S_3 \times \bbZ/4\bbZ$
indeed operates as an even permutation on the ten objects
$\smash{[e_i^{(')}]}$.
\end{ttt}

\begin{prop}
\label{typ1_cl}
Let\/~$G \subseteq (\bbZ/2\bbZ)^4 \!\rtimes\! S_5$
be a subgroup of the type described in~\ref{type1}.

\begin{abc}
\item
Then~the element\/
$[e_1] + [e_2] + [e_3] + 2[e_4] \in \frakP$
is invariant
under\/~$G$
modulo\/~$4\frakP$.
In~particular, one has the~class
$$\alpha := \iota_{\Pic,4}(\overline{[e_1] + [e_2] + [e_3] + 2[e_4]}) \in H^1(G, \frakP) \, .$$
\item
If\/~$G$
surjects
onto\/~$\bbZ/4\bbZ$
under the natural projection then\/
$\alpha$
is a proper\/
\mbox{$4$-torsion}
class.
\end{abc}\smallskip

\noindent
{\bf Proof.}
{\em
a)
The element
$[e_1] + [e_2] + [e_3] + 2[e_4] \in \frakP$
remains fixed
under~$S_3$,
while it is sent
to~$-[e_1] - [e_2] - [e_3] + 2[e_5]$
under~$\tau$.
Hence,
$$\tau([e_1] + [e_2] + [e_3] + 2[e_4]) - ([e_1] + [e_2] + [e_3] + 2[e_4]) = -2[e_1] - 2[e_2] - 2[e_3] - 2[e_4] + 2[e_5] \,,$$
which is an element
of~$4\frakP$.\smallskip

\noindent
b)
Furthermore,~by Lemma~\ref{4_tors}.b), we have
$$2\alpha = \iota_{\Pic}(\overline{[e_1] + [e_2] + [e_3] + 2[e_4]}) = \iota_{\Pic}(\overline{[e_1] + [e_2] + [e_3]}) = \iota_{\Pic}(\overline{[e_4] + [e_5]}) \,.$$
In~order to show that indeed
$2\alpha \neq 0$,
let
$\tau' \in G$
be an element that is mapped under the canonical projection to a generator
of~$\bbZ/4\bbZ$.
Then~$\tau'$
alone makes sure that the
\mbox{$S$-orbit}
$\{4,5\}$
is non-split. In~addition, every
\mbox{$S$-orbit}
contained in
$\{1,2,3\}$
is non-split, as well, so that, in particular, there is at least a second~one. Corollary~\ref{Br_nontriv} shows that
$\smash{\iota_{\Pic}(\overline{[e_4] + [e_5]}) \neq 0}$.
}
\eop
\end{prop}

\begin{defi}
In the situation of Proposition~\ref{typ1_cl}.b), we call
$\alpha \in H^1(G, \frakP)$
a {\em
$4$-torsion
class of type~I}.
\end{defi}

\begin{ttt}[$4$-torsion classes of type~II]
\label{type2}
Assume~that
$G \subseteq (\bbZ/2\bbZ)^4 \!\rtimes\! S_5$
satisfies the following two~conditions.

\begin{iii}
\item
$S = p(G) \subset S_5$
stabilises~$1$.
Moreover,~$S$
is contained in the dihedral
group~$D_4$
respecting the block~\cite[\S1.5]{DM} system
$\{\{2,4\}, \{3,5\}\}$.
I.e.,~for
$\sigma \in G$,
one either has
$p(\sigma)(\{2,4\}) = \{2,4\}$
or
$p(\sigma)(\{2,4\}) = \{3,5\}$.
\item
For~an element
$\sigma \in G$,
its image
$p(\sigma)$
interchanges the two blocks
$\{2,4\}$
and
$\{3,5\}$
if and only if
$\sigma([e_1]) = -[e_1]$.
\end{iii}
\end{ttt}

\begin{exs}
\begin{abc}
\item
Let~$G$
be the cyclic group of
order~$4$
generated by the element depicted~below. Then~conditions i) and~ii) are~satisfied.

\begin{figure}[H]
\begin{center}	
\begin{picture}(200,60)
\linethickness{1pt}
\setlength\unitlength{2pt}
\put(0,3){\circle{1.5}}
\put(0,26){\circle{1.5}}
\put(25,3){\circle{1.5}}
\put(25,26){\circle{1.5}}
\put(50,3){\circle{1.5}}
\put(50,26){\circle{1.5}}
\put(75,3){\circle{1.5}}
\put(75,26){\circle{1.5}}
\put(100,3){\circle{1.5}}
\put(100,26){\circle{1.5}}

\put(-2,29){$e_1$}\put(0,22){\vector(0,-1){17.5}}
\put(23,29){$e_2$}\put(28,26){\vector(1,0){20}}
\put(48,29){$e_3$}\put(46,26){\vector(-1,0){19}}
\put(73,29){$e_4$}\put(77,26){\vector(1,0){21}}
\put(98,29){$e_5$}\put(97.5,24){\vector(-1,-1){20}}
\put(-2,-3){$e'_1$}\put(0,7){\vector(0,1){17.5}}
\put(23,-3){$e'_2$}\put(28,3){\vector(1,0){20}}
\put(48,-3){$e'_3$}\put(46,3){\vector(-1,0){19}}
\put(73,-3){$e'_4$}\put(77,3){\vector(1,0){21}}
\put(98,-3){$e'_5$}\put(97.5,5){\vector(-1,1){20}}
\end{picture}
\end{center}

\caption{The generator of the cyclic group
$G$
of
order~$4$}\vskip-1.2\bigskipamount
\label{fig_5}
\end{figure}
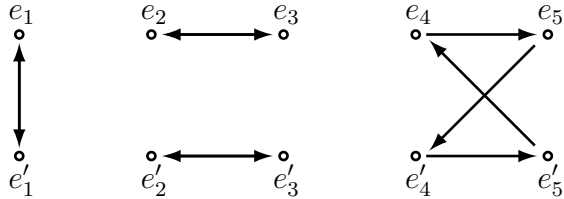
\item
Take
for~$G$
the cyclic group of
order~$8$
generated by the following~element.
Then~$G$
satisfies conditions i) and~ii).

\begin{figure}[H]
\begin{center}	
\begin{picture}(200,60)
\linethickness{1pt}
\setlength\unitlength{2pt}
\put(0,3){\circle{1.5}}
\put(0,26){\circle{1.5}}
\put(25,3){\circle{1.5}}
\put(25,26){\circle{1.5}}
\put(50,3){\circle{1.5}}
\put(50,26){\circle{1.5}}
\put(75,3){\circle{1.5}}
\put(75,26){\circle{1.5}}
\put(100,3){\circle{1.5}}
\put(100,26){\circle{1.5}}

\put(-2,29){$e_1$}\put(0,22){\vector(0,-1){17.5}}
\put(23,29){$e_2$}\put(27,26){\vector(1,0){21}}
\put(48,29){$e_3$}\put(52,26){\vector(1,0){21}}
\put(73,29){$e_4$}\put(77,26){\vector(1,0){21}}
\put(98,29){$e_5$}\put(98,25){\vector(-10,-3){71}}
\put(-2,-3){$e'_1$}\put(0,7){\vector(0,1){17.5}}
\put(23,-3){$e'_2$}\put(27,3){\vector(1,0){21}}
\put(48,-3){$e'_3$}\put(52,3){\vector(1,0){21}}
\put(73,-3){$e'_4$}\put(77,3){\vector(1,0){21}}
\put(98,-3){$e'_5$}\put(98,4){\vector(-10,3){71}}
\end{picture}
\end{center}

\caption{The generator of the cyclic group
$G$
of
order~$8$}\vskip-.5\bigskipamount
\label{fig_6}
\end{figure}
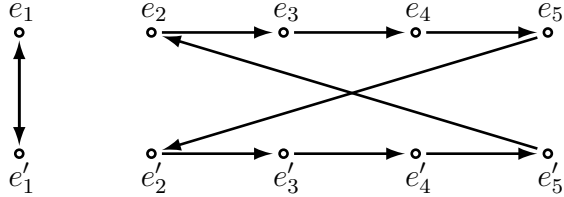
\end{abc}
\end{exs}

\begin{prop}
\label{typ2_cl}
Let\/~$G \subseteq (\bbZ/2\bbZ)^4 \!\rtimes\! S_5$
be a subgroup of the type described in~\ref{type2}.

\begin{abc}
\item
Then~the element\/
$[e_1] + 2[e_2] + 2[e_4] \in \frakP$
is invariant
under\/~$G$
modulo\/~$4\frakP$.
In~particular, one has the~class
$$\alpha := \iota_{\Pic,4}(\overline{[e_1] + 2[e_2] + 2[e_4]}) \in H^1(G, \frakP) \, .$$
\item
If
$\sigma([e_1]) = -[e_1]$
for some\/
$\sigma \in G$
then\/
$\alpha$
is a proper\/
\mbox{$4$-torsion}~class.
\end{abc}\smallskip

\noindent
{\bf Proof.}
{\em
a)
Let~$\sigma \in G$
be an arbitrary~element.
If~$\sigma$
fixes~$[e_1]$
then
$p(\sigma)$
leaves the two blocks in~place. In~this case,
$\sigma([e_1] + 2[e_2] + 2[e_4]) = [e_1] \pm 2[e_2] \pm 2[e_4]$,
so that
$$\sigma([e_1] + 2[e_2] + 2[e_4]) - ([e_1] + 2[e_2] + 2[e_4])$$
is clearly divisible
by~$4$.
Otherwise,
$\sigma([e_1] + 2[e_2] + 2[e_4]) = -[e_1] \pm 2[e_3] \pm 2[e_5]$
and~therefore
$$\sigma([e_1] + 2[e_2] + 2[e_4]) - ([e_1] + 2[e_2] + 2[e_4]) = -2[e_1] - 2[e_2] \pm 2[e_3] - 2[e_4] \pm 2[e_5] \,,$$
which is an element
of~$4\frakP$,~too.\smallskip

\noindent
b)
Lemma~\ref{4_tors}.b) shows that
$2\alpha = \iota_{\Pic}(\overline{[e_1] + 2[e_2] + 2[e_4]}) = \iota_{\Pic}(\overline{[e_1]})$,
which is nonzero according to the criterion given in Corollary~\ref{nontriv}.
}
\eop
\end{prop}

\begin{defi}
In the situation of Proposition~\ref{typ2_cl}.b), we call
$\alpha \in H^1(G, \frakP)$
a {\em
$4$-torsion
class of type~II}.
\end{defi}

\subsubsection*{Occurrence of\/ $4$-torsion}

\begin{conv}
We~say that a proper
\mbox{$4$-torsion}
class
in~$H^1(G, \frakP)$
is {\em of type~I\/} or~{\em II,} if it is so after a suitable permutation of the
indices~$1, \ldots, 5$.
In~this case, we also say that the subgroup\/
$G \subseteq (\bbZ/2\bbZ)^4 \!\rtimes\! S_5$
is {\em of type~I\/} or~{\em II}.
\end{conv}

\begin{rem}[Occurrence of type~I]
\label{rem_typ1}
\begin{abc}
\item
There~are exactly six conjugacy classes of subgroups
of~$(\bbZ/2\bbZ)^4 \!\rtimes\! S_5$
that lead to a
\mbox{$4$-torsion}~class
of type~I. These~correspond one-to-one to the six conjugacy classes of subgroups
of~$S_3 \times \bbZ/4\bbZ$
that surject
onto~$\bbZ/4\bbZ$
under the natural~projection.

\begin{iii}
\item
Among these, there are the groups of the form
$T \times \bbZ/4\bbZ$,
for
$T \subseteq S_3$
any of the four conjugacy classes of~subgroups.
\item
A~fifth conjugacy class is represented by the cyclic subgroup of order four being contained
in~$S_2 \times \bbZ/4\bbZ \subset S_3 \times \bbZ/4\bbZ$
as the kernel~of
$$S_2 \times \bbZ/4\bbZ \xlongrightarrow{(\sgn,\pr)} \bbZ/2\bbZ \times \bbZ/2\bbZ \stackrel{+}{\longrightarrow} \bbZ/2\bbZ \,.$$
Letting~$S_2$
act on
$\{2,3\}$,
a generator of this group operates as~follows.

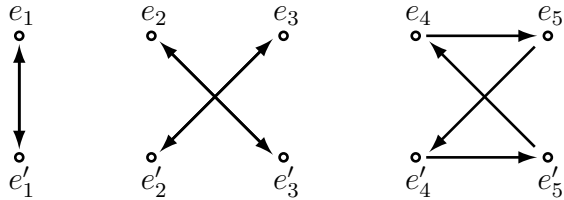
\begin{figure}[H]
\begin{center}	
\begin{picture}(200,60)
\linethickness{1pt}
\setlength\unitlength{2pt}
\put(0,3){\circle{1.5}}
\put(0,26){\circle{1.5}}
\put(25,3){\circle{1.5}}
\put(25,26){\circle{1.5}}
\put(50,3){\circle{1.5}}
\put(50,26){\circle{1.5}}
\put(75,3){\circle{1.5}}
\put(75,26){\circle{1.5}}
\put(100,3){\circle{1.5}}
\put(100,26){\circle{1.5}}

\put(-2,29){$e_1$}\put(0,22){\vector(0,-1){17.5}}
\put(23,29){$e_2$}\put(28,23.5){\vector(1,-1){20}}
\put(48,29){$e_3$}\put(46.5,24){\vector(-1,-1){20}}
\put(73,29){$e_4$}\put(77,26){\vector(1,0){21}}
\put(98,29){$e_5$}\put(97.5,24){\vector(-1,-1){20}}
\put(-2,-3){$e'_1$}\put(0,7){\vector(0,1){17.5}}
\put(23,-3){$e'_2$}\put(28,5.5){\vector(1,1){20}}
\put(48,-3){$e'_3$}\put(46.5,5){\vector(-1,1){20}}
\put(73,-3){$e'_4$}\put(77,3){\vector(1,0){21}}
\put(98,-3){$e'_5$}\put(97.5,5){\vector(-1,1){20}}
\end{picture}
\end{center}\vskip-\smallskipamount

\caption{A generator of a cyclic group of
order~$4$}\vskip-\bigskipamount
\label{fig_5a}
\end{figure}

\item
Similarly,~the sixth conjugacy class is represented by the subgroup of order twelve being the kernel~of
$$S_3 \times \bbZ/4\bbZ \xlongrightarrow{(\sgn,\pr)} \bbZ/2\bbZ \times \bbZ/2\bbZ \stackrel{+}{\longrightarrow} \bbZ/2\bbZ \,.$$
As~an abstract group, this is the dicyclic group
$\Dic_3$.
In~particular, it is non-iso\-mor\-phic to
$A_3 \times \bbZ/4\bbZ$,
which is cyclic of order~twelve.
\end{iii}
\item
The largest subgroup that yields a
\mbox{$4$-torsion}
class of type~I is isomorphic to
$S_3 \times \bbZ/4\bbZ$.
It~is of order 24, thus of index 80
in~$(\bbZ/2\bbZ)^4 \!\rtimes\! S_5$,
and contained in the maximal subgroup of index 10, but not in that of index~5.
\end{abc}
\end{rem}

\begin{rems}[Occurrence of type~II]
\label{rem_typ2}
\begin{abc}
\item
There are exactly eight conjugacy classes of subgroups
of~$(\bbZ/2\bbZ)^4 \!\rtimes\! S_5$
that give rise to a
\mbox{$4$-torsion}~class
of type~II.
\item\looseness-1
The largest of them is of order~64, thus of index 30
in~$(\bbZ/2\bbZ)^4 \!\rtimes\! S_5$.
In~fact, condition \ref{type2}.i) requires an index of at
least~$15$,
which is equivalent to saying that
$G$
be contained in a
\mbox{$2$-Sylow}
subgroup. Furthermore,~condition~ii) enforces another
index~$2$.
It~is contained in the maximal subgroup of index 5, but not in that of index~10.

The groups of type~II are of orders
$4$,
$8$,
$8$,
$16$,
$16$,
$32$,
$32$,
and
$64$,~respectively.
Each~contains an element of order~4 as shown in Figure~\ref{fig_5} or an element of order~8 as shown in Figure~\ref{fig_6}.
\end{abc}
\end{rems}

\begin{theo}
Let\/~$G \subseteq (\bbZ/2\bbZ)^4 \!\rtimes\! S_5$
be a subgroup and assume that\/
$H^1(G, \frakP)$
contains a proper\/
\mbox{$4$-torsion}~element.

\begin{abc}
\item
Then\/~$G$
is of type~I or of type~II.
\item
The~order four group of type~II is of type~I, as~well. There~is no further overlap between these two~types.
\item
One~has\/
$H^1(G, \frakP) \cong \bbZ/4\bbZ \oplus (\bbZ/2\bbZ)^e$,
for\/~$e = 0$,
$1$,
or\/~$2$.
If\/~$G$
is of type~II\/ then\/
$e \leq 1$.
\end{abc}

\noindent
In~particular,
$H^1(G, \frakP)/H^1(G, \frakP)_2$
is of
order\/~$2$,
generated
by\/~$\alpha$
as in Proposition~\ref{typ1_cl} (type~I) or Proposition~\ref{typ2_cl} (type~II).\medskip

\noindent
{\bf Proof.}
{\em
a)
This~is an experimental observation. Running~a loop over all 197 conjugacy classes of subgroups
of~$(\bbZ/2\bbZ)^4 \!\rtimes\! S_5$,
one finds
\mbox{$4$-torsion}
in no case other than those~described.\smallskip

\noindent
b)
The~element shown in Figure~\ref{fig_5} clearly coincides with that depicted in Figure~\ref{fig_5a} after switching
$[e_3]$
with~$[e'_3]$.
This~makes the overlap~visible.

Assume~there would be another overlap. Then~the corresponding group
$G$
must be a
\mbox{$2$-group},
which excludes all type~I cases, except for the groups
$T \times \bbZ/4\bbZ$,
for~$T$
the trivial group
or~$\bbZ/2\bbZ$.
In~these situations,
$S = p(G) \subseteq S_5$
is either a group of
order~$2$,
generated by a
\mbox{$2$-cycle},
or a group of
order~$4$,
generated by two disjoint
\mbox{$2$-cycles}.
Neither~of these occurs as a subgroup
of~$D_4$
that interchanges the two~blocks.\smallskip

\noindent
c)
It~is again an experimental observation that
$H^1(G, \frakP)$
contains only one direct
summand~$\bbZ/4\bbZ$,
for any of the groups~considered. This~includes the overlap case, in which
$H^1(G, \frakP)_2 \cong \bbZ/2\bbZ$
and therefore
$H^1(G, \frakP) \cong \bbZ/4\bbZ$.

In~order to estimate the
exponent~$e$,
we observe that
$H^1(G, \frakP)_2 \cong (\bbZ/2\bbZ)^{e+1}$
and that the number of
\mbox{$S$-orbits},
$\{1,\ldots,5\}$
decomposes into, is at
most~$4$
for type~I and at
most~$3$
for type~II.
Hence,~$e+1 \leq 3$
and
$e+1 \leq 2$,
respectively.
}
\eop
\end{theo}

\begin{rems}
\label{extra_2tors}
\begin{iii}
\item
For~type~I, the groups
$H^1(G, \frakP)$
occurring are
$\bbZ/4\bbZ \oplus (\bbZ/2\bbZ)^2$
in the case of the naive
order~$4$
group (i.e.,
$T$
being the trivial group),
$\bbZ/4\bbZ \oplus \bbZ/2\bbZ$
in the case of the
order~$8$
group (i.e.,
$T \cong \bbZ/2\bbZ$),
and
$\bbZ/4\bbZ$,~otherwise.
\item
For~type~II, the group
$H^1(G, \frakP)$
is
$\bbZ/4\bbZ \oplus \bbZ/2\bbZ$
in two cases and
$\bbZ/4\bbZ$,~otherwise.

\begin{figure}[H]
\begin{center}	
\begin{picture}(165,55)
\linethickness{.8pt}
\setlength\unitlength{1.7pt}
\put(0,3){\circle{1.5}}
\put(0,26){\circle{1.5}}
\put(25,3){\circle{1.5}}
\put(25,26){\circle{1.5}}
\put(50,3){\circle{1.5}}
\put(50,26){\circle{1.5}}
\put(75,3){\circle{1.5}}
\put(75,26){\circle{1.5}}
\put(100,3){\circle{1.5}}
\put(100,26){\circle{1.5}}
\linethickness{1pt}
\put(-2,29){$e_1$}
\put(23,29){$e_2$}\put(25,22){\vector(0,-1){17.5}}
\put(48,29){$e_3$}\put(50,22){\vector(0,-1){17.5}}
\put(73,29){$e_4$}
\put(98,29){$e_5$}
\put(-2,-3){$e'_1$}
\put(23,-3){$e'_2$}\put(25,7){\vector(0,1){17.5}}
\put(48,-3){$e'_3$}\put(50,7){\vector(0,1){17.5}}
\put(73,-3){$e'_4$}
\put(98,-3){$e'_5$}
\end{picture}\hspace{2.6cm}
\begin{picture}(165,55)
\linethickness{0.8pt}
\setlength\unitlength{1.7pt}
\put(0,3){\circle{1.5}}
\put(0,26){\circle{1.5}}
\put(25,3){\circle{1.5}}
\put(25,26){\circle{1.5}}
\put(50,3){\circle{1.5}}
\put(50,26){\circle{1.5}}
\put(75,3){\circle{1.5}}
\put(75,26){\circle{1.5}}
\put(100,3){\circle{1.5}}
\put(100,26){\circle{1.5}}
\linethickness{1pt}
\put(-2,29){$e_1$}
\put(23,29){$e_2$}
\put(48,29){$e_3$}
\put(73,29){$e_4$}\put(75,22){\vector(0,-1){17.5}}
\put(98,29){$e_5$}\put(100,22){\vector(0,-1){17.5}}
\put(-2,-3){$e'_1$}
\put(23,-3){$e'_2$}
\put(48,-3){$e'_3$}
\put(73,-3){$e'_4$}\put(75,7){\vector(0,1){17.5}}
\put(98,-3){$e'_5$}\put(100,7){\vector(0,1){17.5}}
\end{picture}
\end{center}

\caption{Further generators that yield $H^1(G, \frakP) = \bbZ/4\bbZ \oplus \bbZ/2\bbZ$}\vskip-\bigskipamount
\label{fig_7}
\end{figure}
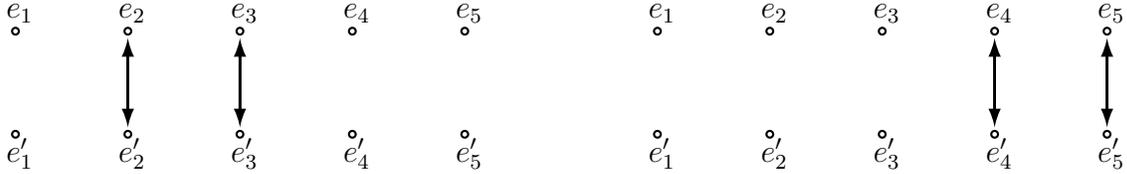

The~two groups that lead to
$H^1(G, \frakP) = \bbZ/4\bbZ \oplus \bbZ/2\bbZ$
are of orders
$8$
and~$16$
and generated by the element of
order~$4$
shown in Figure~\ref{fig_5}, together with the element depicted to the left in Figure~\ref{fig_7} and, in the case of the
order~$16$
group, that shown to the~right.
\end{iii}
\end{rems}

\section{Behaviour under restriction}

\subsubsection*{Proper cubic surfaces--Comparison with a well-known case}\leavevmode

\noindent
For~$\frakP_\cs \cong \bbZ^6$,
the Picard group of a proper cubic surface acted upon by its automorphism group
$W(E_6)$,
there is the following~fact.

\begin{fac}
\label{cub_s}
Let\/~$G' \subseteq G \subseteq W(E_6)$
be subgroups such that neither\/
$H^1(G, \frakP_\cs)$
nor\/
$H^1(G', \frakP_\cs)$
vanishes. Then~the restriction~homomorphism
$$\res_{G'}^G\colon H^1(G, \frakP_\cs) \longrightarrow H^1(G', \frakP_\cs)$$
is~injective.\medskip

\noindent
{\bf Proof.}
{\em
For
\mbox{$2$-torsion},
this is \cite[Corollary~5.9.ii) and~iii)]{EJ10}, while, for
\mbox{$3$-torsion}, 
the same is shown in~\cite[Corollary~3.19.ii) and~iii)]{EJ12}. No~higher torsion occurs in this case~\cite{SwD}.
}
\eop
\end{fac}

\subsubsection*{$2$-torsion}\leavevmode

\noindent
Taking Fact~\ref{cub_s} as a guideline, in the situation of an open degree four del Pezzo surface (Definition~\ref{open_dP4}), one may at least say the~following.

\begin{lem}
Let\/~$G \subseteq (\bbZ/2\bbZ)^4 \!\rtimes\! S_5$
be an arbitrary subgroup,
$G' \subseteq G$
be another, put
$S := p(G)$
and\/
$S' := p(G')$,
and let\/
$\{i_1, \ldots, i_k\} \subseteq \{1, \ldots, 5\}$
be
an\/~\mbox{$S$-orbit}.
Then\/
$\{i_1, \ldots, i_k\}$
is an\/
\mbox{$S'$-invariant}
set, too, and one~has
$$\res_{G'}^G(\iota_{\Pic}^G(\overline{[e_{i_1}] + \cdots + [e_{i_k}]})) = \iota_{\Pic}^{G'}(\overline{[e_{i_1}] + \cdots + [e_{i_k}]}) \,. \smallskip$$

\noindent
{\bf Proof.}
{\em
This is a direct consequence of the compatibility of restriction with the boundary homomorphisms.
}
\eop
\end{lem}

\noindent
In~particular, the two types of
\mbox{$2$-torsion}
classes behave under restriction as~follows.

\begin{iii}
\item
The~restriction of a
\mbox{$2$-torsion}
class of type~I is always a
\mbox{$2$-torsion}
class of type~I
or~$0$.
\item
The~restriction of a
\mbox{$2$-torsion}
class of type~II is a
\mbox{$2$-torsion}
class of type~II or a sum of two
\mbox{$2$-torsion}
classes of type~I. However,~in the latter case, one of the summands or both may degenerate
to~$0$.
\end{iii}

\noindent
These~two observations are certainly a justification for the definition of the two~types.
For~\mbox{$4$-torsion}
classes, there is a stronger~result.

\subsubsection*{$4$-torsion}

\begin{theo}
\label{res_4tors}
Let\/~$G \subseteq (\bbZ/2\bbZ)^4 \!\rtimes\! S_5$
be a subgroup and\/
$G' \subseteq G$
be~another. Suppose~that both\/
$H^1(G, \frakP)$
and\/~$H^1(G', \frakP)$
contain proper\/
\mbox{$4$-torsion}
classes. Then~the homomorphism
$$\overline\res_{G'}^G\colon H^1(G, \frakP)/H^1(G, \frakP)_2 \longrightarrow H^1(G', \frakP)/H^1(G', \frakP)_2 \,,$$
induced by restriction, is an~isomorphism.\medskip

\noindent
{\bf Proof.}
{\em
The groups are of
order~$2$
on either side, so that we only have to show that the homomorphism is nontrivial. For~that, let us recall the commutative diagram of exact~sequences
$$
\xymatrix@R=10pt{
H^1(G, P) \ar@{->}[r]\ar@{->}[d]^{\res}& H^1(G, \frakP) \ar@{->}[r]\ar@{->}[d]^{\res}& H^1(G, \bbZ/2\bbZ) \ar@{->}[d]^{\res} \\
H^1(G', P) \ar@{->}[r] & H^1(G', \frakP) \ar@{->}[r]& H^1(G', \bbZ/2\bbZ) \hsmash{\,,}
}
$$
which, according to Corollary~\ref{surj}, induces the commutative~diagram
$$
\xymatrix@R=10pt{
H^1(G, \frakP)/H^1(G, \frakP)_2 \ar@{^{(}->}[r]^{\;\;\;\;\;\;\;\;\;i_G}\ar@{->}[d]^{\overline\res}& H^1(G, \bbZ/2\bbZ) \ar@{=}[r]\ar@{->}[d]^{\res}& \Hom(G, \bbZ/2\bbZ) \\
H^1(G', \frakP)/H^1(G', \frakP)_2 \ar@{^{(}->}[r]^{\;\;\;\;\;\;\;\;\;\;i_{G'}}& H^1(G', \bbZ/2\bbZ) \ar@{=}[r]& \Hom(G', \bbZ/2\bbZ) \hsmash{\,.}
}
$$
Thus,~the assertion is true, unless
$G'$
is contained in the kernel of the nontrivial homomorphism
$\chi\colon G \to \bbZ/2\bbZ$
in the image
of~$i_G$.
In~order to exclude this possibility, let us explicitly
describe~$\chi$.
We~have to distinguish between the two~types.\smallskip

\noindent
{\em First case.}
$G$
is of type~I.

\noindent
Then, after a suitable permutation of the indices, the situation is as in~\ref{type1} and
$H^1(G, \frakP)/H^1(G, \frakP)_2$
is generated by the class of
$\smash{\alpha = \iota_{\Pic,4}(\overline{[e_1] + [e_2] + [e_3] + 2[e_4]})}$.
The~homomorphism
$\chi\colon G \to \bbZ/2\bbZ$
is then the composition of any
\mbox{$1$-cocycle}
representing
$\alpha \in H^1(G, \frakP)$
with the natural projection
$\frakP \twoheadrightarrow \frakP/P \cong \bbZ/2\bbZ$.
Calculating~the boundary map
$\delta = \iota_{\Pic,4} \colon H^0(G, \frakP/4\frakP) \to H^1(G, \frakP)$
in the usual way, one finds that
$\chi$~is
the restriction
to~$G$
of the homomorphism
$$S_3 \times \bbZ/4\bbZ \to \bbZ/2\bbZ \,, \qquad \sigma\tau^i \mapsto (i \bmod 2), \text{ for any } \sigma \in S_3 \,.$$\vskip.5mm

\noindent
{\em Second case.}
$G$
is of type~II.

\noindent
Here,~after permuting indices if necessary, the situation is as described in~\ref{type2} and
$H^1(G, \frakP)/H^1(G, \frakP)_2$
is generated by the class of
$\smash{\alpha = \iota_{\Pic,4}(\overline{[e_1] + 2[e_2] + 2[e_4]})}$.
Calculating~the boundary map
$\delta = \iota_{\Pic,4} \colon H^0(G, \frakP/4\frakP) \to H^1(G, \frakP)$
and reducing
modulo~$P$,
one finds that
$\chi\colon G \to \bbZ/2\bbZ$
is given~by
$$G \to \bbZ/2\bbZ, \qquad
\left\{
\begin{array}{ll}
\sigma \to \overline{0} & ~\text{ if } \sigma([e_1]) = [e_1] \,, \\
\sigma \to \overline{1} & ~\text{ if } \sigma([e_1]) = -[e_1] \,.
\end{array}
\right.$$\vskip.5mm

\noindent
The assumption that
$G' \subseteq \ker\chi$
therefore yields that one of the
$[e_i]$
is completely fixed
under~$G'$
or that
$\{1,\ldots,5\}$
contains a split orbit of size~three. Either~conclusion is incompatible with the assumption that
$H^1(G', \frakP)$
contains a proper
\mbox{$4$-torsion}
element of either~type.
}
\eop
\end{theo}

Thus,~the restriction
$\res_{G'}^G(\alpha)$
of a proper
\mbox{$4$-torsion}
class~$\alpha$
may
be~$0$
or a proper
\mbox{$2$-torsion}
class, only when
$H^1(G', \frakP)$
does not contain any proper
\mbox{$4$-torsion}
class. Otherwise,~Propositions \ref{typ1_cl} and~\ref{typ2_cl} imply that
$\res_{G'}^G(\alpha)$
is a
\mbox{$4$-torsion}
class of the same type
as~$\alpha$.

\begin{theo}[A~nontriviality criterion for the restriction of a
\mbox{$4$-torsion}
class]
\label{4tors_rest}
Let\/
$G \subseteq (\bbZ/2\bbZ)^4 \!\rtimes\! S_5$
be a subgroup such that there is a proper
\mbox{$4$-torsion}
class\/
$\alpha \in H^1(G, \frakP)$
and let\/
$G' \subseteq G$
be any~subgroup.

\begin{iii}
\item
Suppose~that\/
$\alpha$
is of type~I.
Then\/~$G \subseteq S_3 \times \bbZ/4\bbZ$
and one has the natural surjection\/
$s\colon G \twoheadrightarrow \langle\tau\rangle \cong \bbZ/4\bbZ$.
The~following statements are~true.
\begin{abc}
\item[$\bullet$ ]
If
$s(G') = \bbZ/4\bbZ$
then\/
$\smash{\res_{G'}^G(\alpha)}$
is a proper\/
\mbox{$4$-torsion}
class.
\item[$\bullet$ ]
If
$s(G') = 2\bbZ/4\bbZ$
then\/
$\smash{\res_{G'}^G(\alpha)}$
is a nontrivial\/
\mbox{$2$-torsion}
class.
\item[$\bullet$ ]
If
$s(G') = 0$
then\/
$\smash{\res_{G'}^G(\alpha)}$
is the zero~class.
\end{abc}

\item
If~$\alpha$
is of type~II, suppose that the indices
$1,\ldots,5$
are normalised as in~\ref{type2}. Then,
\begin{abc}
\item[$\bullet$ ]
if\/~$G'$
leaves the orbit\/
$\{1\}$
non-split then
$\smash{\res_{G'}^G(\alpha)}$
is a proper\/
\mbox{$4$-torsion}
class.
\item[$\bullet$ ]
If\/~$G'$
splits the orbit\/
$\{1\}$,
but either of the sets\/
$\{2,4\}$
and\/
$\{3,5\}$
contains at least one non-split orbit then\/
$\smash{\res_{G'}^G(\alpha)}$
is a nontrivial\/
\mbox{$2$-torsion}
class.
\item[$\bullet$ ]
Finally,~if
$G'$
splits the orbit\/
$\{1\}$
and either\/
$\{2,4\}$
or\/
$\{3,5\}$
consists entirely of split orbits then\/
$\smash{\res_{G'}^G(\alpha)}$
is the zero~class.
\end{abc}
\end{iii}\smallskip

\noindent
{\bf Proof.}
{\em
i)
The first assertion follows directly from Proposition~\ref{typ1_cl}.b), together with Theorem~\ref{res_4tors}. Otherwise,~one has, assuming that the indices
$1,\ldots,5$
are normalised as in~\ref{type1},
$\smash{\res_{G'}^G(\alpha) = \iota^{G'}_{\Pic,4}(\overline{[e_1] + [e_2] + [e_3] + 2[e_4]}) = \iota^{G'}_{\Pic,4}(2\overline{[e_4]})}$.
Indeed,~in the situation that
$\smash{s(G') \subsetneqq \bbZ/4\bbZ}$,
the orbits in
$\{1,2,3\}$
are all~split. Lemma~\ref{4_tors}.b) then shows that
$\smash{\res_{G'}^G(\alpha) = \iota^{G'}_{\Pic}(\overline{[e_4]})}$.
Moreover,~$\{4\}$
is a non-split orbit and
$\{5\}$
is another, when 
$s(G') = 2\bbZ/4\bbZ$,
while
$\{4\}$
is a split orbit in the case
that~$s(G') = 0$.
The assertion thus follows from Theorem~\ref{expl_2}.\smallskip

\noindent
ii)
Here,~the first claim directly follows from Proposition~\ref{typ2_cl}.b), together with Theorem~\ref{res_4tors}. On~the other hand, if
$\{1\}$
is a split orbit then Lemma~\ref{4_tors}.b) shows, together with the compatibility of restriction with the boundary homomorphisms,~that
$$\res_{G'}^G(\alpha) = \iota^{G'}_{\Pic,4}(\overline{[e_1] + 2[e_2] + 2[e_4]}) = \iota^{G'}_{\Pic,4}(\overline{2[e_2] + 2[e_4]}) = \iota^{G'}_{\Pic}(\overline{[e_2] + [e_4]}) \,.$$
In~the terminology of the explicit description of
$H^1(G', \frakP)_2$,
given in Theorem~\ref{expl_2},
$\smash{\iota^{G'}_{\Pic}(\overline{[e_2] + [e_4]})}$
corresponds to
$\overline{1} \!\cdot\! \{2,4\}$
or
$\overline{1} \!\cdot\! \{2\} + \overline{1} \!\cdot\! \{4\}$,
depending on whether
$\{2,4\}$
forms one orbit
under~$S' := p(G')$
or~two. The~image of this element in the quotient modulo the split orbits is nonzero if and only if
$\{2,4\}$
is a non-split orbit or, if it consists of two orbits, if at least one of them is non-split. Assume~that this is the~case. Then,~in order to have
$\smash{\iota^{G'}_{\Pic}(\overline{[e_2] + [e_4]}) \neq 0}$,
one needs, in addition, one further non-split orbit, which may only be
$\{3,5\}$,
$\{3\}$,
or~$\{5\}$,
as~required.
}
\eop
\end{theo}

\section{Evaluation of the Brauer classes and examples}

\begin{conv}
In~this section and the next,
$k$
is always assumed to be a number field.
\end{conv}

\begin{defi}
\label{model}
Let~$U$
be a separated scheme that is smooth and of finite type
over~$k$.
Then~by a {\em model\/}
of~$U$,
we mean a separated
scheme~$\calU$
that is of finite type
over~$\calO_k$,
the ring of integers
in~$k$,
such that
$\calU \times_{\Spec \calO_k} \Spec k \cong U$.
\end{defi}

\begin{cau}
In this and the next sections, when considering local evaluation maps
on~$U$,
we always assume that a
model~$\calU$
is~fixed. For~the non-archimedean places, we only consider the restrictions
$\ev_{\alpha,\nu} |_{\calU(\calO_{k,\nu})} \colon \calU(\calO_{k,\nu}) \longrightarrow \bbQ/\bbZ$.
This~is motivated by the results of J.-L.\ Colliot-Th\'el\`ene and F.\ Xu~\cite[\S1]{CX} and the resulting application towards strong approximation. In~particular, we do not care whether perhaps
$\ev_{\alpha,\nu}$
is constant on the whole
of~$U(k_\nu)$,
for certain non-archimedean
places~$\nu$.
\end{cau}

\begin{ex}
Let~$q_1, \ldots, q_r, h \in k[X_0, \ldots, X_N]$
be homogeneous forms such that
$U := V(q_1, \ldots, q_r) \!\setminus\! V(h) \subset \Pb_k^N$
is a hypersurface complement of a complete~intersection. Then,~for forms
$\smash{\widetilde{q}_1, \ldots, \widetilde{q}_r, \widetilde{h} \in \calO_k[X_0, \ldots, X_N]}$
that are just scalar multiples
of~$q_1, \ldots, q_r$,
and~$h$,
respectively, the subscheme
$\smash{\calU := V(\widetilde{q}_1, \ldots, \widetilde{q}_r) \!\setminus\! V(\widetilde{h}) \subset \Pb_{\calO_k}^N}$
is a model
of~$U$.
\end{ex}

\subsubsection*{Constancy of the local evaluation map}\leavevmode

\noindent
Let~$U$
be as above,
$\alpha \in \Br(U)$,
and
$\calU$
be a model
of~$U$.
In~such a situation, based on the fundamental results~\cite{Br07} of M.~Bright, one may often indicate an explicit finite
set~$S_{\calU ,\alpha}$
such~that
$\ev_{\alpha,\nu} |_{\calU(\calO_{k,\nu})} \colon \calU(\calO_{k,\nu}) \longrightarrow \bbQ/\bbZ$
is the zero map for each
place~$\smash{\nu \in \Sigma_k \!\setminus\! S_{\calU ,\alpha}}$.

\begin{theo}[Constancy of the local evaluation map]
\label{const}
Let\/~$X$
be an irreducible scheme being proper and smooth
over\/~$k$
and having the property that\/
$\Pic^0 X_{\overline{k}} = 0$,
(i.e., that
$\Pic X_{\overline{k}} = \NS X_{\overline{k}}$).
Denote~by\/
$l$
the field of definition
of\/~$\NS X_{\overline{k}}$.\smallskip

\noindent
Moreover,~let\/
$\pi\colon B \twoheadrightarrow X$
be an\/
\mbox{$\Ab^n$-bundle},
for some\/
$n \geq 0$,
and\/
$U \subset B$
an open~subscheme. Suppose~that\/
$\Pic U_{\overline{k}}$
is torsion-free.\smallskip

\noindent
Then,~for each prime ideal\/
$\frakp \subset \calO_k$
that is unramified
in\/~$l$,
every model\/
$\calU$
of\/~$U$
that is smooth
above\/~$\frakp$
and whose special
fibre~$\calU_\frakp$
is irreducible, and any algebraic Brauer~class\/
$\alpha \in \Br_1(U)$,
the local evaluation map\/
$\ev_{\alpha,\frakp} |_{\calU(\calO_{k,\frakp})} \colon \calU(\calO_{k,\frakp}) \to \bbQ/\bbZ$
is~constant.\medskip

\noindent
{\bf Proof.}
{\em
One~has
$\Pic B_{\overline{k}} = \Pic X_{\overline{k}}$,
according to \cite[Exercise~II.6.3.a)]{Ha}.
Consequently,~$\Pic U_{\overline{k}}$
is a quotient
of~$\Pic X_{\overline{k}}$.
In~particular,
$\smash{\Pic U_{\overline{k}}}$
is acted upon trivially by
$\smash{\Gal(\overline{k}/l)}$,
which shows that
$$\Br_1(U_l)/\Br_0(U_l) = H^1(\Gal(\overline{k}/l), \Pic U_{\overline{k}}) = \Hom(\Gal(\overline{k}/l), \Pic U_{\overline{k}}) = 0 \, .$$
In~particular, one has
$\alpha |_{U_l} \in \Br_0(U_l)$.
Let~us fix a
prime~$\frakq$
of~$l$
that lies
above~$\frakp$.
Then,~after possibly adding a constant Brauer~class, one may assume that
$\alpha |_{U_{l_\frakq}} = 0$.

Furthermore,~according to our assumptions,
$\smash{\calU_{\calO_{k,\frakp}}}$
fulfils M.~Bright's \cite[page~3]{Br07} Condition~(*). (Note that the arguments in \cite{Br07} work over
$\calO_{k,\frakp}$,
as well as
over~$\bbZ_p$.)\break
Therefore, \cite[Propositions~6 and~3]{Br07} show that the local evaluation map
$\ev_{\alpha,\frakp} |_{\calU(\calO_{k,\frakp})}$
factors via
$\smash{H^2(\Gal(l_\frakq/k_\frakp), \calO_{l_\frakq}^*)}$.
The~latter cohomology group, however, vanishes in the unramified case, according to \cite[Chap.~V, \S2]{Se}.%
}%
\eop
\end{theo}

\begin{rem}[Good reduction implies being unramified]
Let~$X$
be a surface and
$\frakp \subset \calO_k$
be a prime ideal such that there is a model\/
$\calX$
of~$X$
that is proper
over~$\calO_k$
and smooth
above~$\frakp$.
Then~$\frakp$
is unramified in the field of definition
of~$\NS X_{\overline{k}}$.
Indeed,~this follows from the smooth specialisation theorem for \'etale cohomology~\cite[Exp.\ XVI, Corollaire~2.3]{SGA4}. Cf.~\cite[Lemma~2.3.5]{CEJ} for a detailed~argument.
\end{rem}

\begin{rems}
\begin{iii}
\item
A~proper del Pezzo surface of degree four fulfils the general assumptions made on~$X$.

Thus,~let
$X$
be a proper del Pezzo surface of degree four having a model
$\calX$
that is irreducible and has an irreducible special
fibre~$\calX_\frakp$.
Moreover,~let
$\calU \subset \calX$
be an open subscheme that excludes the singular points
of~$\calX_\frakp$.
Then~the local evaluation map
$\ev_{\alpha,\frakp} |_{\calU(\calO_{k,\frakp})} \colon \calU(\calO_{k,\frakp}) \to \bbQ/\bbZ$
is~constant, as long as
$\frakp$
is unramified in the
field~$l$
of definition
of~$\NS X_{\overline{k}}$.
\item
A~non-singular space quadric
$X$
fulfils the general assumptions,~too. Having~taken out the cusp, the cone
above~$X$
is an
\mbox{$\Ab^1$-bundle}, to which Theorem~\ref{const} applies. We~make use of this in Corollary~\ref{quad_const},~below.
\end{iii}
\end{rems}

\subsubsection*{Corestriction}\leavevmode

\noindent
Let~$U' \to U$
be a finite \'etale morphism of~schemes. Then~there is a natural corestriction homomorphism
$\cores_U^{U'}\colon \Br U' \to \Br U$.
In~the case of affine schemes, a construction is described in \cite[Chapter~8]{Sa} and~\cite[Chapter~II, \S1]{Gre}. As~that commutes with arbitrary base change~\cite[Theorem~8.1.d)]{Sa}, the corestriction extends directly to the setting of general~schemes.

\begin{fac}
\label{comp}
In~the case of a finite [separable] field extension, one recovers the usual corestriction in Galois~cohomology.\medskip

\noindent
{\bf Proof.}
{\em
This is~\cite[Chapter~II, Proposition~1.6]{Gre}
}
\eop
\end{fac}

The~case relevant to us is that
$U$
is an open del Pezzo surface of degree four
over~$k$
(Def.~\ref{open_dP4}) and
$U' := U_l$
is a base extension
of~$U$.
There~are the following two~results.

\begin{lem}
\label{cores_Br_H1}
Let\/~$U$
be an open del Pezzo surface of degree four
over\/~$k$
and\/
$l$
a finite extension~field. Then the~diagram
$$
\xymatrix@R=10pt{
\Br_1(U_l)/\Br_0(U_l) \ar@{->}[r]^{\overline\cores} & \Br_1(U)/\Br_0(U) \\
H^1(\Gal(\overline{k}/l), \Pic U_{\overline{k}}) \ar@{->}[r]^\cores \ar@{->}[u]^\cong & H^1(\Gal(\overline{k}/k), \Pic U_{\overline{k}}) \ar@{->}[u]_\cong
}
$$
commutes, in which the upper arrow is induced by corestriction, the lower arrow is the corestriction in Galois cohomology, and the upwards arrows are the isomorphisms induced by the Hochschild--Serre spectral sequence (cf.~Section~\ref{sec_Brauer}).\medskip

\noindent
{\bf Proof.}
{\em
The diagram
$$
\xymatrix@R=10pt{
H^2(\Gal(\overline{k}/l), \overline{k}(U)^*/\overline{k}^*) \ar@{->}[r]^\cores & H^2(\Gal(\overline{k}/k), \overline{k}(U)^*/\overline{k}^*) \\
H^1(\Gal(\overline{k}/l), \Pic U_{\overline{k}}) \ar@{->}[r]^\cores \ar@{->}[u] & H^1(\Gal(\overline{k}/k), \Pic U_{\overline{k}}) \ar@{->}[u] \,,
}
$$
in which the upwards arrows are the boundary homomorphisms associated with the short exact~sequence
$$0 \longrightarrow \overline{k}(U)^*/\overline{k}{}^* \longrightarrow \Div U_{\overline{k}} \longrightarrow \Pic U_{\overline{k}} \longrightarrow 0 \,,$$
commutes. Indeed,~corestriction commutes with boundary homomorphisms. On~the other hand, the natural~diagram
\arraycolsep1pt
$$
\xymatrix@R=10pt@C=20pt{
H^2(\Gal(\overline{k}/l), \overline{k}(U)^*/\overline{k}^*) \ar@{->}[r]^\cores & H^2(\Gal(\overline{k}/k), \overline{k}(U)^*/\overline{k}^*) \\
H^2(\Gal(\overline{k}/l), \overline{k}(U)^*) / H^2(\Gal(\overline{k}/l), \overline{k}^*) \ar@{->}[r]^{\overline\cores} \ar@{->}[u] & H^2(\Gal(\overline{k}/k), \overline{k}(U)^*) / H^2(\Gal(\overline{k}/k), \overline{k}^*) \ar@{->}[u] \\
\Br_1(U_l)/\Br_0(U_l) \ar@{->}[r]^{\overline\cores} \ar@{->}[u] & \Br_1(U)/\Br_0(U) \ar@{->}[u]
}
$$
commutes as~well. For~the lower square, this is Fact~\ref{comp}, while commutativity of the upper square is a particular case of the fact that corestriction is a morphism of~functors.
The~assertion immediately follows from this, in view of the elementary descriptions of the isomorphism
$H^1(\Gal(\overline{k}/k), \Pic U_{\overline{k}}) \to \Br_1(U)/\Br_0(U)$
and its analogue
over~$l$,
given in Facts~\ref{Lichtenbaum}.i) and~ii).
}
\eop
\end{lem}

\begin{lem}
\label{cores_eval}
Let\/~$U$
be any scheme
over\/~$k$
and\/
$l$
a finite extension field. Then,~for every\/
$\alpha \in \Br U_l$,
every place\/
$\nu$
of\/~$k$,
and every
$x \in U(k_\nu)$,
one~has
$$\ev_{\cores_U^{U_l}(\alpha), \nu}(x) = \sum_{w|\nu} \ev_{\alpha, w}(x_{l_w}) \,,$$
the sum running over all places
$w$
of\/~$l$
that lie
above\/~$\nu$.\medskip

\noindent
{\bf Proof.}
{\em
There is the Cartesian~diagram\vspace{-2mm}
$$
\xymatrix@R=10pt{
\smash{\coprod\limits_{w|\nu}} \Spec l_w \smash{\ar@{->}[r]^{\;\;\;\;\;\;(x_{l_w})_w}}\ar@{->}[d]& U_l \ar@{->}[d] \\
\Spec k_\nu \ar@{->}[r]^{\;\;\;x} & U \hsmash{\, .}
}
$$
As~the corestriction commutes with base change,
$x^* \cores_U^{U_l}(\alpha) = \smash{\sum\limits_{w|\nu}} \cores_{k_\nu}^{l_w}(x_{l_w}^* \alpha)$.
Consequently,
\begin{align*}
\ev_{\cores_U^{U_l}(\alpha), \nu}(x) = \inv_{k_\nu} x^* \!\cores_U^{U_l}(\alpha) = \sum_{w|\nu} \inv_{k_\nu} \cores_{k_\nu}^{l_w}(x_{l_w}^* \alpha) = \sum_{w|\nu} \inv_{l_w} x_{l_w}^* \alpha& \\[-2mm]
= \sum_{w|\nu} \ev_{\alpha, w}&(x_{l_w}) \,.
\end{align*}
Here,~the next-to-last equality is one of the fundamental properties of the invariant of a Brauer class of a local field \cite[Chapitre~XI, \S2, Proposition~1.ii)]{Se}.
}
\eop
\end{lem}

\subsubsection*{Cyclic classes}\leavevmode

\noindent
A class
$c \in H^1(\Gal(\overline{k}/k), \Pic U_{\overline{k}})$
is called {\em cyclic\/} if there is a Galois extension
$l/k$
with cyclic Galois group such that the restriction
of~$c$
to
$\smash{H^1(\Gal(\overline{k}/l), \Pic U_{\overline{k}})}$
vanishes. In the abstract setting of section~\ref{sec_Brauer}, this means that a class in
$H^1(G, \frakP)$
is cyclic if there is a normal subgroup
$H \triangleleft G$
with cyclic quotient
$G/H$
such that the restriction to
$H^1(H, \frakP)$~vanishes.

\begin{prop}[Manin's class field theoretic method for evaluation]
Let\/~$U$
be an open del Pezzo surface of degree four
over\/~$k$.
Assume~that\/
$U$
has an adelic~point. Moreover,~let\/
$c \in H^1(\Gal(\overline{k}/k), \Pic U_{\overline{k}})$
be a cyclic class and let\/
$l/k$
be a cyclic Galois extension that
annihilates\/~$c$.

\begin{abc}
\item
Then\/~$c$
is the image under inflation of a class\/
$\widetilde{c} \in H^1(\Gal(l/k), \Pic U_l)$.
\item
There is a natural isomorphism
$$\iota\colon H^{-1}(\Gal(l/k), \Pic U_l) \longrightarrow [N_{l/k} \Div U_l \cap \Div_0 U] / N_{l/k} \Div_0 U_l \, ,$$
where\/
$\Div$
denotes the Galois module of all divisors,
$\Div_0$
the submodule of principal divisors, and\/
$N$
the~norm.
\item
Fix a generator\/
$\sigma \in \Gal(l/k)$
and put\/
$c^- := \cl_\sigma \cup\; \widetilde{c} \in H^{-1}(\Gal(l/k), \Pic U_l)$,
for\/
$\cl_\sigma \in H^{-2}(\Gal(l/k), \bbZ)$
the fundamental~class. Moreover,~let\/
$f \in k(U)$
be a rational function representing a principal divisor in the residue class\/
$\iota(c^-)$.

Then there is a lift\/
$\alpha \in \Br(U)$
of\/~$c$
such that, for
every place\/
$\nu$
of\/~$k$
and
every\/
$x \in U(k_\nu)$,
at which\/
$f$
is defined and nonzero, 
$\smash{\ev_{\alpha,\nu}(x) = (\frac{i}{[l_w \!:\! k_\nu]} \bmod \bbZ)}$
if and only~if
$$(f(x), l_w/k_v) = \sigma^i \, .$$
Here,~$w$
is any place
of\/~$l$
lying
above\/~$\nu$
and\/
$( . , l_w/k_v)$
denotes the norm residue symbol\/~\cite[Chapitre~XIII, \S4]{Se}.
\end{abc}\smallskip

\noindent
{\bf Proof.}
{\em
a)
A direct application of the inflation-restriction sequence yields that
$c$
is the image of a class from
$\smash{H^1(\Gal(l/k), (\Pic U_{\overline{k}})^{\Gal(\overline{k}/l)})}$.
However,~as
$U$
has an adelic point, one has
$\smash{(\Pic U_{\overline{k}})^{\Gal(\overline{k}/l)} = \Pic U_l}$,
according to~\cite[Proposition~2.21]{Br02}.\smallskip

\noindent
b)
Yu.\,I.\ Manin~\cite[Proposition~31.3]{Ma} was the first, who came up with such an isomorphism. The~non-proper case, which is not substantially different,  is considered, for example, in~\cite[Chapter~III, Lemma~8.19]{Ja}.\smallskip

\noindent
c)
This, finally, is Manin's original method for the evaluation of a Brauer class \mbox{\cite[\S45.2]{Ma}}. Cf.~\cite[Chapter~IV, Section~4]{Ja} for a few more~details.
}
\eop
\end{prop}

\begin{rems}
\begin{iii}
\item
Under the canonical isomorphism
$$H^{-2}(\Gal(l/k), \bbZ) \stackrel{\cong}{\longrightarrow} \Gal(l/k)^\ab = \Gal(l/k) \,,$$
the fundamental
class~$\cl_\sigma$
is that mapped
to~$\sigma$.
\item
A~class in 
$H^{-1}(\Gal(l/k), \Pic U_l)$
is represented by an element
$\calL \in \Pic U_l$
of trivial~norm. I.e.,~by a divisor
$L \in \Div U_l$
whose norm is a principal divisor,
$N_{l/k}L = \div f \in \Div_0 U$.
One~then indeed has that
$$\iota(\overline\calL) = \overline{\div f} \,.$$
This follows immediately from the construction
of~$\iota$
given in the proof of \cite[Proposition~31.3]{Ma}.
\item
The ambiguity that the rational function
$f$
itself is determined
by~$\div f$
only up to a constant factor is absorbed by the possible choices of a lift of
$c \in H^1(\Gal(\overline{k}/k), \Pic U_{\overline{k}}) \cong \Br_1(U)/\im\Br(k)$
to~$\Br_1(U)$.
\item
At~a
point~$x$
where
$f$
is undefined due to a pole or at which
$f(x) = 0$,
one can often nevertheless determine 
$\ev_{\alpha,\nu}(x)$
using the continuity
of~$\ev_{\alpha,\nu}$.
Moreover,~there exists a rational function
$g \in l(U_l)$
such that
$f \!\cdot\! N_{l/k}(g)$
is defined and nonzero
at~$x$.
Indeed,~$U_l$
is a regular scheme and therefore the divisor
$(-L)$
is locally~principal.
\end{iii}
\end{rems}

\subsubsection*{$2$-torsion classes of type~I}\leavevmode

\begin{lem}
\label{1_cocyc_type1}
Let\/~$G \subseteq (\bbZ/2\bbZ)^4 \!\rtimes\! S_5$
be a~subgroup. Then

\begin{iii}
\item
a\/~\mbox{$2$-torsion}
class\/
$\smash{\iota_{\Pic}(\overline{[e_i]}) \in H^1(G, \frakP)}$
of type~I is always~cyclic. In~fact, let\/
$G' \subseteq G$
be the stabiliser
of\/~$[e_i]$.
Then\/~$G/G' \cong \bbZ/2\bbZ$
and\/
$\smash{\iota_{\Pic}(\overline{[e_i]})}$
is the image under inflation of the class\/
$\smash{\widetilde{c}\in H^1(G/G', \bbZ \!\cdot\! [e_i])}$
of the\/
\mbox{$1$-cocycle}
$$
\sigma \mapsto \left\{
\begin{array}{rl}
     0, & \text{ for\/ } \sigma \in G/G' \text{ the neutral element\,,} \\
-[e_i], & \text{ for\/ } \sigma \in G/G' \text{ the nontrivial element\,.}
\end{array}
\right.
$$
\item
The class\/
$\cl_\sigma \cup\; \widetilde{c} \in H^{-1}(G/G', \bbZ \!\cdot\! [e_i]))$
is represented by\/
$-[e_i]$.
\end{iii}\smallskip

\noindent
{\bf Proof.}
{\em
i)
Observe at first that
$\bbZ \!\cdot\! [e_i] \subseteq \frakP^{G'}$,
according to the definition
of~$G'$.
The~assertion then follows from the explicit description of
$\iota_{\Pic}$,
given in Lemma~\ref{H1_expl}.a.ii), together with a calculation in~cocycles.\smallskip

\noindent
ii)~is a consequence of~i) and the explicit formula for the cup product, given in \cite[Chapitre~XI, Annexe, Lemme~3]{Se}.
}
\eop
\end{lem}

\begin{ttt}[Evaluation--Planes tangent to a degenerate quadric]
\label{eval_tang}
Let~$U \subset X = V(q_1, q_2)$
be an open del Pezzo surface of degree four
over~$k$
(cf.~Definition~\ref{open_dP4}), on which there is an algebraic
\mbox{$2$-torsion}
Brauer
class~$\alpha$
of type~I.
Then~$\smash{\alpha = \iota_{\Pic}(\overline{[e_i]})}$,
for a certain linear system
$e_i$
of~conics. In~particular, one of the five degenerate quadrics in the pencil
$\smash{(\lambda q_1 + \mu q_2)_{(\mu:\nu) \in \Pb^1}}$,
that containing the entire
system~$e_i$,
is
\mbox{$k$-rational}.
Let~us say that such a quadric {\em induces\/} the
class~$\alpha$.

Without~restriction, assume that
$V(q_1)$
induces the Brauer
class~$\alpha$.
We~then take a hyperplane
$V(t)$
that is tangent
to~$V(q_1)$.
Over~a quadratic extension field
$\smash{l = k(\sqrt{d})}$,
the intersection
$V(q_1) \cap V(t)$
splits into two planes, so that
$\div t$
decomposes into two components, a
conic~$C_i$ from the class
$[e_i]$
and its conjugate
from~$\smash{[e'_i]}$.
Hence,~$\div t = N_{l/k}C_i$
and~$\div 1/t = N_{l/k}(-C_i)$.
The evaluation of the Brauer class is therefore given by the norm residue symbols
$$\big( 1/t(x), k_\nu(\sqrt{d})/k_\nu \big) = \big( t(x), k_\nu(\sqrt{d})/k_\nu \big) \,.$$
Note~that when the splitting occurs already
over~$k$
then
$\{i\}$
is a split orbit and
$\smash{\iota_{\Pic}(\overline{[e_i]}) = 0}$,
according to Theorem~\ref{expl_2}.
\end{ttt}

\begin{rems}
\begin{iii}
\item
(Some kind of normal form){\bf .}
One~may write the
\mbox{rank-$4$}
quadric in the form
$l_1l_2 - l_3^2 + dl_4^2$,
for four linearly independent linear forms
$l_1, \ldots, l_4$.
The~evaluation is then given by the norm residue symbols
$\smash{(l_1(x), k_\nu(\sqrt{d})/k_\nu)}$.
\item
This type of Brauer class has been evaluated before in exactly the same~way. Cf., e.g., \cite[Example~8.1]{JS}, where two linearly independent
\mbox{$2$-torsion}
classes of type~I are~occurring. Furthermore,~the next lemma shows that
\mbox{$2$-torsion}
classes of type~I are restrictions of Brauer classes on open
rank~$4$
quadrics, which have been evaluated in the same way by J.-L.\ Colliot-Th\'el\`ene and F.~Xu \cite[\S5.8]{CX},~already.
\end{iii}
\end{rems}

\begin{lem}
\label{eval_quadr}
Let\/~$U = X \!\setminus\! H \subset X \subset \Pb^4_k$,
for\/
$H \subset \Pb^4_k$
a hyperplane, be an open del Pezzo surface of degree four
over\/~$k$
(cf.~Definition~\ref{open_dP4}) having a\/
\mbox{$2$-torsion}
class\/
$\alpha \in \Br_1(U)/\Br_0(U)$
of type~I. Suppose~that\/
$\smash{V(q) =: Q \supset X}$
is a\/
\mbox{$k$-rational}
rank\/~$4$
quadric inducing the
class\/~$\alpha$,
and write\/
$\smash{x_0 \in Q}$
for its~cusp.

\begin{abc}
\item
Then\/
$\Pic (Q \!\setminus\! (H \!\cup\! \{x_0\}))_{\overline{k}} \cong \bbZ$.
\item
There~is a class\/
$\smash{\widetilde\alpha \in \Br_1(Q \!\setminus\! (H \!\cup\! \{x_0\}))/\Br_0(Q \!\setminus\! (H \!\cup\! \{x_0\}))}$,
the restriction
to\/~$U$
of which
is\/~$\alpha$.
\end{abc}\smallskip

\noindent
{\bf Proof.}
{\em
a)
One has
$\smash{\Pic (Q \!\setminus\! \{x_0\})_{\overline{k}} = \Pic Q'_{\overline{k}}}$,
for
$Q'$
the non-singular space quadric, the cone of which is
$\smash{(Q \!\setminus\! \{x_0\})_{\overline{k}}}$.
Hence,~$\Pic (Q \!\setminus\! \{x_0\})_{\overline{k}} = \bbZ^2$.
Under~this isomorphism, the class of the hyperplane section
$H$
is mapped
to~$(1,1)$,
which implies that
$\smash{\Pic (Q \!\setminus\! (H \!\cup\! \{x_0\}))_{\overline{k}} = \langle e\rangle}$,
for~$e$
the class of one of the two linear systems of~planes. Moreover,~the generator
$e$
is of infinite~order. Indeed,~its restriction
to~$U$
is the class of one of the ten linear systems of~conics.\smallskip

\noindent
b)
Let~$l$
be the quadratic extension field that splits the two linear systems of planes
on~$Q$.
Then~the nontrivial element
of~$\Gal(l/k)$
maps
$e$
to~$(-e)$,
so~that
$$
\sigma \mapsto \left\{
\begin{array}{rl}
 0, & \text{ if\/ } \sigma \in \Gal(\overline{k}/k) \text{ induces the identity on } l \,, \\
-e, & \text{ if\/ } \sigma \in \Gal(\overline{k}/k) \text{ induces the conjugation on } l \,,
\end{array}
\right.
$$
is a
\mbox{$1$-cocycle}
defining a class
$\smash{\widetilde\alpha \in H^1(\Gal(\overline{k}/k), \Pic (Q \!\setminus\! (H \!\cup\! \{x_0\}))_{\overline{k}})}$.
In~view of Lemma \ref{1_cocyc_type1}.i), it is evident that
$\alpha$
is the restriction
of~$\smash{\widetilde\alpha}$
to~$U$.
}
\eop
\end{lem}

\begin{coro}
\label{quad_const}
Let\/~$q_1$
and\/~$q_2$
be two quadratic forms
and\/~$h$
a linear form in five variables
over\/~$\calO_k$,
$\smash{\calX = V(q_1, q_2) \subset \Pb^4_{\calO_k}}$
the associated scheme, and\/
$X \subset \Pb^4_k$
its generic~fibre. Suppose~that\/
$\calU = \calX \!\setminus\! V(h)$
is a model of an open del Pezzo
surface\/~$U$
of degree~four of the kind that there is a\/
\mbox{$2$-torsion}
class\/
$\alpha \in \Br_1(U)/\Br_0(U)$
of type~I. Let\/
$\smash{V(q_1) \supset X}$
be a degenerate quadric that induces the
class\/~$\alpha$.\smallskip

\noindent
Let\/
$\frakp \subset \calO_k$
be a prime ideal such that the reduction
of\/~$q_1$
modulo\/~$\frakp$
is of
rank\/~$4$
and suppose that the cusp\/
$x_0 \in V(q_1)_\frakp$
is not a point on the special fibre\/
$\calU_\frakp$.
Then~the local evaluation map\/
$\ev_{\alpha,\frakp} |_{\calU(\calO_{k,\frakp})} \colon \calU(\calO_{k,\frakp}) \to \bbQ/\bbZ$
is~constant.\medskip

\noindent
{\bf Proof.}
{\em
Let~$\calC$
be the Zariski closure
in~$\smash{\Pb_{\calO_k}^4}$
of the cusp on
$V(q_1)_k$,
put
$\calW := V(q_1) \!\setminus\! (V(h) \!\cup\! \calC)$,
and let
$W$
be the generic fibre
of~$\calW$.
Then~$\calW \supset \calU$.
Moreover,~according to Lemma~\ref{eval_quadr},
$\alpha \in \Br_1(U)$
is the restriction of an algebraic Brauer class
$\widetilde\alpha \in \Br_1(W)$.
Thus,~for every
$x \in \calU(\calO_{k,\frakp})$,
one has
$\ev_{\alpha}(x) = \ev_{\widetilde\alpha}(x)$,
so that it suffices to show that
$\ev_{\widetilde\alpha} |_{\calW(\calO_{k,\frakp})} \colon \calW(\calO_{k,\frakp}) \to \bbQ/\bbZ$
is~constant.

And~this, in fact, is a direct consequence of Theorem~\ref{const}. Indeed,~for the generic fibre, one knows that
$\smash{\big(V(q_1) \!\setminus\! \calC\big)_k}$
is an
\mbox{$\Ab^1$-bundle}
over a non-singular space~quadric.
Furthermore,~according to Lemma~\ref{eval_quadr}.a),
$\Pic W_{\overline{k}} \cong \bbZ$,
which is torsion-free.
Finally,~by assumption,
$\calW$
is smooth
above~$\frakp$
and its special fibre
$\calW_\frakp$
is~irreducible.%
}%
\eop
\end{coro}

\subsubsection*{$2$-torsion classes of type~II}\leavevmode

\begin{figure}[H]
\begin{center}
\begin{picture}(50,55)
\linethickness{0.8pt}
\setlength\unitlength{1.7pt}
\put(0,3){\circle{1.5}}
\put(0,26){\circle{1.5}}
\put(25,3){\circle{1.5}}
\put(25,26){\circle{1.5}}
\linethickness{1pt}
\put(-2,29){$e_i$}\put(2,26){\vector(1,0){21}}
\put(23,29){$e_j$}\put(22.5,24){\vector(-1,-1){20}}
\put(-2,-3){$e'_i$}\put(2,3){\vector(1,0){21}}
\put(23,-3){$e'_j$}\put(22.5,5){\vector(-1,1){20}}
\end{picture}\hspace{2.6cm}
\begin{picture}(50,55)
\linethickness{0.8pt}
\setlength\unitlength{1.7pt}
\put(0,3){\circle{1.5}}
\put(0,26){\circle{1.5}}
\put(25,3){\circle{1.5}}
\put(25,26){\circle{1.5}}
\linethickness{1pt}
\put(-2,29){$e_i$}\put(0,22){\vector(0,-1){17.5}}
\put(23,29){$e_j$}
\put(-2,-3){$e'_i$}\put(0,7){\vector(0,1){17.5}}
\put(23,-3){$e'_j$}
\end{picture}
\end{center}\vskip-\medskipamount

\caption{Two generators for the dihedral group $D_4$}\vskip-\medskipamount
\label{fig_8}
\end{figure}
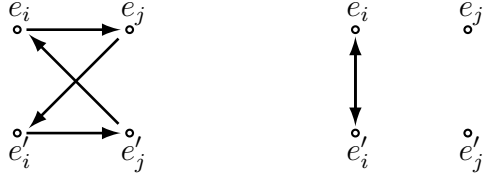

\noindent
\mbox{$2$-torsion}
classes of type~II are, in general, non-cyclic. Indeed,~let
$G \subseteq (\bbZ/2\bbZ)^4 \!\rtimes\! S_5$
be a subgroup and let
$\smash{\alpha := \iota_{\Pic}(\overline{[e_i]+[e_j]})}$
be such a~class. The~quotient
group~$G'$
which 
$G$
induces on the orbit
$\{[e_i], [e_j], [e'_i], [e'_j]\}$~\cite[\S1.6]{DM}
is a subgroup of the dihedral
group~$D_4$,
generated by the two elements depicted in Figure~\ref{fig_8}.

When~exactly this group occurs then there are two nontrivial subgroups
of~$G'$
splitting the~orbit. Either~is of
order~$2$.
They~are generated by the elements shown~below.
One~readily sees that the two are conjugate to each other
in~$D_4$,
which shows in particular that they are non-normal. Thus,~the only normal subgroup
of~$D_4$
that splits the orbit is the trivial~group.
However,~$D_4$
is clearly not cyclic, not even~abelian.

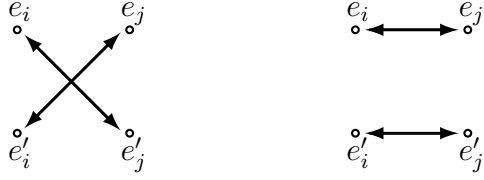
\begin{figure}[H]
\begin{center}
\begin{picture}(50,55)
\linethickness{0.8pt}
\setlength\unitlength{1.7pt}
\put(0,3){\circle{1.5}}
\put(0,26){\circle{1.5}}
\put(25,3){\circle{1.5}}
\put(25,26){\circle{1.5}}
\linethickness{1pt}
\put(-2,29){$e_i$}\put(3,23.5){\vector(1,-1){20}}
\put(23,29){$e_j$}\put(21.5,24){\vector(-1,-1){20}}
\put(-2,-3){$e'_i$}\put(3,5.5){\vector(1,1){20}}
\put(23,-3){$e'_j$}\put(21.5,5){\vector(-1,1){20}}
\end{picture}\hspace{2.6cm}
\begin{picture}(50,55)
\linethickness{0.8pt}
\setlength\unitlength{1.7pt}
\put(0,3){\circle{1.5}}
\put(0,26){\circle{1.5}}
\put(25,3){\circle{1.5}}
\put(25,26){\circle{1.5}}
\linethickness{1pt}
\put(-2,29){$e_i$}\put(3,26){\vector(1,0){20}}
\put(23,29){$e_j$}\put(21,26){\vector(-1,0){19}}
\put(-2,-3){$e'_i$}\put(3,3){\vector(1,0){20}}
\put(23,-3){$e'_j$}\put(21,3){\vector(-1,0){19}}
\end{picture}
\end{center}\vskip-\medskipamount

\caption{Subgroups of $D_4$ splitting the orbit}\vskip-\medskipamount
\end{figure}

\begin{lem}
\label{cores_norm}
Let\/~$G \subseteq (\bbZ/2\bbZ)^4 \!\rtimes\! S_5$
and let\/
$\smash{\alpha := \iota_{\Pic}(\overline{[e_i]+[e_j]}) \in H^1(G, \frakP)}$
be a\/
\mbox{$2$-torsion}
class of type~II.
Put\/~$H \subset G$
to be the stabiliser
of\/~$i$
under the operation
of\/~$p(G)$
on\/~$\{1,\ldots,5\}$.
Then

\begin{iii}
\item
$H \subset G$
is a subgroup of
index\/~$2$. 
\item
One~has\/
$\smash{\alpha = \cores_H^G \iota_{\Pic}^H(\overline{[e_i]})}$.
\end{iii}\smallskip

\noindent
{\bf Proof.}
{\em
i)
As~$\iota_{\Pic}(\overline{[e_i]\!+\![e_j]})$
is a
\mbox{$2$-torsion}
class of type~II, the subset
$\{i,j\} \!\subset\! \{1,\ldots,5\}$
is a
\mbox{$G$-orbit}.\smallskip

\noindent
ii)
This~is an immediate consequence of Lemma~\ref{H1_expl}.b). Indeed,
$N_H^G(\overline{[e_i]}) = \overline{[e_i]+[e_j]}$,
for
$N_H^G\colon (\frakP/2\frakP)^H \to (\frakP/2\frakP)^G$
the norm~map.
}
\eop
\end{lem}

\begin{ttt}[Evaluation--Corestriction]
\label{eval_cores}
Let~$U \subset X = V(q_1, q_2)$
be an open del Pezzo surface of degree four
over~$k$
(cf.~Definition~\ref{open_dP4}), on which there is an algebraic
\mbox{$2$-torsion}
Brauer
class~$\alpha$
of type~II.
Then~$\smash{\alpha = \iota_{\Pic}(\overline{[e_i]\!+\![e_j]})}$,
for
$e_i$
and~$e_j$
two conjugate linear systems of~conics. In~particular, two of the five degenerate quadrics in the pencil
$\smash{(\lambda q_1 + \mu q_2)_{(\lambda:\mu) \in \Pb^1},}$
that containing the entire
system~$e_i$
and that
containing~$e_j$,
form a
\mbox{$\smash{\Gal(\overline{k}/k)}$-orbit}.
Let~us say that such a pair of quadrics {\em induces\/} the
class~$\alpha$.

Without~restriction, assume that the pair
$\{V(q_1), V(q_2)\}$
induces~$\alpha$.
Then,~after scaling by constants,
$q_1$
and~$q_2$
are defined over a quadratic extension
field~$k'$ 
and conjugate to each~other. Take~a hyperplane
$V(t)$,
for~$t$
a linear form defined
over~$k'$,
that is tangent
to~$V(q_1)$.
Over~a further quadratic extension
$\smash{l = k'(\sqrt{d})}$,
the intersection
$V(q_1) \cap V(t)$
then splits into two~planes.

According~to Lemma~\ref{cores_norm}.ii) together with Lemma~\ref{cores_Br_H1}, the Brauer
class~$\alpha$
is the corestriction of the
\mbox{$2$-torsion}
Brauer class
on~$U_{k'}$
induced by the degenerate
quadric~$V(q_1)$.
Lemma~\ref{cores_eval} shows that the evaluation of the Brauer class at a
place~$\nu$
is hence given~by
$$\sum_{w|\nu} (t(x), k'_w(\sqrt{d})/k'_w) \,.$$
Note~that, as
$k'/k$
is a quadratic extension, the sum is just one summand when
$\nu$
is inert or ramified
in~$k'$.
In~the split case, it has two~summands.
\end{ttt}

\begin{ex}
Let~$X \subset \Pb^4_\bbQ$
be given by the system of~equations
\begin{align*}
X_0^2 + 2X_0X_2 - 4X_0X_4 - 2X_1X_2 + 4X_1X_4 - X_1X_3 - X_2^2 &= 0 \,, \\
-2X_0X_4 + X_1X_4 + X_2X_3 - 2X_4^2 &= 0 \,,
\end{align*}
$\calX \subset \Pb^4_\bbZ$
the subscheme that is defined by the same system of equations
as~$X$,
and put
$U := X \!\setminus\! H$
and
$\calU := \calX \!\setminus\! \calH$,
for
$\smash{H := V(X_0) \subset \Pb^4_\bbQ}$
and~$\smash{\calH := V(X_0) \subset \Pb^4_\bbZ}$.

A~point search using a variant of Elkies' method~\cite{El} provides almost 25\,000
\mbox{$\bbQ$-rational}
points of height up to
$1000$,
among which 27 are
\mbox{$\bbZ$-integral},
for example
$(1\!:\!1\!:\!0\!:\!1\!:\!0)$,
$(1\!:\!1\!:\!-2\!:\!-3\!:\!-2)$,
and
$(1\!:\!331\!:\!49\!:\!900\!:\!252)$.

The~16 lines
on~$X$
are defined over the Galois hull of
$\smash{l' := \bbQ\big(\!\sqrt{2\sqrt{5}+4}\big)}$,
which is a field of
degree~$8$
with Galois
group~$D_4$.
The~second of the equations
defining~$X$
gives rise to a quadric of
rank~$4$,
the linear systems of planes on which are
\mbox{$\bbQ$-rational}.
The~other four degenerate quadrics in the pencil are defined
over~$\smash{k' = \bbQ(\sqrt{5})}$
and form two length two orbits under the operation
of~$\smash{\Gal(\bbQ(\sqrt{5})/\bbQ)}$.
Thus,~having chosen the indices suitably, the absolute Galois group fixes
$[e_1]$
and~$[e'_1]$,
and has the two further orbits
$\smash{\{[e_2], [e'_2], [e_3], [e'_3]\}}$
and~$\smash{\{[e_4], [e'_4], [e_5], [e'_5]\}}$,
on which it acts in a totally coupled~manner. The~operation on either orbit is via the
full~$D_4$.

Consequently,~one has
$\Br_1(U)/\Br_0(U) \cong \bbZ/2\bbZ$,
generated by the
\mbox{$2$-torsion}
Brauer class
$\smash{\alpha = \iota_{\Pic} (\overline{[e_2] + [e_3]})}$
of type~II.
Indeed,~there is no proper
\mbox{$4$-torsion}
occurring, as the operation
of~$D_4$
as given above does not coincide with any of the subgroups
of~$(\bbZ/2\bbZ)^4 \!\rtimes\! S_5$,
described in Remarks~\ref{rem_typ1} and~\ref{rem_typ2}. Thus,~the claim directly follows from Theorem~\ref{expl_2}.

One~of the degenerate quadrics in the pencil defining
$X$~is
\begin{align*}
V\big((-\sqrt{5}\!-\!3)X_0^2 + (-2\sqrt{5}\!-\!6)X_0X_2 + (4\sqrt{5}\!+\!4)X_0X_4 + (2\sqrt{5}\!+\!6)X_1X_2 &\\[-1.5mm]
{} + (\sqrt{5}\!+\!3)X_1X_3 + (-4\sqrt{5}\!-\!8)X_1X_4 + (\sqrt{5}\!+\!3)X_2^2 + 4X_2X_3 &- 8X_4^2\big) \,.
\end{align*}
The~two linear systems of planes on this quadric are defined
over~$l'$
and conjugate to each other
under~$\smash{\Gal(l'/\bbQ(\sqrt{5}))}$.
The~tangent hyperplane at the point
$(0\!:\!0\!:\!0\!:\!1\!:\!0)$
is
$\smash{V\big((\sqrt{5}\!+\!3)X_1 + 4X_2\big)}$.
The~evaluation
of~$\alpha$
at a place
$\nu$
of~$\bbQ$
is hence given~by
\begin{equation}
\label{eval_ex518}
\sum_{w|\nu} \textstyle \big(\frac{(\sqrt{5}+3)X_1 + 4X_2}{X_0}, k'_w(\sqrt{2\sqrt{5}+4})/k'_w\big) ,
\end{equation}
where
$w$
runs over the places
of~$\smash{k' = \bbQ(\sqrt{5})}$
lying
above~$\nu$.

The~surface
$X$
has bad reduction only at the primes
$2$
and~$5$,
whereas the unique prime above
$5$
in~$\smash{\bbQ(\sqrt{5})}$
splits
in~$l'$.
On~the other hand, there is exactly one prime
of~$\smash{\bbQ(\sqrt{5})}$
lying
above~$2$,
too, and that actually ramifies
in~$l'$.
Finally,~there are two primes
of~$\smash{\bbQ(\sqrt{5})}$
lying above the infinite place, but one of them splits
in~$l'$.
Thus,~the sum (\ref{eval_ex518}) is effectively only one term in either case,
$\nu = 2$
or~$\infty$.
In~particular, the local evaluation map
$\ev_{\alpha,\infty}$
just tests the sign of the expression
$\smash{\frac{(-\sqrt{5}+3)X_1 + 4X_2}{X_0}}$.

An~experiment shows that
$\ev_{\alpha,2}$
takes both values on
$\calU(\bbZ_2)$
and that the same is true for
$\ev_{\alpha,\infty}$
on~$U(\bbQ)$.
Thus,~strong approximation is~violated. Both~combinations
$(0,0)$
and~$\smash{(\frac12, \frac12)}$
occur, in fact, for
\mbox{$\bbZ$-integral}~points.
\end{ex}

\begin{rem}[Some kind of normal form]
Let~$U$
be an open
degree~$4$
del Pezzo surface
over~$k$
carrying a
\mbox{$2$-torsion}
Brauer class of type~II.
Then~$U$
is as a hyperplane complement of a projective surface that may be defined in the~form
\begin{align*}
l_1l_2 - l_3^2 + dl_4^2 &= 0 \,, \\
l^\sigma_1 l^\sigma_2 - (l^\sigma_3)^2 + d^\sigma(l^\sigma_4)^2 &= 0 \,,
\end{align*}
for a quadratic extension
field~$k'$,
an element
$d \in k' \!\setminus\! (k')^2$,
and linearly independent linear forms
$l_1,\ldots,l_4$
over~$k'$.
Here,~$\sigma$
denotes the conjugation
of~$k'$.
The~evaluation is then given by
$\sum\limits_{w|\nu} (l_1(x), k'_w(\sqrt{d})/k'_w)$.
\end{rem}

\begin{rem}[The cyclic case]
There are clearly
\mbox{$2$-torsion}
classes of type~II that are~cyclic. Just~assume that the quotient
group~$G'$
induced
by~$G$
on the orbit
$\{[e_i], [e_j], [e'_i], [e'_j]\}$
is cyclic of
order~$4$.
It~is then generated by the
\mbox{$4$-cycle}~$\sigma$
depicted to the left in~Figure~\ref{fig_8}. In~this situation,
$l/k$
is a cyclic field extension of
degree~$4$
and
$k'$
is its intermediate quadratic~field.

Moreover,~the method for evaluation described in \ref{eval_cores} may then be directly translated into class field~theory. One~easily sees that the evaluation of the Brauer class is given by the norm residue~symbols
$$\textstyle (\frac1{t(x)t^\sigma(x)}, l_w/k_\nu) = (t(x)t^\sigma(x), l_w/k_\nu) \,,$$
for
$w$
any place
of~$l$
that lies
above~$\nu$.
Note~that, although
$l/k$
is of
degree~$4$,
the concrete norm residue symbols above take no values other than
$\sigma^2$
and the neutral~element.
Indeed,~$t(x)t^\sigma(x)$
is a norm from the intermediate
field~$k'$.
\end{rem}

\subsubsection*{$4$-torsion classes of type~I}\leavevmode

\begin{lem}
Let\/~$G \subseteq (\bbZ/2\bbZ)^4 \!\rtimes\! S_5$
be a subgroup of the kind that\/
$H^1(G, \frakP)$
contains a\/
$4$-torsion
class of type~I. Then,

\begin{iii}
\item
in the notation of~\ref{type1},
$G$
allows a natural surjection\/
$s\colon G \twoheadrightarrow \langle\tau\rangle \cong \bbZ/4\bbZ$
and\/
$\smash{\iota_{\Pic,4}(\overline{[e_1] + [e_2] + [e_3] + 2[e_4]}) \in H^1(G, \frakP)}$
is the image under inflation of a class\/
$\smash{\widetilde{c}\in H^1(\langle\tau\rangle, \frakP^{S_3})}$.
In~particular,
a\/~\mbox{$4$-torsion}
class of type~I is always~cyclic.

The~class\/
$\smash{\widetilde{c}}$
is that of the [unique]\/
\mbox{$1$-cocycle}
such that\/
$\smash{\tau \mapsto \frac{-[e_1]-[e_2]-[e_3]-[e_4]+[e_5]}2}$.
\item
The class\/
$\cl_\tau \cup\; \widetilde{c} \in H^{-1}(\langle\tau\rangle, \frakP^{S_3})$
is represented by\/
$\smash{\frac{-[e_1]-[e_2]-[e_3]-[e_4]+[e_5]}2}$.
\end{iii}\smallskip

\noindent
{\bf Proof.}
{\em
i)
One clearly has
$\ker s \subseteq S_3$,
hence
$\frakP^{S_3} \subseteq \frakP^{\ker s}$.
The~assertion itself follows from the description of
$\iota_{\Pic,4}$,
given in Lemma~\ref{4_tors}.a), together with an explicit calculation in~cocycles.\smallskip

\noindent
ii)~is, once again, a consequence of~i) and the explicit formula for the cup product, given in \cite[Chapitre~XI, Annexe, Lemme~3]{Se}.
}
\eop
\end{lem}

\begin{ttt}[Evaluation--$k$-rational quadrilaterals]
The conjugates of\vspace{.2mm}
$\smash{\frac{-[e_1]-[e_2]-[e_3]-[e_4]+[e_5]}2}$
under
$\smash{\langle\tau\rangle}$
are\!
$\smash{\frac{[e_1]+[e_2]+[e_3]-[e_4]-[e_5]}2}$,\!
$\smash{\frac{-[e_1]-[e_2]-[e_3]+[e_4]-[e_5]}2}$,\!
and\!
$\smash{\frac{[e_1]+[e_2]+[e_3]+[e_4]+[e_5]}2}$.
According~to Lemma~\ref{Pic_open}.iii), up to sign, these are the classes
in~$\Pic U_{\overline{k}}$
of four of the 16~lines. In~the blown-up model (cf.\ Paragraph~\ref{blowup}), we can say more precisely that the classes are
$-[E_5]$,
$-[L_{45}]$,
$-[E_4]$,
and~$-[C]$.

It~is not hard to see that these four lines form the edges of a~quadrilateral. I.e.,~in cyclic order, each line intersects its two neighbours, but not the~third. Such~a quadrilateral is certainly non-planar, as then the lines would pairwise meet each~other. On~the other hand, the four points of intersection determine a three-dimensional linear subspace
of~$\Pb^4_k$.
In~other words, there is a
\mbox{$k$-rational}
hyperplane
$V(t)$
cutting the quadrilateral out of the degree four del Pezzo
surface~$X$.

Finally,~let
$\smash{\widetilde{l}}$
be the field of definition of the 16~lines.
Then~$\smash{\Gal(\widetilde{l}/k) \cong G}$
and there is the intermediate field
$l$,
corresponding
to~$\ker s \subseteq G$.
Then~$\Gal(l/k) \cong \langle\tau\rangle$.
I.e.,~$l/k$
is a cyclic extension of
degree~$4$.
The~four edges of the quadrilateral are defined
over~$l$.
Moreover,~$\Gal(l/k)$
permutes them cyclically. We~therefore have
$\div t = N_{l/k}L$,
for~$L$
one of the edges of the quadrilateral, and
$\div 1/t = N_{l/k}(-L)$.
The evaluation of the Brauer class is hence given by the norm residue symbols
$$(1/t(x), l_w/k_\nu) = -(t(x), l_w/k_\nu) \,.$$
\end{ttt}

\begin{rem}
It~is not hard to determine the
\mbox{$k$-rational}
quadrilaterals algorithmically. The~reader might compare Algorithm~\ref{alg_Pic} below, which, however, does by far~more.
\end{rem}

\begin{ex}
Let~$X \subset \Pb^4_\bbQ$
be given by the system of~equations
\begin{align*}
-X_0^2 \!+\! 8X_0X_1 \!-\! 4X_0X_2 \!-\! 10X_0X_3 \!+\! 4X_1^2 \!-\! 6X_1X_2 \!-\! 8X_1X_3 \!+\! 2X_2^2 \!+\! 3X_2X_3 \!-\! X_3^2 & \\[-1.5mm]
{} \!-\! X_3X_4 &= 0 \,, \\
7X_0^2 \!+\! 9X_0X_1 \!+\! 6X_0X_3 \!+\! 7X_1^2 \!+\! X_2X_4 \!+\! 3X_3^2 &= 0 \,,
\end{align*}
$\calX \subset \Pb^4_\bbZ$
the subscheme that is defined by the same system of equations
as~$X$,
and put
$U := X \!\setminus\! H$
and
$\calU := \calX \!\setminus\! \calH$,
for
$H := V(X_0) \subset \Pb^4_\bbQ$
and~$\calH := V(X_0) \subset \Pb^4_\bbZ$.

A~point search provides almost 1000
\mbox{$\bbQ$-rational}
points of height up to
$1000$,
among which
$(1\!:\!0\!:\!2\!:\!-1\!:\!-2)$,
$(1\!:\!1\!:\!2\!:\!-1\!:\!-10)$,
$(1\!:\!-20\!:\!-32\!:\!-9\!:\!88)$,
and
$(1\!:\!-80\!:\!-62\!:\!11\!:\!718)$
are
\mbox{$\bbZ$-integral}.

The~16 lines
on~$X$
are defined over a number field of degree~twelve. The Galois group operating on the 10 linear systems of conics is exactly the dicyclic
group~$\Dic_3$,
described in Remark \ref{rem_typ1}.a.iii).
Thus,~one has
$\Br_1(U)/\Br_0(U) \cong \bbZ/4\bbZ$,
generated by a
\mbox{$4$-torsion}
class~$\alpha$
of type~I.
Note~here,
$p(\Dic_3) \cong S_3$
acts intransitively with two orbits of sizes
$3$
and~$2$.
Thus,~Theorem~\ref{expl_2} shows that no further
\mbox{$2$-torsion}
occurs.

There is exactly one
\mbox{$\bbQ$-rational}
hyperplane cutting a quadrilateral out
of~$X$,
namely
$V(35X_0 + 50X_1 - 19X_2 - 2X_3 + 5X_4)$.
The~four edges are defined
over~$\bbQ(\zeta_5)$
and acted upon cyclically by
$\Gal(\bbQ(\zeta_5)/\bbQ)$.
The~evaluation
of~$\alpha$
is hence given by the norm residue~symbols
$$\textstyle -\big(\frac{35X_0 + 50X_1 - 19X_2 - 2X_3 + 5X_4}{X_0}, \bbQ_\nu(\zeta_5)/\bbQ_\nu \big) .$$

The~surface
$X$
has bad reduction only at the primes
$2$,
$5$,
$31$,
and~$251$,
among which
$31$
and~$251$
completely split
in~$\bbQ(\zeta_5)$.
Moreover,~the local evaluation map
$\ev_{\alpha,2} |_{\calU(\bbZ_2)} \colon \calU(\bbZ_2) \to \bbQ/\bbZ$
is constant,~too.
Indeed,~$2$
is an inert prime. Thus,~otherwise there would exist a
\mbox{$\bbZ_2$-integral}
point
on~$\calU$
such that the linear form above takes an even~value. This,~however, yields that
$X_2$
and~$X_4$
must have different~parity. From~this, the second equation
defining~$X$
shows that
$X_3$
must be odd, which contradicts the first~equation.

Finally,~an experiment shows that
$\ev_{\alpha,5}$
takes all four values on
$\calU(\bbZ_5)$
and that
$\ev_{\alpha,\infty}$
takes the values
$0$
and~$\smash{\frac12}$
on~$U(\bbQ)$.
In~particular, on a
\mbox{$\bbZ$-integral}
point,
$\ev_{\alpha,5}$
may not take the values
$\smash{\frac14}$
or~$\smash{\frac34}$,
which alone shows that strong approximation is~violated.
\end{ex}

\begin{prop}[A normal form]
Let\/~$U$
be an open del Pezzo surface of degree four
over\/~$k$
that carries a\/
\mbox{$4$-torsion} 
Brauer class of type~I.
Then\/~$U$
is a hyperplane complement of a projective surface
$X$
that may be given by equations of the~form
\begin{align*}
l_1 l_2 &= m_0 z \,, \\
m_1 m_2 &= l_0 z \,,
\end{align*}
where the linear forms to the left are defined over a cyclic degree four extension\/~$l$
and acted upon by a generator
of\/~$\Gal(l/k)$
via the rule\/
$l_1 \mapsto m_1 \mapsto l_2 \mapsto m_2 \mapsto l_1$.
Moreover,
$l_0$~and\/~$m_0$
are linear forms defined over the intermediate quadratic field and conjugate to each other, and\/
$z$
is a linear form defined
over\/~$k$.\medskip

\noindent
{\bf Proof.}
{\em
We~take
$z$
to be a linear form that cuts out
of~$X$
a quadrilateral inducing the Brauer~class.
Next,~let
$L_1, \ldots, L_4$
be the four edges of the quadrilateral, in cyclic~order. Then
$L_1 \cup L_2$
determines a
plane~$V(z,l_1)$.
We~denote the Galois conjugates of the linear form
$l_1$,
in this order, by
$m_1$,
$l_2$,
and~$m_2$.
Then~$L_2 \cup L_3$
is contained in the plane
$V(z,m_1)$~etc.
Consequently,~the quadrilateral is given
in~$V(z) \cong \Pb^3$
by the system of equations
\begin{align*}
l_1 l_2 &= 0 \,, \\
m_1 m_2 &= 0 \,.
\end{align*}
Here,~the forms
$l_1 l_2$
and~$m_1 m_2$
are global sections of the coherent sheaf
$\calI_X(2) |_{V(z)}$.
In~order to complete the proof, it has to be shown that both can be lifted to sections of
$\calI_X(2)$
on the whole
of~$\Pb^4$.
As~there is the exact sequence
$$0 \longrightarrow \calI_X(1) \xlongrightarrow{\cdot z} \calI_X(2) \longrightarrow \calI_X(2) |_{V(z)} \longrightarrow 0 \,,$$
this indeed follows from the fact that
$H^1(\Pb^4, \calI_X(1)) = 0$~\cite[Remark~1.8.44]{La}.
}
\eop
\end{prop}

\begin{rem}
The surface above appears when blowing down the line
$V(l_0,m_0)$
on the cubic surface
$$C\colon l_0 l_1 l_2 = m_0 m_1 m_2$$
that is given in the classical Cayley--Salmon~form. The~nine obvious lines
on~$C$
are clearly defined over the
field~$l$.
The~combinatorics of the 27~lines on a cubic surface then implies that the other 18~lines may be defined at most over an
\mbox{$S_3$--extension}
of~$l$~\cite[Proposition~4.6]{EJ11}.
This~shows again that the field of definition of the 16~lines
on~$X$
has a Galois group that allows an injection
into~$S_3 \times \bbZ/4\bbZ$.
\end{rem}

\subsubsection*{$4$-torsion classes of type~II}\leavevmode

\begin{ttt}
When~$G \subseteq (\bbZ/2\bbZ)^4 \!\rtimes\! S_5$
is the cyclic group of
order~$4$,
generated by the element shown in Figure~\ref{fig_5}, then this is the overlap~case. One~has
$H^1(G, \frakP) \cong \bbZ/4\bbZ$,
the
\mbox{$4$-torsion}
class being~cyclic.
The~group
$G$
operates on the 16~lines such that there are four orbits of length four~each. Two~of these orbits form the edges of quadrilaterals, while the others consist of mutually skew~lines.
\end{ttt}

\begin{ttt}
When~$G \subseteq (\bbZ/2\bbZ)^4 \!\rtimes\! S_5$
is the cyclic group of
order~$8$,
generated by the element shown in Figure~\ref{fig_6}, then there is a
\mbox{$4$-torsion}
class of type~II
in~$H^1(G, \frakP)$,
and this, too, is clearly a cyclic~class.

Moreover,~one easily sees that the orbit of
$\smash{\frac{-[e_1]-[e_2]-[e_3]-[e_4]-[e_5]}2}$
consists of eight distinct elements
of~$\frakP$.
The~remaining eight elements
of~$\frakP$
that are of the type
$\smash{\frac{\pm[e_1]\pm[e_2]\pm[e_3]\pm[e_4]\pm[e_5]}2}$
and involve an even number of plus signs form a single orbit,~too, and either orbit has sum~zero.

Thus,~when
$U \subset X$
is an open del Pezzo surface of degree four
over~$k$
such that
$\Gal(\overline{k}/k)$
operates on the 16~lines exactly
via~$G$
then the lines form two orbits of length eight~each. As~one has
$\Pic U_{\overline{k}} = \Pic X_{\overline{k}}/\langle H\rangle$,
the sum over an orbit, considered
in~$\Pic X_{\overline{k}}$,
must be a multiple
of~$H$,
and since the degree
is~$8$,\pagebreak[3]
it
is~$2H$.
Moreover,~a short calculation in coherent cohomology shows that the restriction homomorphism
$\Gamma(\Pb^4, \calO_{\Pb^2}(2)) \to \Gamma(X, \calO_{\Pb^2}(2) |_X)$
is surjective, cf.~\cite[Exercise~III.5.5.a)]{Ha}. Hence, either of the two configurations of eight lines is cut out
of~$X$
by a quadric
in~$\Pb^4$.
\end{ttt}

\begin{ex}
Let~$X \subset \Pb^4_\bbQ$
be given by the system of~equations
\begin{align*}
X_0^2 - X_0X_3 - X_1^2 + X_1X_2 + X_1X_3 + X_2X_3 + X_3^2 - X_3X_4 &= 0 \,, \\
-X_0^2 - X_0X_1 - X_0X_3 + X_1^2 + X_2X_4 &= 0 \,,
\end{align*}
$\calX \subset \Pb^4_\bbZ$
the subscheme that is defined by the same system of equations
as~$X$,
and put
$U := X \!\setminus\! H$
and
$\calU := \calX \!\setminus\! \calH$,
for
$H := V(X_0) \subset \Pb^4_\bbQ$
and~$\calH := V(X_0) \subset \Pb^4_\bbZ$.

A~point search provides more than 6000
\mbox{$\bbQ$-rational}
points of height up to
$1000$,
among which
$(1\!:\!-1\!:\!0\!:\!1\!:\!-1)$,
$(1\!:\!1\!:\!0\!:\!-1\!:\!-1)$,
$(1\!:\!0\!:\!1\!:\!1\!:\!2)$,
$(1\!:\!2\!:\!0\!:\!1\!:\!-1)$,
$(1\!:\!0\!:\!-2\!:\!1\!:\!-1)$,
$(1\!:\!0\!:\!0\!:\!-1\!:\!-3)$,
$(1\!:\!0\!:\!3\!:\!-1\!:\!0)$,
$(1\!:\!42\!:\!-221\!:\!-47\!:\!8)$,
and
$(1\!:\!42\!:\!9\!:\!-277\!:\!-222)$
are
\mbox{$\bbZ$-integral}.

The~16 lines
on~$X$
are defined
over~$l := \bbQ(\zeta_{17}+\zeta_{17}^{-1})$.
The~Galois group
$\Gal(l/\bbQ)$
is cyclic of
order~$8$.
A~generator operates on the 10 linear systems of conics, when they are suitably indexed, exactly as the
element~$g$
shown in Figure~\ref{fig_6}. Thus,
$\Br_1(U)/\Br_0(U) \cong \bbZ/4\bbZ$,
generated by a
\mbox{$4$-torsion}
Brauer class
$\alpha$
of type~II.
Note~here that
$p(\langle g\rangle) \cong \bbZ/4\bbZ$
acts with two orbits of sizes
$4$
and~$1$,
so that Theorem~\ref{expl_2} excludes any further
\mbox{$2$-torsion}.

One finds that the quadrics
$V(q_1)$
and~$V(q_2)$,
for
\begin{align*}
q_1 &:= 3X_0X_2 \!-\! 2X_0X_3 \!+\! 12X_0X_4 \!-\! 4X_1^2 \!+\! 3X_1X_2 \!-\! 2X_1X_3 \!+\! 5X_1X_4 \!+\! X_2^2 \!+\! 17X_2X_3 \\[-1.5mm]
& \hspace{7.2cm} {}- 22X_2X_4 \!+\! 15X_3^2 \!-\! 13X_3X_4 \!-\! X_4^2 \,,\, \hspace{0.3cm}\text{and} \\
q_2 &:= 3X_0X_2 \!-\! 19X_0X_3 \!-\! 5X_0X_4 \!-\! 4X_1^2 \!+\! 20X_1X_2 \!+\! 15X_1X_3 \!+\! 5X_1X_4 \!+\! X_2^2 \!+\! 17X_2X_3 \\[-1.5mm]
& \hspace{8.5cm} {}+ 12X_2X_4 \!+\! 15X_3^2 \!-\! 13X_3X_4 \!-\! X_4^2
\end{align*}
cut the two configurations of eight lines out
of~$X$.
The~evaluation
of~$\alpha$
is hence given by the norm residue~symbols
$$- \big( q_1, \bbQ_\nu(\zeta_{17}+\zeta_{17}^{-1})/\bbQ_\nu \big) ,$$
and when one takes the quadratic
form~$q_2$,
the answer is the~same.

The~surface
$X$
has bad reduction only at the primes
$2$
and~$17$.
The~local evaluation map
$\ev_{\alpha,\infty}$
is constant, since the field
$\bbQ(\zeta_{17}+\zeta_{17}^{-1})$
is totally~real. Moreover,~the local evaluation map
$\ev_{\alpha,2} |_{\calU(\bbZ_2)} \colon \calU(\bbZ_2) \to \bbQ/\bbZ$
is constant,~too.
Indeed,~$2$
is an inert prime. Thus,~otherwise there would exist a
\mbox{$\bbZ_2$-integral}
point
on~$\calU$
such that both quadratic forms
$q_1$
and~$q_2$
evaluate to an even~number. This~leads to a contradiction as~follows. First,~the equations defining
$X$
imply that
$X_3$
must be~odd. Then~the first equation shows that
$X_1X_2 + X_2 + X_4$
must be odd,~too. Using~this, the assumptions about
$q_1$
and~$q_2$
reduce to
$X_1X_2 + X_1X_4 + X_2$
being odd and
$X_2 + X_4 + X_1 + X_1X_4$
being~even. These~conditions, however, turn out to be fulfilled only when all coordinates are odd, which contradicts the second equation
defining~$X$.

Finally,~an experiment shows that
$\ev_{\alpha,17}$
takes all four values on
$\calU(\bbZ_{17})$.
In~particular, on a
\mbox{$\bbZ$-integral}
point,
$\ev_{\alpha,17}$
may not take the values
$\smash{\frac14}$,
$\smash{\frac12}$,
or~$\smash{\frac34}$,
which shows that strong approximation is~violated.
\end{ex}

\begin{rem}
The~six further conjugacy classes of subgroups
of~$(\bbZ/2\bbZ)^4 \!\rtimes\! S_5$
that yield a
\mbox{$4$-torsion}
Brauer~class of type~II are non-cyclic. Even~worse, it turns out that no normal subgroup, except for the trivial group, annihilates the
\mbox{$4$-torsion}~class.
As~all other naive approaches, such as that to use a corestriction, do not apply either, it is our conclusion that only a generic algorithm~helps.
\end{rem}

\section{The generic algorithm}

\subsubsection*{The computer-algebraic framework}\leavevmode

\noindent
We implemented a generic algorithm in {\tt magma} relying on some of the data types and functionality already~existing. The data types~include

\begin{iii}
\item
Finite groups. I.e.,~permutation groups described by a sequence of~generators.
\item
Finitely generated abelian groups, described in the form
$\bbZ/a_1\bbZ \oplus \cdots \oplus \bbZ/a_m\bbZ$,
for a sequence of non-negative integers such that
$a_m|a_{m-1}| \ldots |a_2|a_1$.
\item
Finitely generated
\mbox{$G$-modules}.
This means that a finitely generated abelian group, as above, is given together with a sequence of matrices describing the operation of the generators
of~$G$.
The~data structure of a finitely generated
\mbox{$G$-module}
includes the full information about the underlying finite
group~$G$. 
\end{iii}

{\tt magma} allows to perform fundamental operations in these categories, including the computation of the cokernel
$\pr\colon N \to \coker \varphi$
of a homomorphism
$\varphi\colon M \to N$
of finitely generated abelian~groups.

Also,~for a
\mbox{$G$-module}
$M$
as above,
$H^i(G, M)$
may be computed for
$i = 0$,
$1$,
or~$2$.
This~means that an isomorphism is asserted to some
$H := \bbZ/a_1\bbZ \oplus \cdots \oplus \bbZ/a_m\bbZ$.
In~addition, for every
\mbox{$i$-cocycle}
with values
in~$M$,
its cohomology class may be expressed as an element
of~$H$.
And~vice versa, there is function returning for every element
of~$H$
a
representing~\mbox{$i$-cocycle}.

Furthermore,~let us mention that there is some more functionality that we take to work for granted, such as that to compute Gr\"obner bases, to compute the unit group of a finite ring, and to compute the Galois group of a polynomial over a number field, as well as the decomposition and inertia groups contained~within.

Finally,~an open del Pezzo surface of
degree~$4$
over a number
field~$k$
is, for us, just given by two quadratic forms in five variables over the maximal
order~$\calO_k$.
The~hyperplane section taken out is supposed to be the vanishing locus of the first~variable.

\subsubsection*{Computing the Picard group as a Galois module}

\begin{lem}
Let\/~$X$
be a proper del Pezzo surface of degree four over an algebraically closed field and\/
$U \subset X$
an open degree four del Pezzo surface (as in Definition~\ref{open_dP4}).

\begin{abc}
\item
Then\/~$X$
contains exactly 40~quadrilaterals.
\item
Let\/~$D$
be the free abelian group over the set of the 16~lines
on\/~$X$
and\/~$D_0 \subset D$
the submodule generated by all formal sums
$1L_1 + 1L_2 + 1L_3 + 1L_4$,
for\/
$L_1, \ldots, L_4$
the edges of a~quadrilateral. Then~there is a canonical isomorphism\/
$D/D_0 \cong \Pic U$.
\end{abc}\smallskip

\noindent
{\bf Proof.}
{\em
a)
We work in the blown-up~model. As~the lines
$E_i$
are mutually disjoint, every quadrilateral must contain a line of
type~$L_{ij}$.
Say,~a quadrilateral
contains~$L_{12}$.
The~lines intersecting
$L_{12}$
are
$E_1$,
$E_2$,
and~$L_{kl}$,
for~$3 \leq k < l \leq 5$.
Thus,~the are the~quadrilaterals

\begin{iii}
\item[$\bullet$ ]
$[E_1, L_{12}, E_2, C]$
and analogous, which~are
$(\atop{5}{2}) = 10$~items,
and
\item[$\bullet$ ]
$[E_1, L_{12}, L_{45}, L_{13}]$
and analogous. These~are
$5 \!\cdot\! (\atop{4}{2}) = 30$~items.
\end{iii}\smallskip

\noindent
b)
By~Lemma~\ref{Pic_open}, one knows that
$\Pic U \cong \bbZ^5$.
On~the other hand, the 16~lines generate
$\Pic U$.
I.e.,~the canonical homomorphism
$D \to \Pic U$
is~surjective.
Moreover,~$D_0$
is contained in the kernel, so that a homomorphism
$D/D_0 \twoheadrightarrow \Pic U$
gets~induced. It~is therefore sufficient to show that
$D/D_0$
may be generated by five~elements. Indeed,~every surjection from such a group
onto~$\bbZ^5$
is~bijective.

For~this, again let us work in the blown-up~model. First~of all, there are the relations
$L_{ij} \equiv -E_i - E_j - C \pmod {D_0}$,
which show that
$D/D_0$
is generated
by~$E_1, \ldots, E_5$,
and~$C$.
Furthermore, one~has
\begin{align*}
0 &\equiv  (E_1 + L_{12} + E_2 + C)
         + (E_1 + L_{13} + E_3 + C)
         + (E_4 + L_{45} + E_5 + C) \\[-1.5mm]
  & \hspace{8.0cm} - (E_1 + L_{12} + L_{45} + L_{13}) \pmod {D_0} \\
  &\equiv   E_1 + E_2 + E_3 + E_4 + E_5 + 3C \pmod {D_0} \,,
\end{align*}
implying that
$E_1, \ldots, E_4$,
and~$C$
form a generating system, as~required.
}
\eop
\end{lem}

\begin{alg}[Computing
$\Pic U_{\overline{k}}$
as a Galois module]
\label{alg_Pic}
Given an open
degree~$4$
del Pezzo surface
$U$
over a number field
$k$,
this algorithm computes the Picard group
$\Pic U_{\overline{k}}$
as a Galois~module.

\begin{iii}
\item
Represent lines in general position
in~$\Pb^4$
by parametrisations of the form
$(1 : t : (k_1\!+\!k_2t) : (k_3\!+\!k_4t) : (k_5\!+\!k_6t))$.
Write down a system of equations in six variables that encodes the containment of such a line in the
surface~$U$.
Then~calculate a Gr\"obner basis
over~$k$,
in order to obtain a univariate polynomial
$g$
of
degree~$16$,
the zeroes of which correspond one-to-one to the 16~lines
on~$U$.

If the Gr\"obner basis turns out not to be of the expected form then apply an automorphism
of~$\Pb^4$,
chosen at random, and try~again.
\item
Put~$l$
to be the splitting field
of~$g$
and calculate the Galois
group~$G$
of~$l$.
Make~sure that the operation
of~$G$,
i.e.\ the homomorphism
$i\colon G \to \Aut_k l$,
is stored,~too.
\item
Using Gr\"obner bases again, this time
over~$l$,
determine the 16~lines
on~$U$
explicitly. Store~them into a list, numbered from 1 to~16.

For~each generator
$\sigma \in G$,
use the automorphism
$i(\sigma)$
of~$l$
to determine the permutation of the 16~lines that
$\sigma$~induces.
Form~the associated permutation matrices, thereby transforming
$\bbZ^{16}$
into a
\mbox{$G$-module}
$\Div$.
\item
For~each pair of lines, determine using linear algebra
over~$l$
whether they intersect on the proper surface or~not. The result is a
$16\times16$-intersection
matrix.
\item
For each quadruple
$\{i_1, \ldots, i_4\} \subset \{1, \ldots, 16\}$,
determine whether it represents a quadrilateral. It~suffices to check that exactly four of the six intersection numbers are equal
to~$1$.
Store~the 40 quadrilaterals found into a
$16\times40$-matrix~$Q$
over~$\bbZ$,
each row of which contains four ones and zeros~otherwise.

Also~determine, using linear algebra
over~$l$,
for each quadrilateral a linear form cutting it out of the~surface. Store~these linear forms into a
list~$L$.
\item
Calculate
$\Pic U_{\overline{k}}$
as the cokernel of the homomorphism
$q\colon \bbZ^{40} \to \Div = \bbZ^{16}$
given
by~$Q$.
Determine the
\mbox{$G$-module}
structure
on~$\Pic U_{\overline{k}}$
that is induced from that
on~$\Div$
via the canonical epimorphism
$\pr\colon \Div \twoheadrightarrow \Pic U_{\overline{k}}$.
\item
Return~$\Pic U_{\overline{k}}$.
As~further values, return the homomorphism
$\pr\colon \Div \twoheadrightarrow \Pic U_{\overline{k}}$
[including the
\mbox{$G$-module}
$\Div$],
the homomorphism
$q\colon \bbZ^{40} \to \Div$,
the number
field~$l$,
and the
list~$L$.
\end{iii}
\end{alg}

\begin{rem}
Having~run Algorithm~\ref{alg_Pic}, it is just one call of a {\tt magma} intrinsic to compute
$H^1(G, \Pic U_{\overline{k}})$
in the form
$\bbZ/a_1\bbZ \;[\oplus \bbZ/a_2\bbZ]$.
Just~another call returns a
\mbox{$1$-cocycle}
for any [base] element~chosen.
\end{rem}

\subsubsection*{Computing a\/ $2$-cocycle with values in\/ $l(U)^*$}

\begin{alg}[Computing a $2$-cocycle with values in $l(U)^*/l^*$]
\label{alg_2coc}
Given a
\mbox{$1$-cocycle}
$\varphi\colon G \to \Pic U_{\overline{k}}$,
as well as the homomorphisms
$\pr\colon \Div \twoheadrightarrow \Pic U_{\overline{k}}$
and
$q\colon \bbZ^{40} \to \Div$
[together with the
list~$L$],
this algorithm computes a corresponding
\mbox{$2$-cocycle}
with values
in~$l(U)^*/l^*$.

\begin{iii}
\item
Lift~$\varphi$
along~$\pr$
to a cochain
$\widetilde\varphi\colon G \to \Div$.
\item
Calculate the coboundary
$\delta\widetilde\varphi\colon G \times G \to \Div$
according to the formula
$$\delta\widetilde\varphi(\sigma,\tau) := \varphi(\tau)^\sigma - \varphi(\sigma\tau) + \varphi(\sigma) \,.$$
\item
Lift~$\delta\widetilde\varphi$
along~$q$
to a function
$\psi\colon G \times G \to \bbZ^{40}$.
This~means to find a particular solution of a
$16\times40$
linear system of equations, for each pair
$(\sigma,\tau) \in G \times G$.
\item
Return the function
$\psi\colon G \times G \to \bbZ^{40}$.
\end{iii}
\end{alg}

\begin{rem}
A vector
$\psi(\sigma,\tau) = (e_1,\ldots,e_{40})$
encodes the rational function
$L_1^{e_1} \cdots L_{40}^{e_{40}}$
modulo~constants. For~performance reasons, we do not multiply the product out at this stage, but do so only after the evaluation at a~point.
\end{rem}

\begin{lem}[Lifting to a $2$-cocycle with values in $l(U)^*$]
\label{lift_coc}
Let\/~$U$
be an integral scheme over a
field\/~$k$
having a
\mbox{$k$-rational}
point\/
$x \in U(k)$,
$l/k$
a finite Galois extension,
$G = \Gal(l/k)$
its Galois group, and\/
$\smash{c\colon G \times G \to l(U)^*/l^*}$
a\/
\mbox{$2$-cocycle}.\smallskip

\noindent
Suppose~that all the images\/
$c(\sigma, \tau) \in l(U)^*/l^*$,
for\/~$\sigma, \tau \in G$,
are classes of functions modulo scalars that are defined and nonzero
at\/~$x$.
Then~the unique lift\/
$\smash{\widetilde{c}\colon G \times G \to l(U)^*}$
of\/~$c$
so that\/
$\widetilde{c}(\sigma, \tau)(x) = 1$,
for
all\/~$\sigma, \tau \in G$,
is a\/
\mbox{$2$-cocycle}.\medskip

\noindent
{\bf Proof.}
{\em
By assumption,
$\delta c = 0$.
Therefore,
$\delta\widetilde{c}\colon G \times G \times G \to l(U)^*$
vanishes
modulo~$l^*$.
In~other words,
$\delta\widetilde{c}$
takes only constant~functions. Thus,~in order to exactly determine
$\delta\widetilde{c}$,
it suffices to calculate the values of the
$\delta\widetilde{c}(\sigma, \tau, \upsilon)$
at the
point~$x$.
But,~according to the assumption, we have
$\widetilde{c}(\sigma, \tau)(x) = 1$
for
all\/~$\sigma, \tau \in G$.
Hence,
$\delta\widetilde{c}(\sigma, \tau, \upsilon)(x) = 1$,
for
all\/~$\sigma, \tau, \upsilon \in G$,
which completes the~proof.
}
\eop
\end{lem}

\begin{rem}
In~\cite[Section~3.3]{Be}, a different approach is presented on how to lift
\mbox{$2$-cocycles},
which applies to the situation of an affine Ch\^atelet~surface. 
\end{rem}

\begin{ttt}
\label{loc_cyc}
Given~a
\mbox{$2$-cocycle}
$c\colon \Gal(l/k) \times \Gal(l/k) \to l(U)^*$
representing the Brauer class
$\alpha \in \Br(U)$,
in order to evaluate
$\alpha$
at a point
$x \in U(k_\nu)$,
two things have to be done as the next~step.

\begin{abc}
\item
One has to restrict the
\mbox{$2$-cocycle}
to the local Galois group
$\Gal(lk_\nu/k_\nu) \subseteq \Gal(l/k)$.
\item
One~has to replace the functions
$c(\sigma,\tau)$
by their
values~$c(\sigma,\tau)(x) \in lk_\nu$.
\end{abc}
In~practice, one usually approximates the coordinates of a
\mbox{$k_\nu$-rational}
point by elements
of~$k$.
It~is possible to remain in the realm of number fields by doing what follows, instead of a) and~b). Fix~a
prime~$w$
of~$l$
lying
above~$\nu$~and

\begin{abc}
\item[a') ]
restrict the
\mbox{$2$-cocycle}
to the decomposition group
$D_w \subseteq \Gal(l/k)$.
\item[b') ]
Store the values
$c(\sigma,\tau)(x)$
as elements
of~$l$.
In~further calculations, they can be considered as elements
of~$l_w \;[\cong lk_\nu]$.
\end{abc}
The result is a
\mbox{$2$-cocycle}
representing a Brauer class of the local
field~$l_w$.
\end{ttt}

\subsubsection*{Computing the invariant of a Brauer class of a local field from a given\/
\mbox{$2$-cocycle}}

\begin{rem}
This~problem has been tackled before, in a very thorough manner, by T.~Preu~\cite{Pr}. We~are, however, not aware of any implementation realising Preu's approach. In~our implementation, which is sufficient for our purposes, we followed~\cite{Pr} only very partially, and took many shortcuts that make use of features built into~{\tt magma}.
Note~that the archimedean case is trivial, as then the local Galois group is of order at
most~$2$.
\end{rem}

\begin{lem}
\label{Gal_coh_liminv}
Let\/~$l/k$
be a finite Galois extension of non-archimedean local~fields.

\begin{abc}
\item
Then~there is a natural isomorphism
$$H^2(\Gal(l/k), l^*) \cong \varprojlim\limits_{\raisebox{2.7mm}{\scriptsize{n}}} H^2(\Gal(l/k), l^*\!/(1+\frakm_{\calO_l}^n)) \, . \vspace{-5mm}$$
\item
If~$\#\Gal(l/k)$
is relatively prime to the residue characteristic
of\/~$l$
then
$$\smash{H^2(\Gal(l/k), l^*) \cong H^2(\Gal(l/k), l^*\!/(1+\frakm_{\calO_l})) \, .}$$
\end{abc}\smallskip

\noindent
{\bf Proof.}
{\em
a)
Clearly,
$\smash{\calO_l^* \cong \varprojlim_n \calO_l^*\!/(1+\frakm_{\calO_l}^n)}$
and~$\smash{l^* \cong \varprojlim_n l^*\!/(1+\frakm_{\calO_l}^n)}$.
\mbox{Moreover,} for arbitrary natural
numbers
$i$
and~$n$,
the cohomology
$$\smash{H^i(\Gal(l/k), \calO_l^*\!/(1+\frakm_{\calO_l}^n))}$$
is a finite abelian~group.
Indeed,~$\Gal(l/k)$
is a finite group and
$\smash{\calO_l^*\!/(1+\frakm_{\calO_l}^n)}$
a finite
$\Gal(l/k)$-module.
Consequently,~\cite[Corollary~2.7.6]{NSW} applies and shows~that
$$H^i(\Gal(l/k), \calO_l^*) \cong \varprojlim\limits_{\raisebox{2.7mm}{\scriptsize{$n$}}} H^i(\Gal(l/k), \calO_l^*\!/(1+\frakm_{\calO_l}^n)) \,, \vspace{-5mm}$$
for
every~$i$.
Observe~here that, in our case, continuous group cohomology agrees with ordinary group cohomology when we equip the finite group
$\Gal(l/k)$
with the discrete~topology.

Furthermore,~for
each~$n$,
the commutative diagram
$$
\xymatrix@R=10pt{
0 \ar@{->}[r]& \calO_l^* \ar@{->}[r]\ar@{->}[d]& l^* \ar@{->}[r]^{\nu_l}\ar@{->}[d]& \bbZ \ar@{->}[r]\ar@{=}[d]& 0 \\
0 \ar@{->}[r]& \calO_l^*\!/(1+\frakm_{\calO_l}^n) \ar@{->}[r]& l^*\!/(1+\frakm_{\calO_l}^n) \ar@{->}[r]^{\;\;\;\;\;\;\;\;\;\;\;\overline\nu_l}& \bbZ \ar@{->}[r]& 0
}
$$
of short exact sequences induces a commutative diagram
{\tiny
\begin{equation}
\label{diag}
\xymatrix@R=10pt@C=.85em{
0 \ar@{->}[r]& H^2(\Gal(l/k), \calO_l^*) \ar@{->}[r]\ar@{->}[d]& H^2(\Gal(l/k), l^*) \ar@{->}[r]\ar@{->}[d]& H^2(\Gal(l/k), \bbZ) \ar@{->}[r]\ar@{=}[d]& H^3(\Gal(l/k), \calO_l^*)\ar@{->}[d] \\
0 \ar@{->}[r]& H^2(\Gal(l/k), \calO_l^*\!/(1\!+\!\frakm_{\calO_l}^n)) \ar@{->}[r]& H^2(\Gal(l/k), l^*\!/(1\!+\!\frakm_{\calO_l}^n)) \ar@{->}[r]^{\;\;\;\;\;\;\;\;\;\;\overline{\overline\nu}_l}& H^2(\Gal(l/k), \bbZ) \ar@{->}[r]& H^3(\Gal(l/k), \calO_l^*\!/(1\!+\!\frakm_{\calO_l}^n))
}
\end{equation}}
of exact~sequences.
Indeed,~$H^1(\Gal(l/k), \bbZ) = 0$.

Putting~$C_n := \im \overline{\overline\nu}_l$,
the~lower exact sequence may be split into two,
$$0 \to H^2(\Gal(l/k), \calO_l^*\!/(1\!+\!\frakm_{\calO_l}^n)) \to H^2(\Gal(l/k), l^*\!/(1\!+\!\frakm_{\calO_l}^n)) \to C_n \to 0$$
and
$$0 \to C_n \to H^2(\Gal(l/k), \bbZ) \to H^3(\Gal(l/k), \calO_l^*\!/(1\!+\!\frakm_{\calO_l}^n)) \,.$$
Either~of these sequences remains exact after applying the inverse image
functor~$\smash{\varprojlim_n}$.
Indeed,~this is clear for the second one as
$\smash{\varprojlim_n}$
is left-exact.
Moreover,~since all the abelian groups
$H^2(\Gal(l/k), \calO_l^*\!/(1\!+\!\frakm_{\calO_l}^n))$
are finite, they form an inverse system satisfying the Mittag--Leffler condition. The~claim on the first exact sequence therefore results from~\cite[Proposition~2.7.3]{NSW}.

Consequently,~the application of
$\smash{\varprojlim_n}$
to~(\ref{diag}) yields a commutative diagram of exact sequences, in which all vertical homomorphisms, except possibly for the second from the left, are isomorphisms. A~standard diagram argument shows that the latter must be an isomorphism, too, as~required.\smallskip

\noindent
b)
In~this situation, for
each~$n \geq 1$,
the quotient
$(1+\frakm_{\calO_l}^n)/(1+\frakm_{\calO_l}^{n+1})$
is an abelian
\mbox{$p$-group},
for~$p$
the residue characteristic
of~$l$,
and hence cohomologically~trivial. From~this, the long exact cohomology sequence associated~with
$$0 \longrightarrow (1+\frakm_{\calO_l}^n)/(1+\frakm_{\calO_l}^{n+1}) \longrightarrow l^*\!/(1+\frakm_{\calO_l}^{n+1}) \longrightarrow l^*\!/(1+\frakm_{\calO_l}^n) \longrightarrow 0$$
immediately implies that the sequence
$H^2(\Gal(l/k), l^*\!/(1+\frakm_{\calO_l}^n))$
is stationary for
$n \geq 1$.
}
\eop
\end{lem}

\begin{rems}
\begin{iii}
\item
Part~b) may also be shown directly, without relying on~a), by applying \cite[Chapitre~XII, \S3, Lemme~3]{Se} 
for~$q := 2$
and~$3$
to the
\mbox{$\Gal(l/k)$-module}
$M := 1+\frakm_{\calO_l}$
and its submodules
$M_n := 1+\frakm_{\calO_l}^n$,
for~$n \geq 1$.
\item
In~the application to open del Pezzo surfaces of degree four carrying an algebraic
\mbox{$4$-torsion}
Brauer class of type~II,
$\Gal(l/k)$
is always
a~\mbox{$2$-group}.
Hence,~Lemma~\ref{Gal_coh_liminv}.b) is applicable for the evaluation at all non-archimedean primes, except for those of residue
characteristic~$2$.

For~simplicity of the presentation, we restrict ourselves to the case that
$\#G$
is a power
of~$2$,
although our approach works in~general.
\end{iii}
\end{rems}

\begin{alg}[for the case that
$\#G$
is a power
of~$2$,
the residue characteristic
is~$\neq \!2$,
and the Brauer class is at most
\mbox{$4$-torsion}]%
\label{eval_loc}\leavevmode

\noindent
Given~a
\mbox{$2$-cocycle}
$\varphi\colon G \times G \to l^*$,
for
$l/k$
a finite Galois extension of non-archimedean local fields and
$G \cong \Gal(l/k)$,
satisfying the conditions listed, this algorithm computes the invariant
in~$\bbQ/\bbZ$
of the class in
$H^2(G, l^*)$
represented
by~$\varphi$.

\begin{iii}
\item
Put~$q := \#(\calO_l/\frakm_{\calO_l})$
and~$\smash{f := \log q/\log \#(\calO_k/\frakm)}$.
If~$f=1$
or~$2$
then execute the modi\-fication described below in~\ref{mod_inertia}.c).
Calculate the inertia group
$E \subseteq G$
and
put~$e := \#E$.

Run in a loop
through~$G$
and identify a Frobenius element
$\Frob \in G$
by testing the condition that
$\nu_l(\Frob(x) - x^{\#(\calO_k/\frakm_{\calO_k})}) > 0$,
for
$x \in \calO_l$
an element such that
$\overline{x} \in \calO_l/\frakm_{\calO_l}$
generates that field as an extension
of~$\calO_k/\frakm_{\calO_k}$.

Run, once again, in a loop
through~$G$
and create a map
$\expo\colon G \to \{0,1,\ldots,f-1\}$
of the kind that
$\smash{\Frob^{\expo(g)} \equiv g \pmod E}$,
for every
$g \in G$.
\item
Fix~a uniformiser
$\pi_l$
of~$l$
and a multiplicative
generator~$u \in \calO_l/\frakm_{\calO_l}$
of the residue~field. 
Then~create the
\mbox{$G$-module}
$l^*\!/(1+\frakm_{\calO_l})$
as~follows. Put, at~first
$M := \bbZ \oplus \bbZ/(q-1)\bbZ$.
In~addition, associate the~matrices
$$A_i :=
\left(
\begin{array}{cc}
1 & 0 \\
\smash{\log_u \!\frac{g_i(\pi_l)}{\pi_l}} & \#(\calO_k/\frakm)^{\expo(g_i)}
\end{array}
\right) ,
$$
with the generators
$g_i \in G$,
in order to transform
$M$~into
a
\mbox{$G$-module}.
Here,
$\log_u\colon (\calO_l/\frakm_{\calO_l})^* \to \bbZ/(q-1)\bbZ$
denotes the discrete logarithm of
base~$u$.
\item
Define the standard
\mbox{$2$-cocycle}
$\st\colon G \times G \to M$
of invariant
$(1/f \bmod 1)$~by
$$(\sigma,\tau) \mapsto \left\{
\begin{array}{cl}
(0,0)   & \text{ if } \expo(\sigma) + \expo(\tau) < f \,, \\
(e, \log_u \!\frac{\pi}{\pi_l^e}) & \text{ if } \expo(\sigma) + \expo(\tau) \geq f \,,
\end{array}
\right.
$$
for~$\pi$
a uniformiser
of~$k$.
\item
Calculate the image of the input cocycle by~putting
$$\textstyle \widetilde\varphi(\sigma, \tau) := \big( \nu_l(\varphi(\sigma, \tau)), \log_u \!\frac{\varphi(\sigma, \tau)}{\pi_l^{\nu_l(\varphi(\sigma, \tau))}}\big) ,$$
for every
$(\sigma, \tau) \in G \times G$.
\item
Calculate the abelian group
$H^2(G, M)$
and, as its elements, the classes
$c_\st$
and~$c_{\widetilde\varphi}$
of the
\mbox{$2$-cocycles}
$\st$
and~$\widetilde\varphi$.
The~cohomology class
$c_\st$~is
a proper
\mbox{$f$-torsion}
element and
$c_{\widetilde\varphi}$
is one of its~multiples. Find~a natural number
$m$
such that
$c_{\widetilde\varphi} = m c_\st$,
and output
$(m/f \bmod 1) \in \bbQ/\bbZ$
as the desired~invariant.
\end{iii}
\end{alg}

\begin{ttt}
\label{mod_inertia}
\begin{abc}
\item
The~idea behind Algorithm~\ref{eval_loc} is to compare the given
\mbox{$2$-cocycle},
which is for the extension
field~$l$,
with the standard
\mbox{$2$-cocycle}
for the maximal unramified subextension
$l \cap k^\nr$
of
degree~$f$.
Thus,~it may work only for invariants being an integral multiple
of~$(1/f \bmod 1)$.
This~is why in step~i) we had to suppose
that~$f \neq 1,2$.
\item
If~$f = 1$
or~$2$
then one might replace
$l$
by~$lk^{\nr,4}$,
for~$k^{\nr,4}$
the unramified
\mbox{degree-$4$}
extension
of~$k$,
extend the
group~$G$
accordingly, and then run Algorithm~\ref{eval_loc} as~described.
\item
However,~the standard
\mbox{$2$-cocycle}
has values only in
$k \subseteq l$
and the 
\mbox{$2$-cocycle}
to be evaluated has values only
in~$l$.
Moreover,~$\pi_l$
is still a uniformiser
for~$lk^{\nr,4}$.
Thus,~we can get by without explicitly introducing the larger
field~$lk^{\nr,4}$.
Instead,

\begin{iii}
\item[i') ]
put~$G' := G \times \bbZ/4\bbZ$
if
$f=1$
and,
if~$f=2$,
$$G' := \ker(G \times \bbZ/4\bbZ \xlongrightarrow{(\expo \!\bmod 2) \circ \pr_1 - \;(\pr_2 \!\bmod 2)} \bbZ/2\bbZ) \,.$$
Moreover,~let
$\varphi'\colon G' \times G' \to l^*$
be the composition
of~$\pr_1 \times \pr_1$
with~$\varphi$.
Define
$\expo'\colon G' \twoheadrightarrow \{0,\ldots,3\}$
to be the lift of the homomorphism
$G \to \bbZ/4\bbZ$
induced by the projection to the second~factor.
Put~$q' := q^{4/f}$
and~$f' := 4$.
\item[ii') ]
Put~$\smash{\log'_u := \frac{q'-1}{q-1} \!\cdot\! \log_u}$,
at least
on~$(\calO_l/\frakm_{\calO_l})^*$.
Then~put
$M' := \bbZ \oplus \bbZ/(q'-1)\bbZ$
and transform
$M'$
into a
\mbox{$G'$-module}
by means of
matrices~$A'_i$
that are defined as in~iv), but using
$\log'_u$
and~$\expo'$
instead of
$\log_u$
and~$\expo$.

Finally,~execute steps iii), iv), and~v),~accordingly.
\end{iii}
\end{abc}
\end{ttt}

To~summarise, the generic algorithm runs as~follows.

\begin{alg}[The generic algorithm]
\label{alg_gen}
Given an open
degree~$4$
del Pezzo
surface~$U$
over a number
field~$k$
and a sequence of
\mbox{$\calO_{k,\nu}$-integral}
points, for
$\nu$
a prime
of~$k$,
this algorithm computes the evaluations of the points given, for any algebraic Brauer class
on~$U$.

\begin{iii}
\item
Run~Algorithm~\ref{alg_Pic} to compute the Galois module
$\Pic U_{\overline{k}}$.
Extract~the underlying
group~$G$,
which is the Galois group operating on the 16~lines
on~$U$.
Moreover,~store the
field~$l$,
over which the lines are~defined.
\item
Call the {\tt magma} intrinsic to compute
$H^1(G, \Pic U_{\overline{k}})$.
Choose~the element to be~evaluated.
\item
Run~Algorithm~\ref{alg_2coc} to compute a
\mbox{$2$-cocycle}
$c$
with values in
$l(U)^*/l^*$
representing the cohomology class~chosen.
Lift~$c$
to a
\mbox{$2$-cocycle}
$\widetilde{c}$
with values
in~$l(U)^*$
as described in Lemma~\ref{lift_coc}, using the first of the given~points.
\item
Compute~the decomposition group
$D_w \subseteq G$,
for~$w$
any prime
of~$l$
lying
above~$\nu$.
Then~localise the
\mbox{$2$-cocycle}
$\widetilde{c}$
at the
prime~$\nu$
as described in~\ref{loc_cyc}.a'). The~result is a
\mbox{$2$-cocycle}
$c_\nu$
for the operation of the decomposition
group~$D_w \subseteq G$.
\item
Run~in a loop over the given sequence of
\mbox{$\calO_{k,\nu}$-integral}
points. For~each of them do the~following.
\item[$\bullet$ ]
Evaluate~the
\mbox{$2$-cocycle}
$c_\nu$
at the current
point~$\xi$,
as described in~\ref{loc_cyc}.b'). This~yields a
\mbox{$2$-cocycle}
$c_{\nu,\xi}$
representing a cohomology class
in~$H^2(D_w, l_w^*)$.
\item[$\bullet$ ]
Run~Algorithm~\ref{eval_loc} to determine the invariant of the
\mbox{$2$-cocycle}
$c_{\nu,\xi}$.
Store~the value into a~list.
\item
Output~the list of invariants~found.
\end{iii}
\end{alg}

\begin{rem}
We~assume, in practice, that the
\mbox{$\calO_{k,\nu}$-integral}
points are given as
\mbox{$k$-rational}
points lying
on~$U$.
One~might have the idea to work with
\mbox{$k$-rational}
approximations instead that are not exactly located
on~$U$.
In~such a case, one would have to take care, in addition, about the quality of the approximation.
\end{rem}

\begin{ex}[requiring the generic algorithm]
\label{fl176}
Let~$X \subset \Pb^4_\bbQ$
be given by the system of~equations
\begin{align*}
X_0^2 + 2X_0X_1 - 3X_0X_3 + X_1^2 - 3X_1X_3 - X_2^2 - X_2X_4 + 2X_3^2 - X_3X_4 &= 0 \,, \\
-2X_0^2 - X_0X_1 + 2X_0X_3 - 2X_1^2 + 2X_1X_3 + X_2X_4 - X_3^2 &= 0 \,,
\end{align*}
$\calX \subset \Pb^4_\bbZ$
the subscheme that is defined by the same system of equations
as~$X$,
and put
$U := X \!\setminus\! H$
and
$\calU := \calX \!\setminus\! \calH$,
for
$H := V(X_0) \subset \Pb^4_\bbQ$
and~$\calH := V(X_0) \subset \Pb^4_\bbZ$.

A~point search provides more than 600
\mbox{$\bbQ$-rational}
points of height up to
$1000$,
among which
$(1\!:\!1\!:\!-1\!:\!2\!:\!-1)$,
$(1\!:\!15\!:\!-5\!:\!4\!:\!-71)$,
$(1\!:\!20\!:\!15\!:\!-6\!:\!74)$,
$(1\!:\!-9\!:\!-2\!:\!1\!:\!-86)$,
$(1\!:\!41\!:\!15\!:\!-6\!:\!263)$,
$(1\!:\!223\!:\!-229\!:\!308\!:\!-247)$,
$(1\!:\!299\!:\!-213\!:\!312\!:\!-419)$,
and
$(1\!:\!-96\!:\!-53\!:\!22\!:\!-434)$
are
\mbox{$\bbZ$-integral}.

The~16 lines
on~$X$
are defined over a number
field~$l$
of
degree~$64$.
The~Galois group
$\Gal(l/\bbQ)$
operates on the 10 linear systems of conics as the largest
subgroup~$G$
of~$(\bbZ/2\bbZ)^4 \!\rtimes\! S_5$
fulfilling conditions~\ref{type2}.i) and~ii).
Thus,~$\Br_1(U)/\Br_0(U) \cong \bbZ/4\bbZ$,
generated by a
\mbox{$4$-torsion}
Brauer class
$\alpha$
of type~II.
Note~here that
$p(G) \cong D_4$
acts intransitively with two orbits of sizes
$4$
and~$1$,
so that Theorem~\ref{expl_2} excludes any further
\mbox{$2$-torsion}.

First~of all, it turns out that the local evaluation map
$\ev_{\alpha,\infty}\colon U(\bbR) \to \frac12\bbZ/\bbZ$
at the infinite place is constantly~zero. Indeed,~the base extension
of~$\alpha$
is an element of
$\Br_0(U_\bbR)$,
although
$\Br_1(U_\bbR)/\Br_0(U_\bbR) \cong \bbZ/2\bbZ$.
In~order to explain this, let us note that exactly three of the five degenerate quadrics in the pencil
defining~$X$
are~real, as is shown by a simple calculation. In~particular, the decomposition group
$D_{w_\infty}$
is indeed of
order~$2$.
Among~the degenerate quadrics, only one has its linear systems of planes defined
over~$\bbR$.
Numbering~the five degenerate quadrics from
$1$
to~$5$
as in~\ref{type2},
$\smash{D_{w_\infty}}$
must split the
orbit~$\{1\}$.
This~is implied by Theorem~\ref{4tors_rest}.ii), since an order two group never gives rise to a
\mbox{$4$-torsion}
class. More~precisely, the operation of the nonzero element
of~$\smash{D_{w_\infty}}$
is necessarily of the form indicated in~Figure~\ref{fig_11} below, which implies the claim, again according to Theorem~\ref{4tors_rest}.ii).\medskip

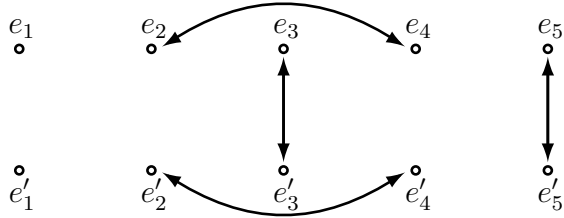
\begin{figure}[H]
\begin{center}	
\begin{picture}(200,60)
\linethickness{1pt}
\setlength\unitlength{2pt}
\put(0,3){\circle{1.5}}
\put(0,26){\circle{1.5}}
\put(25,3){\circle{1.5}}
\put(25,26){\circle{1.5}}
\put(50,3){\circle{1.5}}
\put(50,26){\circle{1.5}}
\put(75,3){\circle{1.5}}
\put(75,26){\circle{1.5}}
\put(100,3){\circle{1.5}}
\put(100,26){\circle{1.5}}

\put(-2,29){$e_1$}
\put(23,29){$e_2$}
\put(48,29){$e_3$}\put(50,22){\vector(0,-1){17.5}}
\put(73,29){$e_4$}\cbezier(28,27)(43,37)(57,37)(72,27)\put(28,27){\vector(-10,-6){1}}\put(72,27){\vector(10,-6){1}}
\put(98,29){$e_5$}\put(100,22){\vector(0,-1){17.5}}
\put(-2,-3){$e'_1$}
\put(23,-3){$e'_2$}
\put(48,-3){$e'_3$}\put(50,7){\vector(0,1){17.5}}
\put(73,-3){$e'_4$}\cbezier(28,2)(43,-8)(57,-8)(72,2)\put(28,2){\vector(-10,6){1}}\put(72,2){\vector(10,6){1}}
\put(98,-3){$e'_5$}\put(100,7){\vector(0,1){17.5}}
\end{picture}
\end{center}\medskip

\caption{The operation of the decomposition group at infinity}
\label{fig_11}
\end{figure}\vskip-\medskipamount

The~surface
$X$
has bad reduction at the primes
$3$,
$5$,
$7$,
and~$19$.
The~local evaluation map
$\ev_{\alpha,19} |_{\calU(\bbZ_{19})} \colon \calU(\bbZ_{19}) \to \bbQ/\bbZ$
is~constant. Indeed,~direct calculations show the following. Three~of the five degenerate quadrics in the pencil
defining~$X$
are defined
over~$\smash{\bbQ_{19}}$,
the other two only over the ramified extension
$\smash{\bbQ_{19}(\sqrt{-19})}$.
The~linear systems of planes on the latter two are defined over
$\smash{\bbQ_{19}(\sqrt{-19})}$,~already.
Finally,~two of the former three have their linear systems of planes defined only over the unramified extension
field~$\smash{\bbQ_{19}(\sqrt{-1})}$,
while on the final one, they are defined
over~$\bbQ_{19}$.
Thus,~the decomposition group is isomorphic to the Klein four group, and, again following the conventions of~\ref{type2}, two of its generators necessarily operate as shown below, which implies the claim according to Theorem~\ref{4tors_rest}.ii).

\begin{figure}[H]
\begin{center}	
\begin{picture}(165,55)
\linethickness{.8pt}
\setlength\unitlength{1.7pt}
\put(0,3){\circle{1.5}}
\put(0,26){\circle{1.5}}
\put(25,3){\circle{1.5}}
\put(25,26){\circle{1.5}}
\put(50,3){\circle{1.5}}
\put(50,26){\circle{1.5}}
\put(75,3){\circle{1.5}}
\put(75,26){\circle{1.5}}
\put(100,3){\circle{1.5}}
\put(100,26){\circle{1.5}}
\linethickness{1pt}
\put(-2,29){$e_1$}
\put(22.5,29){$e_2$}
\put(48,29){$e_3$}
\put(73,29){$e_4$}\cbezier(28,27)(43,38)(57,38)(72,27)\put(28,27){\vector(-10,-6){1}}\put(72,27){\vector(10,-6){1}}
\put(98,29){$e_5$}
\put(-2,-3){$e'_1$}
\put(22.5,-3){$e'_2$}
\put(48,-3){$e'_3$}
\put(73,-3){$e'_4$}\cbezier(28,2)(43,-9)(57,-9)(72,2)\put(28,2){\vector(-10,6){1}}\put(72,2){\vector(10,6){1}}
\put(98,-3){$e'_5$}
\end{picture}\hspace{2.6cm}
\begin{picture}(165,55)
\linethickness{0.8pt}
\setlength\unitlength{1.7pt}
\put(0,3){\circle{1.5}}
\put(0,26){\circle{1.5}}
\put(25,3){\circle{1.5}}
\put(25,26){\circle{1.5}}
\put(50,3){\circle{1.5}}
\put(50,26){\circle{1.5}}
\put(75,3){\circle{1.5}}
\put(75,26){\circle{1.5}}
\put(100,3){\circle{1.5}}
\put(100,26){\circle{1.5}}
\linethickness{1pt}
\put(-2,29){$e_1$}
\put(23,29){$e_2$}
\put(48,29){$e_3$}\put(50,22){\vector(0,-1){17.5}}
\put(73,29){$e_4$}
\put(98,29){$e_5$}\put(100,22){\vector(0,-1){17.5}}
\put(-2,-3){$e'_1$}
\put(23,-3){$e'_2$}
\put(48,-3){$e'_3$}\put(50,7){\vector(0,1){17.5}}
\put(73,-3){$e'_4$}
\put(98,-3){$e'_5$}\put(100,7){\vector(0,1){17.5}}
\end{picture}
\end{center}\smallskip

\caption{Generators of the decomposition group at $p=19$}\vskip-\bigskipamount
\end{figure}
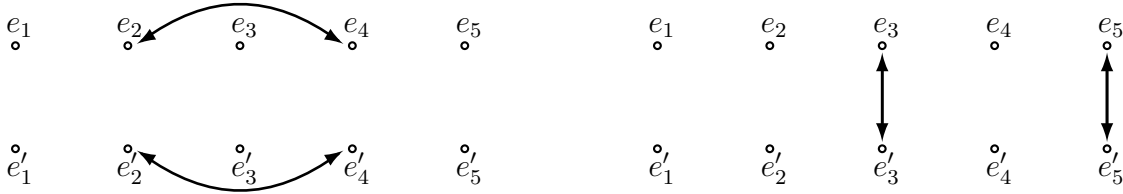

The~Brauer
class~$2\alpha$
is a
\mbox{$2$-torsion}
class of type~I. Indeed,~using the conventions as above, Theorem~\ref{expl_2} shows that
$\smash{2\alpha = \iota_{\Pic}(\overline{[e_1]})}$.
Moreover,~the evaluation
of~$2\alpha$
at a
prime~$p$
is given by the norm residue~symbol
$\smash{\big(\frac{X_2+X_3}{X_0}, \bbQ_p(\sqrt{5})\big)}$.
This~is in fact obtained by a direct application of the method, described in~\ref{eval_tang}. Note~here that
$V(q_1)$,
for
$q_1$
the term to the left in the first equation
of~$X$,
is the only quadric in the pencil
defining~$X$
that is defined
over~$\bbQ$
and has
rank~$4$.

As~a quadratic form of
rank~$4$,
$q_1$~is
of discriminant
exactly~$5$.
Thus,~Corollary~\ref{quad_const} shows that the local evaluation~map
$\ev_{2\alpha,p} |_{\calU(\bbZ_p)}$
is constant, for
$p = 3$
and~$7$.
Note~here that
$(1\!:\!(-1)\!:\!0\!:\!0\!:\!0)$
is the cusp
of~$V(q_1)$.
Its~reduction
modulo~$7$
does not lie
on~$V(q_2)_7$.
On~the other hand, its reduction
modulo~$3$
lies
on~$V(q_2)_3$,
but as a singular point, so that no
\mbox{$\bbQ_3$-rational}
point
on~$V(q_2)$
can reduce to~it.
Consequently,~$\ev_{\alpha,p} |_{\calU(\bbZ_p)}$
may take at most two values for either of these~primes.

We~ran Algorithm~\ref{alg_gen} for each of the four
primes~$p = 3$,
$5$,
$7$,
and~$19$.
As~our sample of
\mbox{$\bbZ_p$-integral}
points, we used the set of all of
\mbox{$\bbQ$-rational}
points, found by the point search, that are
\mbox{$p$-adically}~integral.
The~outcome is as~follows.

\begin{iii}
\item
Over~the prime
$p=3$,
the decomposition group is cyclic of
order~$8$
and the inertia group is its subgroup of
order~$2$.
In~particular, the modification of Algorithm~\ref{eval_loc} does not get~called. Among~the
\mbox{$\bbQ$-rational}
points found, 510 are 
\mbox{$3$-adically}~integral.
The~lift
$\alpha \in \Br_1(U)$
of one of the two proper
\mbox{$4$-torsion}
classes in
$\smash{H^1(\Gal(\overline\bbQ/\bbQ), \Pic U_{\overline\bbQ})}$,
chosen running the algorithm, evaluates 208~times
to~$0$
and 202~times
to~$\smash{\frac12}$.
\item
Over~the prime
$p=5$,
the decomposition group is cyclic of
order~$4$
and coincides with the inertia~group. I.e., the decomposition field is totally~ramified. The~modification of Algorithm~\ref{eval_loc} runs and enlarges the decomposition group
to~$\bbZ/4\bbZ \times \bbZ/4\bbZ$.
Among~the
\mbox{$\bbQ$-rational}
points found, 452 are 
\mbox{$5$-adically}~integral.
The~Brauer class 
$\alpha \in \Br_1(U)$
evaluates
163~times
to~$0$,
99~times
to~$\smash{\frac14}$,
79~times
to~$\smash{\frac12}$,
and 111~times
to~$\smash{\frac34}$.
\item
Over~the prime
$p=7$,
the decomposition group is isomorphic
to~$\bbZ/4\bbZ \times \bbZ/2\bbZ$
and the inertia group is a subgroup of
order~$2$
with a cyclic~quotient. The~modification of Algorithm~\ref{eval_loc} is not~called. Among~the
\mbox{$\bbQ$-rational}
points found, 576 are 
\mbox{$7$-adically}~integral.
The~Brauer
class~$\alpha$
evaluates
341~times
to~$0$
and 235~times
to~$\smash{\frac12}$.
\item
Over~the prime
$p=19$,
the algorithm detects that the decomposition group is isomorphic to the Klein four group and that the inertia group is a subgroup of
order~$2$.
The~modification runs and enlarges the decomposition group
to~$\bbZ/4\bbZ \times \bbZ/2\bbZ$.
Among~the
\mbox{$\bbQ$-rational}
points found, 582 are 
\mbox{$19$-adically}~integral
and the Brauer
class~$\alpha$
is found to evaluate
to~$0$
at each of~them, in accordance with our theoretical~understanding.
\end{iii}

\noindent
On~the
surface~$U$,
strong approximation is~violated. For~example, the Brauer
class~$2\alpha$
already yields that
\mbox{$\bbZ_5$-integral}
points of the kind that
$\smash{\ev_{\alpha,5}(\xi) = \frac14}$
or~$\smash{\frac34}$
cannot be approximated by
\mbox{$\bbZ$-integral}~points.
The~\mbox{$4$-torsion}
class
$\alpha$
enforces additional constraints, for instance that adelic points
$\bbx$
such that
$\smash{\ev_{\alpha,5}(\bbx_5) = \frac12}$,
$\ev_{\alpha,3}(\bbx_3) = 0$,
and~$\ev_{\alpha,7}(\bbx_7) = 0$
may not be approximated,~either.

Finally,~among the
\mbox{$\bbZ$-integral}
points within our range, only
$(1\!:\!-9\!:\!-2\!:\!1\!:\!-86)$
has
evaluation~$\smash{\frac12}$
at~$5$.
In~particular, we see only three of the four admissible combinations of~values.
\end{ex}

\begin{rems}
\begin{iii}
\item
The~calculations described in Example~\ref{fl176} took around 2 hours and 12 minutes, running {\tt magma}, version 2.23.4, on one core of an AMD Phenom II X4 955~processor. In~comparison with the more specific and more advanced class field theoretic methods, which do not apply here, this is, of course, disappointingly~slow.

First~of all, the point search using Elkies' method took more than three~minutes. Then~steps i) to~iii) of Algorithm~\ref{alg_gen} play the role of an initialisation procedure. Their~running time is dominated by the size of the
field~$l$
of definition of the 16~lines, which is a number field of degree~64 for the surface~above. They~take around two~minutes and 20~seconds, the lion's share of which accounts for the Gr\"obner base calculation
over~$l$
in step~iii) of Algorithm~\ref{alg_Pic}.

For~each of the four primes, step~iv) of Algorithm~\ref{alg_gen} is another initialisation~step. It~takes about 90~seconds, each time, which is mainly the time to compute the decomposition group. Then~points are evaluated at a frequency of about one per almost 10~seconds,
for~$p=5$,
and at about one per
$1.7$~seconds,
for the other~primes. A~considerable part of this time is spent on the determination of the
$\#D_w^2$
values of the
\mbox{$2$-cocycle}
$c_{\nu,\xi}$,
which means to multiply 40~elements
of~$l$,
each~time.
\item
There is some {\tt magma} code available from the first author's web page that shows our implementation of Algorithm~\ref{alg_gen}. It includes the concrete calculations for Example~\ref{fl176}, as well as the output that is generated.
\end{iii}
\end{rems}


\subsubsection*{The modification necessary for the case that the residue
characteristic divides the group order}\leavevmode

\noindent
In~our application, this is simply the case when the residue characteristic
is~$2$.
Our~approach to this situation is even more generic and in a certain sense experimental, since Lemma~\ref{Gal_coh_liminv}.a) does not give any indication on how large the integer
$n$
has to be chosen in~practice.
For~$n$
fixed,
we rely on the short exact~sequence
$$0 \longrightarrow {\calO_l}^*\!/(1+\frakm_{\calO_l}^n) \longrightarrow l^*\!/(1+\frakm_{\calO_l}^n) \stackrel{\overline\nu_l}{\longrightarrow} \bbZ \longrightarrow 0$$
and on the canonical isomorphism
${\calO_l}^*\!/(1+\frakm_{\calO_l}^n) \cong ({\calO_l}/\frakm_{\calO_l}^n)^*$.

\begin{alg}[for the case that
$\#G$
is a power
of~$2$,
and the Brauer class is of at most
\mbox{$4$-torsion}]\leavevmode
\label{alg_gen2}

\noindent
Given~a positive
integer~$n$
and a
\mbox{$2$-cocycle}
$\varphi\colon G \times G \to l^*$,
for
$l/k$
a finite Galois extension of local fields and
$G \cong \Gal(l/k)$,
satisfying the conditions listed, this algorithm computes the invariant
in~$\bbQ/\bbZ$
of the class in
$H^2(G, l^*)$
represented
by~$\varphi$,
as soon as
$n$
is sufficiently~large.

\begin{iii}
\item
Run step~i) of Algorithm~\ref{eval_loc}.
If~$f=1$
or~$f=2$
then run a totally naive modification, which just replaces
$l$
by~$lk^{\nr,4}$.
Then~$f$
is a multiple
of~$4$.
\item
Fix~a uniformiser
$\pi_l$
of~$l$
and calculate the unit group of the finite ring
$\calO_l/\frakm_{\calO_l}^n$,
together with a partial map
$\pr\colon \calO_l \dashrightarrow \bbZ/a_1\bbZ \oplus \cdots \oplus \bbZ/a_m\bbZ$.

Then~create the
\mbox{$G$-module}
$l^*\!/(1+\frakm_{\calO_l}^n) $
as~follows.
Put~$M := \bbZ \oplus \bbZ/a_1\bbZ \oplus \cdots \oplus \bbZ/a_m\bbZ$.
For~each generator
$g_i \in G$,
construct an
\mbox{$(m+1)\times(m+1)$-matrix}
$A_i$
in the following~manner.

Put~the first row to be
$(1,0,\ldots,0)$.
Complete~the first column by the coefficients
of~$\smash{\pr(\frac{g_i(\pi_l)}{\pi_l})}$.
For~$j = 1,\ldots,n$,
complete~the
\mbox{$(j+1)$-st}
column
by~$\pr(g_i(\widetilde{e}_j))$,
for
$\widetilde{e}_j$
a lift of the
\mbox{$j$-th}
generator
$e_j$
of the unit group
$\bbZ/a_1\bbZ \oplus \cdots \oplus \bbZ/a_m\bbZ$
to~$\calO_l$.

Finally,~transform
$M$~into
a
\mbox{$G$-module}
using the
matrices~$A_i$.
\item
Define the standard
\mbox{$2$-cocycle}
$\st\colon G \times G \to M$~by
$$(\sigma,\tau) \mapsto \left\{
\begin{array}{cl}
(0,0,\ldots,0)   & \text{ if } \expo(\sigma) + \expo(\tau) < f \,, \\
(e, \pr(\frac{\pi}{\pi_l^e})) & \text{ if } \expo(\sigma) + \expo(\tau) \geq f \,,
\end{array}
\right.
$$
for~$\pi$
a uniformiser
of~$k$.
\item
Calculate the image of the input cocycle by~putting
$$\textstyle \widetilde\varphi(\sigma, \tau) := \big( \nu_l(\varphi(\sigma, \tau)), \pr (\frac{\varphi(\sigma, \tau)}{\pi_l^{\nu_l(\varphi(\sigma, \tau))}})\big) ,$$
for every
$(\sigma, \tau) \in G \times G$.
\item
Calculate the abelian group
$H^2(G, M)$
and, as its elements, the classes
$c_\st$
and~$c_{\widetilde\varphi}$
of the
\mbox{$2$-cocycles}
$\st$
and~$\widetilde\varphi$.

If~$2c_\st = 0$
then output a message saying that a larger value
of~$n$
has to be chosen and terminate~immediately. Otherwise,~find a natural number
$m$
such that
$c_{\widetilde\varphi} = m c_\st$,
and output
$(m/f \bmod 1) \in \bbQ/\bbZ$
as the desired~invariant.
\end{iii}
\end{alg}

\begin{rem}
We~successfully ran Algorithm~\ref{alg_gen2} for several~examples. In~these, the minimal value
of~$n$,
for which the calculations went through, varied between
$3$
and~$9$.
\end{rem}

\setlength\parindent{0mm}
\end{document}